\documentclass{amsart}%
\usepackage{float}
\usepackage{amsmath}
\usepackage{graphicx}
\usepackage{latexsym}
\usepackage{amsfonts}
\usepackage{amssymb}%
\theoremstyle{plain}
\newtheorem{theorem}{Theorem}

\newtheorem{lemma}[theorem]{Lemma}
\newtheorem{proposition}[theorem]{Proposition}
\theoremstyle{definition}

\newtheorem{remark}[theorem]{Remark}
\newtheorem*{remark*}{Remark}

\begin{document}
\title[Fluctuations under a stable catalyst]{Hydrodynamic limit fluctuations\vspace{4pt}\\of super-Brownian motion\vspace{4pt}\\with a stable catalyst\vspace{40pt}}
\author[Fleischmann]{Klaus Fleischmann}
\address{Weierstrass Institute for Applied Analysis and Stochastics, Mohrenstr.\ 39,
D--10117 Berlin, Germany}
\email{fleischm@wias-berlin.de}
\urladdr{http://www.wias-berlin.de/\symbol{126}fleischm}
\author[M\"{o}rters]{Peter M\"{o}rters}
\address{University of Bath, Department of Mathematical Sciences, Claverton Down, Bath
BA2 7AY, United Kingdom.}
\email{maspm@bath.ac.uk}
\urladdr{http://www.bath.ac.uk/\symbol{126}maspm}
\author[Wachtel]{Vitali Wachtel}
\address{Weierstrass Institute for Applied Analysis and Stochastics, Mohrenstr.\ 39,
D--10117 Berlin, Germany}
\email{vakhtel@wias-berlin.de}
\urladdr{http://www.wias-berlin.de/\symbol{126}vakhtel}
\thanks{K.F. and V.W.: Supported by the DFG}
\thanks{P.M.: Supported by EPSRC grant EP/C500229/1 and an Advanced Research Fellowship}
\thanks{WIAS preprint No.\ 1052 of August 17, 2005,\quad ISSN 0946\thinspace
--\thinspace8633,\quad bounds13.tex}
\thanks{Corresponding author: Klaus Fleischmann}
\maketitle

\thispagestyle{empty}
\setcounter{page}{0}\newpage

\begin{quotation}
\noindent\textsc{Abstract. }{\small We consider the behaviour of a continuous
super-Brownian motion catalysed by a random medium with infinite overall
density under the hydrodynamic scaling of mass, time, and space. We show that,
in supercritical dimensions, the scaled process converges to a macroscopic
heat flow, and the appropriately rescaled random fluctuations around this
macroscopic flow are asymptotically bounded, in the sense of log-Laplace
transforms, by generalised stable Ornstein-Uhlenbeck processes. The most
interesting new effect we observe is the occurrence of an index-jump from a
`Gaussian' situation to stable fluctuations of index }$1+\gamma${\small ,
where }$\gamma\in(0,1)$ {\small is an index associated to the medium. }
\end{quotation}

{\small
\tableofcontents
}

\section{Introduction and main results\label{S.intr}}

\subsection{Motivation and background\label{SS.mot}}

In order to describe the long-term behaviour of infinite interacting spatial
particle systems with mass preservation on average, limit theorems under
mass-time-space rescaling are an established tool. A typical feature that can
be captured by this means is the clumping behaviour of spatial branching
processes in \emph{low dimensions}: In some models, for a critical scaling one
can observe convergence to a nontrivial field of isolated mass clumps. The
spatial contraction allows to get hold of large mass clumps in remote
locations, and the index of mass-rescaling serves as a measure of the strength
of the clumping effect, quantifying the degree of intermittency. In some of
these results a macroscopic time dependence can be retained, giving insight in
the long-time developments of the clumps. For a recent result in this
direction, see Dawson et al.\ \cite{DawsonFleischmannMoerters2002.sclump}.

In \emph{higher dimensions} one does not expect to observe clumping under
mass-time-space rescaling, but convergence to a non-random mass flow, the
\emph{hydrodynamic limit}. In this case one can hope to get a deeper
understanding from the investigation of fluctuations around this limit. Such
fluctuations were studied by Holley and Stroock \cite{HolleyStroock1978} and
Dawson \cite{Dawson1978.k}, and their results were later refined and extended,
e.g.\ by Dittrich \cite{Dittrich1987}. There is also a large body of
literature on hydrodynamic limits of interacting particle systems, see
e.g.\ \cite{MasiPresutti1991,KipnisLandim1999,Spohn1991}. Our main motivation
behind this paper is to investigate the possible effects on fluctuations
around the hydrodynamic limit if the original process is influenced by a
\emph{random medium}, which in our model acts as a catalyst for the local
branching rates.

In Dawson et al.\ \cite{DawsonFleischmannGorostiza1989.hydro}, fluctuations
under mass-time-space rescaling were derived for a class of spatial infinite
branching particle systems in $\,\mathbb{R}^{d}$\thinspace\ (with symmetric
$\alpha$--stable motion and $(1+\beta)$--branching) in supercritical
dimensions in a random medium with \emph{finite}\/ overall density. This leads
to \emph{generalized Ornstein-Uhlenbeck processes}\/ which are the same as for
the model in the constant (averaged) medium. In other words, for the
log-Laplace equation the governing effect is homogenization: After rescaling,
the equation approximates an equation with homogeneous branching rate, the
medium is simply averaged out. The nature of the fluctuations for the case of
a medium with \emph{infinite}\/ overall density remained unresolved over the
years.\vspace{2pt}

The purpose of the present paper is to get progress in this direction. Our
main result shows that a medium with an \emph{infinite}\/ overall density can
have a drastic effect on the fluctuation behaviour of the model under critical
rescaling in supercritical dimensions, and homogenization is no longer the
effect governing the macroscopic behaviour. In fact, despite the infinite
overall density of the medium, we still have a law of large numbers under a
certain mass-time-space rescaling. But under this scaling, the variances
(given the medium) blow up, and the related fluctuations do \emph{not}\/ obey
a central limit theorem. However, fluctuations can be described to some degree
by a stable process.

To be more precise, we start with a branching system with \emph{finite}\/
variance given the medium, considered as a branching process with a random
law, where this randomness of the laws comes from the randomness of the medium
(quenched approach). Under a mass-time-space rescaling, the random laws of the
fluctuations are asymptotically bounded from above and below by the laws of
constant multiples of a generalized Ornstein-Uhlenbeck process with
\emph{infinite}\/ variance. Here the ordering of random laws is defined in
terms of the random Laplace transforms. The generalized Ornstein-Uhlenbeck
process is the same as the fluctuation limit of a super-Brownian motion with
infinite variance branching in the case of a constant medium. In fact, the
branching mechanism is $(1+\gamma)$--branching, where $\,\gamma\in
(0,1)$\thinspace\ is the index of the medium. Altogether, the present result
is a big step towards an affirmative answer to the old open problem of
understanding fluctuations in the case of a random medium with infinite
overall density. It also leads to random medium effects which are in line with
experiences concerning the clumping behaviour in subcritical dimensions as in
\cite{DawsonFleischmannMoerters2002.sclump}.

\subsection{Preliminaries: notation}

For $\lambda\in\mathbb{R}$, introduce the reference function%
\begin{equation}
\phi_{\lambda}(x)\;:=\;\mathrm{e}^{-\lambda|x|}\quad\text{for}\ \,x\in
\mathbb{R}^{d}. \label{ref.fct}%
\end{equation}
For $\,f:\mathbb{R}^{d}\rightarrow\mathbb{R},$ set%
\begin{equation}
|f|_{\lambda}\;:=\;\Vert f/\phi_{\lambda}\Vert_{\infty} \label{norm.l}%
\end{equation}
where $\,\Vert\cdot\Vert_{\infty}$\thinspace\ refers to the supremum norm.
Denote by $\,\mathcal{C}_{\lambda}$\thinspace\ the separable Banach space of
all continuous functions $\,f:\mathbb{R}^{d}\rightarrow\mathbb{R}$%
\thinspace\ such that $\,|f|_{\lambda}$\thinspace\ is finite and that
$\,f(x)/\phi_{\lambda}(x)$\thinspace\ has a finite limit as $\,|x|\rightarrow
\infty.$\thinspace\ Introduce the space
\begin{equation}
\mathcal{C}_{\exp}\;=\;\mathcal{C}_{\exp}(\mathbb{R}^{d})\;:=\;\bigcup
_{\lambda>0}\,\mathcal{C}_{\lambda\,} \label{C.spaces}%
\end{equation}
of (at least) \emph{exponentially decreasing}\/ continuous test functions on
$\mathbb{R}^{d}$. An index $+$ as in $\mathbb{R}_{+}$ or $\mathcal{C}_{\exp
}^{+}$ refers to the corresponding non-negative members.

Let $\,\mathcal{M}=\mathcal{M}(\mathbb{R}^{d})$\thinspace\ denote the set of
all (non-negative) Radon measures $\,\mu$\thinspace\ on $\,\mathbb{R}^{d}%
$\thinspace\ and $\,\mathrm{d}_{0}$\thinspace\ a complete metric on
$\,\mathcal{M}$\thinspace\ which induces the vague topology. Introduce the
space $\,\mathcal{M}_{\mathrm{tem}}=\mathcal{M}_{\mathrm{tem}}(\mathbb{R}%
^{d})$\thinspace\ of all measures $\,\mu$\thinspace\ in $\,\mathcal{M}%
$\thinspace\ such that $\,\left\langle \mu,\phi_{\lambda}\right\rangle
:=\int\!{d}\mu\;\phi_{\lambda}<\infty,$\thinspace\ for all $\,\lambda
>0.$\thinspace\ We topologize this set $\,\mathcal{M}_{\mathrm{tem}}%
$\thinspace\ of \emph{tempered}\/ measures by the metric
\begin{equation}
\mathrm{d}_{\mathrm{tem}}(\mu,\nu)\;:=\;\mathrm{d}_{0}(\mu,\nu)+\sum
_{n=1}^{\infty}2^{-n}\left(  |\mu-\nu|_{1/n}\,_{\!_{\!_{\,}}}\wedge\,1\right)
\!\quad\text{for}\,\ \mu,\nu\in\mathcal{M}_{\mathrm{tem\,}}.
\end{equation}
Here $\,|\mu-\nu|_{\lambda}$\thinspace\ is an abbreviation for $\,\big
|\langle\mu,\phi_{\lambda}\rangle-\langle\nu,\phi_{\lambda}\rangle
\big|.$\thinspace\ Note that $\,\mathcal{M}_{\mathrm{tem}}$\thinspace\ is a
Polish space (that is, $\left(  \mathcal{M}_{\mathrm{tem\,}},\mathrm{d}%
_{\mathrm{tem}}\right)  $ is a complete separable metric space), and that
$\,\mu_{n}\rightarrow\mu$\thinspace\ in $\,\mathcal{M}_{\mathrm{tem}}%
$\thinspace\ if and only if%
\begin{equation}
\left\langle \mu_{n\,},\varphi\right\rangle \;\underset{n\uparrow\infty
}{\longrightarrow}\;\left\langle \mu,\varphi\right\rangle \quad\text{for
}\,\varphi\in\mathcal{C}_{\exp\,}.
\end{equation}

Probability measures will be denoted as $\,\mathsf{P},\mathbb{P},\mathcal{P}%
,$\ whereas $\,\mathsf{E},\mathbb{E},\mathcal{E}$\thinspace\ and
$\,\mathsf{V}\!$ar$,\mathbb{V}\!$ar$,\mathcal{V}\!$ar\thinspace\ refer to the
corresponding expectation and variance symbols.

Let $p$ denote the standard heat kernel in $\mathbb{R}^{d}$ given by
\begin{equation}
p_{t}(x)\;:=\;(2\pi t)^{-d/2}\,\exp\!\Big[-\frac{|x|^{2}}{2t}\,\Big ]\quad
\text{for}\,\ t>0,\,\ x\in\mathbb{R}^{d}. \label{notation.heat.kernel}%
\end{equation}
Write $W=\bigl(W,\,(\mathcal{F}_{t})_{t\geq0\,},\,\mathcal{P}_{x}%
,\,x\in\mathbb{R}^{d}\bigr)$ for the corresponding (standard) Brownian motion
in $\mathbb{R}^{d}$ with natural filtration, and $S=\left\{  S_{t}%
:\,t\geq0\right\}  $ for the related semigroup. Quantities depending on time
$t,$ as $p_{t},S_{t}$ or solutions $u(t,$\thinspace$\cdot\,),$ are formally
set to $0$ if $t<0.$

Let $\,\ell$\thinspace\ denote the Lebesgue measure on $\,\mathbb{R}^{d}%
.$\thinspace\ Write $\,B(x,r)$\thinspace\ for the closed ball\thinspace
\ around $\,x\in\mathbb{R}^{d}$\thinspace\ with radius $\,r>0.$\thinspace\ In
this paper, $G$ denotes the Gamma function.

With $c=c(q)$\thinspace\ we always denote a positive constant which (in the
present case) might depend on a quantity $q$ and might also change from place
to place. Moreover, an index on $c$ as $c_{(\mathrm{\#})}$ or $c_{\mathrm{\#}%
}$ will indicate that this constant first occurred in formula line (\#) or
(for instance) Lemma~\#, respectively. We apply the same labelling rules also
to parameters like $\,\lambda$\thinspace\ and $\,k.$

\subsection{Modelling of catalyst and reactant\label{SS.mod}}

Of course, there is some freedom in choosing the model we want to work with.
To avoid unnecessary limit procedures, we work on $\,\mathbb{R}^{d}$%
\thinspace\ and with continuous-state branching as the branching system,
namely with continuous super-Brownian motion, which is a spatial version of
Feller's branching diffusion. The branching rate of an intrinsic
`particle'\ varies in space and in fact is selected from a random field to be
specified. In this context, it is convenient to speak also of the random field
as the \emph{catalyst, }and of the branching system given the random medium as
the \emph{reactant.}

First we want to specify the catalyst. In our context, a very natural way is
to start from a \emph{stable random measure}\/ $\Gamma$ on $\mathbb{R}^{d}$
with index $\,\gamma\in(0,1)$ determined by its log-Laplace functional%
\begin{equation}
-\log\mathsf{E}\exp\left\langle \Gamma,-\varphi\right\rangle \;=\;\int
\!\!{d}z\;\varphi^{\gamma}(z)\quad\text{for}\,\ \varphi\in\mathcal{C}%
_{\mathrm{\exp}\,}^{+}. \label{logL.Gamma}%
\end{equation}
(The letter $\,\mathsf{P}$\thinspace\ always stands for the law of the
catalyst, whereas $\,\mathbb{P}$\thinspace\ is reserved for the law of the
reactant given the catalyst.) See, for instance, \cite[Lemma 4.8]%
{DawsonFleischmann1992.equ} for background concerning $\,\Gamma.$%
\thinspace\ Clearly, $\,\Gamma$\thinspace\ is a spatially homogeneous random
measure with independent increments and infinite expectation. $\,\Gamma
$\thinspace\ has a simple \emph{scaling property},%
\begin{equation}
\Gamma(k\,dz)\ \overset{\mathcal{L}}{=}\ k^{d/\gamma}\,\Gamma(dz)\quad
\text{for}\,\ k>0, \label{scal.Gamma}%
\end{equation}
where $\,\overset{\mathcal{L}}{=}$\thinspace\ refers to equality in law.
However, $\,\Gamma$\thinspace\ is a purely atomic measure, hence, its atoms
cannot be hit by a Brownian path or a super-Brownian motion in dimensions
$\,d\geq2.$\thinspace\ Thus, $\,\Gamma$\thinspace\ cannot serve directly as a
catalyst for a non-degenerate reaction model based on Brownian particles in
higher dimensions. Therefore we look at the density function after smearing
out $\,\Gamma$\thinspace\ by the (non-normalized) function $\,\vartheta
_{1\,},$\thinspace\ where $\,\vartheta_{r}:=\mathsf{1}_{B(0,r)\,},$%
\thinspace\ $r>0,$\thinspace\ that is,%
\begin{equation}
\Gamma^{1}(x)\;:=\;\int\!\Gamma({d}z)\;\vartheta_{1}(x-z)\quad\text{for}%
\,\ x\in\mathbb{R}^{d}. \label{not.Gamma.eps}%
\end{equation}
In the sequel, the unbounded function $\,\Gamma^{1}$\thinspace\ with infinite
overall density will play the r\^{o}le of the \emph{random medium}: It will
act as a catalyst that determines the spatially varying branching rate of the
reactant. Once again, smoothing\ is needed, since otherwise the medium will
not be hit by an intrinsic Brownian reactant particle. In our proofs, the
independence and scaling properties of $\,\Gamma$\thinspace\ will be
advantageous, though one would expect analogous results to hold for quite
general random media with infinite overall density.

Consider now the \emph{continuous super-Brownian motion}\/ $\,X=X[\Gamma^{1}%
]$\thinspace\ in $\,\mathbb{R}^{d},$\thinspace\ $d\geq1,$\thinspace
\ \emph{with random catalyst}\/ \thinspace$\Gamma^{1}.$\thinspace\ More
precisely, $\,$\emph{for almost all samples}\/ $\,\Gamma^{1},\,\ $this is a
continuous time-homogeneous Markov process $\,X=X[\Gamma^{1}]=(X,\,\mathbb{P}%
_{\mu},\,\mu\in\mathcal{M}_{\mathrm{tem}})$\thinspace\ with log-Laplace
transition functional
\begin{equation}
-\log\mathbb{E}_{\mu}\!\exp\left\langle X_{t\,},-\varphi\right\rangle
\;=\;\left\langle \mu,u(t,\,\cdot\,)_{\!_{\!_{{}}}}\right\rangle
\quad\text{for}\,\ \varphi\in\mathcal{C}_{\exp\,}^{+},\;\,\mu\in
\mathcal{M}_{\mathrm{tem\,}}, \label{loglap}%
\end{equation}
where $\,u=u[\varphi,\Gamma^{1}]=\left\{  u(t,x):\,t\geq0,\;x\in\mathbb{R}%
^{d}\right\}  $ is the unique mild non-negative solution of the reaction
diffusion equation%
\begin{equation}
\frac{\partial}{\partial t}u(t,x)\;=\;\tfrac{1}{2}\Delta u(t,x)-\varrho
\,\Gamma^{1}(x)\,u^{2}(t,x)\quad\text{for}\,\ t\geq0,\;\,x\in\mathbb{R}^{d},
\label{logLap.equ}%
\end{equation}
with initial condition $\,u(0,\,\cdot\,)=\,\varphi.\,\ $Here\thinspace
\ $\varrho>0$\thinspace\ is an additional parameter (for scaling purposes).
For background on super-Brownian motion we recommend \cite{Dawson1993},
\cite{Etheridge2000}, or \cite{Perkins2002}, and for a survey on catalytic
super-Brownian motion, see e.g.\ \cite{DawsonFleischmann2002.survasc} or
\cite{Klenke2000}.

\textrm{F}rom Dawson and Fleischmann
\cite{DawsonFleischmann1983.env1,DawsonFleischmann1985.env2} the following
dichotomy concerning the long-term behaviour of $\,X$\thinspace\ is basically
known (although there the phase space is $\,\mathbb{Z}^{d}$\thinspace\ and the
processes are in discrete time): Starting from the Lebesgue measure
$\,X_{0}=\ell$, $\,$the process $\,X$\thinspace\ dies locally in law as
$\,t\uparrow\infty$\thinspace\ if $\,d\leq2/\gamma\,\ $(recall that
$\,0<\gamma<1\,\ $is the index of the random medium $\,\Gamma^{1}),$%
\thinspace\ whereas in all higher dimensions one has persistent convergence in
law to a non-trivial limit state denoted by $\,X_{\infty\,}$. \emph{From now
on, we restrict our attention to (supercritical) dimensions} $\,d>2/\gamma.$

We are interested in the large scale behaviour of $\,X.$

\subsection{Main results of the paper}

Introduce the \emph{scaled processes}\/ $\,X^{k},$\thinspace\ $k>0,$%
\thinspace\ defined by%
\begin{equation}
X_{t}^{k}(B)\;:=\;k^{-d}\,X_{k^{2}t}(kB)\quad\text{for}\,\ t\geq
0,\;\,B\subseteq\mathbb{R}^{d}\text{ \thinspace Borel.} \label{scale}%
\end{equation}
This \emph{hydrodynamic rescaling} leaves the underlying Brownian motions
invariant (in law), and the expectation of the scaled process is the heat
flow:%
\begin{equation}
\mathbb{E}_{\mu}X_{t}^{k}\;=\;S_{t}\,\mu^{k}\quad\text{for}\,\ X_{0}=\mu
\in\mathcal{M}_{\mathrm{tem\,}}. \label{expect}%
\end{equation}
In particular, if $\,X$\thinspace\ is started with the Lebesgue measure
$\,\ell,$\thinspace\ the expectation is preserved in time. We also define the
\emph{critical scaling index}%
\begin{equation}
\varkappa_{\mathrm{c}}\;:=\;\frac{\gamma d-2}{1+\gamma}\;>\;0. \label{not.ski}%
\end{equation}

\begin{theorem}
[\textbf{Refined law of large numbers}]\label{C.LLN}Suppose $\,d>2/\gamma
.$\thinspace\ Start $\,X$\thinspace\ with $k$--dependent initial states
$\,X_{0}=\mu_{k}\in\mathcal{M}_{\mathrm{tem}}$\thinspace\ such that $X_{0}%
^{k}=\mu\in\mathcal{M}_{\mathrm{tem}}$ \thinspace for $\,k>0.$\thinspace\ If
$\,\varkappa\in\lbrack0,\varkappa_{\mathrm{c}}),$\thinspace\ then
\begin{equation}
k^{\varkappa}\!\left(  X_{t}^{k}-S_{t}\mu\right)  \;\underset{k\uparrow\infty
}{\Longrightarrow}\;0\quad\text{in}\,\ \mathsf{E}\mathbb{P}_{\mu_{k}%
}\text{--law.} \label{LLN}%
\end{equation}

\end{theorem}

The refined law of large numbers is actually a by-product of the proofs of our
main result, as will be explained immediately after
Proposition~\ref{L.exp.conv}.

In contrast to \cite{DawsonFleischmannGorostiza1989.hydro}, in the present
paper we use Laplace transforms instead of Fourier transforms although
fluctuations we are interested in are signed objects. This is possible since
these fluctuations themselves are deviations from non-negative $\,X^{k}%
,$\thinspace\ and related stable limiting quantities have skewness parameter
$\,\beta=-1$,\thinspace\ for which Laplace transforms are meaningful.

For $\,x\in\mathbb{R}^{d}$\thinspace\ we put%
\begin{equation}
\mathrm{en}(x)\ :=\ \left\{
\begin{array}
[c]{ll}%
\log^{+}\!\left(  |x|^{-1}\right)  & \text{if}\,\ d=4,\vspace{6pt}\\
\,|x|^{4-d} & \text{if}\,\ d\geq5,
\end{array}
\right.  \label{en}%
\end{equation}
and for $\,\mu\in\mathcal{M}_{\mathrm{tem}\,},$ and $\lambda>0$,
\begin{equation}
\mathrm{En}_{\lambda}(\mu)\ :=\ \int\!\!\mu(dx)\,\phi_{\lambda}(x)\int
\!\!\mu(dy)\,\phi_{\lambda}(y)\ \mathrm{en}(x-y). \label{En}%
\end{equation}
Note that $\,\mathrm{En}_{\lambda}(\delta_{x})\equiv\infty$\thinspace\ if
$\,d>3.$

\begin{theorem}
[\textbf{Asymptotic fluctuations}]\label{T.fluct}Suppose $\,d>2/\gamma
.$\thinspace\ Start $\,X$\thinspace\ with $k$--dependent initial states
$\,X_{0}=\mu_{k}\in\mathcal{M}_{\mathrm{tem}}$\thinspace\ such that $X_{0}%
^{k}=\mu\in\mathcal{M}_{\mathrm{tem}}$ \thinspace for $\,k>0.$\thinspace\ In
the case\/ $\,d>3,$\thinspace\ suppose additionally that $\,\mu$\thinspace\ is
a measure of finite energy in the sense that $\,\mathrm{En}_{\lambda}%
(\mu)<\infty$ for all $\lambda>0$. If $\,\varkappa=\varkappa_{\mathrm{c}\,}%
,$\thinspace\ then there exists constants $\,\overline{c}>\underline{c}%
>0$\thinspace\ such that for any\/ $\,\varphi_{1},\ldots,\varphi_{n}%
\in\mathcal{C}_{\mathrm{exp}}^{+}$\thinspace\ and\/ $\,0=:t_{0}\leq t_{1}%
\leq\dots\leq t_{n\,},$\thinspace\ in $\mathsf{P}$--probability,%
\begin{equation}%
\begin{array}
[c]{l}%
\displaystyle
\limsup_{k\rightarrow\infty}\mathbb{E}_{\mu_{k}}\!\exp\!\Big[\sum_{i=1}%
^{n}k^{\varkappa}\big\langle X_{t_{i}}^{k}-S_{t_{i}}\mu,\,-\varphi
_{i}\big\rangle\Big]\vspace{6pt}\\%
\displaystyle
\qquad\leq\ \exp\!\Bigg[\overline{c}\,\bigg\langle\mu,\sum_{i=1}^{n}%
\int_{t_{i-1}}^{t_{i}}\!dr\ S_{r}\Big(\Big(\sum_{j=i}^{n}S_{t_{j}-r}%
\varphi_{j}\Big)^{1+\gamma}\Big)\bigg\rangle\Bigg]
\end{array}
\label{fluct1}%
\end{equation}
and%
\begin{equation}%
\begin{array}
[c]{l}%
\displaystyle
\liminf_{k\rightarrow\infty}\mathbb{E}_{\mu_{k}}\!\exp\!\Big[\sum_{i=1}%
^{n}k^{\varkappa}\big\langle X_{t_{i}}^{k}-S_{t_{i}}\mu,\,-\varphi
_{i}\big\rangle\Big]\vspace{6pt}\\%
\displaystyle
\qquad\geq\ \exp\!\Bigg[\underline{c}\,\bigg\langle\mu,\sum_{i=1}^{n}%
\int_{t_{i-1}}^{t_{i}}\!dr\ S_{r}\Big(\Big(\sum_{j=i}^{n}S_{t_{j}-r}%
\varphi_{j}\Big)^{1+\gamma}\Big)\bigg\rangle\Bigg].
\end{array}
\label{fluct1'}%
\end{equation}

\end{theorem}

Explicit values of $\,\overline{c}$\thinspace\ and $\,\underline{c}$%
\thinspace\ are given in (\ref{c.fluct}) and (\ref{c.u}), respectively.

\begin{remark}
[\textbf{Generalized Ornstein-Uhlenbeck process}]\label{R.limit.process}%
\hspace{-1.2pt}The right-hand sides of (\ref{fluct1}) and (\ref{fluct1'}) are
the Laplace transforms of the finite-dimensional distributions of different
multiples of a process $Y$ taking values in the Schwartz space of tempered
distributions. This process $Y$ can be called a \emph{generalized
Ornstein-Uhlenbeck process} as it solves the \emph{generalized Langevin
equation},
\begin{equation}
dY_{t}\ =\ \tfrac{1}{2}\Delta Y_{t}\,dt+dZ_{t}\quad\text{for}\,\ t\geq
0,\,\ Y_{0}=0,
\end{equation}
where $dZ_{t}/dt$\thinspace\ is a $(1+\gamma)$--stable noise, i.e.\ $Z$%
\thinspace\ is the process with independent increments with values in the
Schwartz space such that, for $0\leq s\leq t\,\ $and $\,\varphi\in
\mathcal{C}_{\mathrm{exp}}^{+}$,
\begin{equation}
E\mathrm{e}^{-\left\langle Z_{t}-Z_{s},\varphi\right\rangle }\ =\ \exp
\!\Big[\int_{s}^{t}\!dr\,\left\langle S_{r}\mu,\varphi^{1+\gamma}\right\rangle
\Big ].
\end{equation}
$Y$ is described in detail in \cite[Section 4]%
{DawsonFleischmannGorostiza1989.hydro}, where it appeared as the hydrodynamic
fluctuation limit process corresponding to super-Brownian motion with finite
mean branching rate, but with infinite variance $(1+\gamma)$--branching.
Recall that the Markov process $Y$ has log-Laplace transition functional
\begin{equation}
-\log E\left\{  \exp\left\langle Y_{t\,},-\varphi\right\rangle \;\big|\;Y_{0}%
\right\}  \;=\;\left\langle Y_{0},S_{t}\varphi\right\rangle +\left\langle
\mu,v(t,\,\cdot\,)_{\!_{\!_{{}}}}\right\rangle \quad\text{for}\,\ \varphi
\in\mathcal{C}_{\exp\,}^{+},
\end{equation}
where $\,v=v[\varphi]=\left\{  v(t,x):\,t\geq0,\;x\in\mathbb{R}^{d}\right\}  $
solves$\ $%
\begin{equation}%
\begin{array}
[c]{c}%
\displaystyle
\frac{\partial}{\partial t}v(t,x)\;=\;\tfrac{1}{2}\Delta v(t,x)\,+\,(S_{t}%
\varphi)^{1+\gamma}\,(x)\vspace{6pt}\\
\text{with initial condition }\,v(0,\,\cdot\,)=0.
\end{array}
\label{equ.lim}%
\end{equation}

In particular, in our limit procedure the finite variance property of the
original process given the medium is lost and, by a subtle averaging effect,
an index jump of size $1-\gamma>0$ occurs. \hfill$\Diamond$
\end{remark}

\begin{remark}
[\textbf{Ordering}]\hspace*{5pt}The stochastic ordering of the random laws in
our \mbox{asymptotic} bounds in (\ref{fluct1}) and (\ref{fluct1'}) is
well-known in queueing and risk theory, see \cite{MullerStoyan2002} for
background.\hfill$\Diamond$
\end{remark}

\begin{remark}
[\textbf{Existence of a fluctuation limit}]Theorem~\ref{T.fluct} leaves
\emph{open}, whether a fluctuation limit exists in $\mathsf{P}$--probability
and whether it is a generalised Ornstein-Uhlenbeck process as described
above.\hfill$\Diamond$
\end{remark}

\begin{remark}
[\textbf{Variance considerations}]\label{R.variance}In the case $\,\mu
_{k}\equiv\ell,$\thinspace\ for $\,\varphi\in\mathcal{C}_{\exp\,},$%
\thinspace\ the $\mathsf{P}$--random variance
\begin{align}
\mathbb{V}\mathrm{ar}_{\ell}  &  \left[  k^{\varkappa}\langle X_{t}^{k}%
-S_{t}\mu,\varphi\rangle\right]  \;=\;k^{2\varkappa}\,\mathbb{V}%
\!\mathrm{ar}_{\ell}\langle X_{t}^{k},\varphi\rangle\label{variance}\\
&  =\;2\varrho\,k^{2\varkappa-2d}\int_{0}^{k^{2}t}\!{d}s\int\!\!{d}%
x\;\Gamma^{1}(x)\,\big[S_{k^{2}t-s}\varphi\,(k^{-1}\cdot\,)\big]^{2}%
(x)\nonumber
\end{align}
equals (by scaling) approximately%
\begin{equation}
2\varrho\,k^{2\varkappa-2d+2+d/\gamma}\int_{0}^{t}\!{d}s\int\!\Gamma
({d}z)\,[S_{s}\varphi]^{2}(z)\quad\text{as}\,\ k\uparrow\infty.
\end{equation}
Hence, for $\,\varkappa$\thinspace\ satisfying
\begin{equation}
0\;\leq\;\varkappa\;<\;\varkappa_{\mathrm{var}}\;:=\;\frac{(2\gamma
-1)d-2\gamma}{2\gamma}\,, \label{kappa.var}%
\end{equation}
implying $\,\gamma\in(\tfrac{1}{2},1)$\thinspace\ and $\,d>2\gamma
/(2\gamma-1),$\thinspace\ the random variances (\ref{variance}) converge to
zero as $\,k\uparrow\infty,$\thinspace\ yielding the refined law of large
numbers (\ref{LLN}), whereas for $\,\varkappa>\varkappa_{\mathrm{var}}%
$\thinspace\ these variances explode. Note that $\,\varkappa_{\mathrm{var}%
}<\varkappa_{\mathrm{c\,}},$\thinspace\ since $\,(\gamma-1)(d-2\gamma
)<0$.\thinspace\ Therefore a quenched variance consideration as in
(\ref{variance}) can only imply statement (\ref{LLN}) in the restricted case
(\ref{kappa.var}). Of course, \emph{annealed} variances are infinite already
for fixed $\,k,$\thinspace\ which follows from (\ref{variance}).\hfill
$\Diamond$
\end{remark}

\subsection{Heuristics, concept of proof, and outline\label{SS.heuristics}}

For this discussion we first focus on the case $n=1$ in Theorem~\ref{T.fluct}.
{F}rom (\ref{loglap}), (\ref{logLap.equ}), and scaling,%
\begin{equation}
\log\mathbb{E}_{\mu_{k}}\!\exp\!\left[  -k^{\varkappa}\langle X_{t}^{k}%
-S_{t}\mu,\varphi\rangle\right]  \ =\ \left\langle \mu,k^{\varkappa}%
S_{t}\varphi-u_{k}(t,\,\cdot\,)_{\!_{\!_{\,}}}\right\rangle \!,
\label{scaling}%
\end{equation}
where $u_{k}$\thinspace\ solves the (scaled) equation%
\begin{equation}%
\begin{array}
[c]{c}%
\displaystyle
\frac{\partial}{\partial t}u_{k}(t,x)\;=\;\tfrac{1}{2}\Delta u_{k}%
(t,x)\,-\,k^{2-d}\,\varrho\,\Gamma^{1}(kx)\,u_{k}^{2}(t,x)\vspace{6pt}\\
\text{with initial condition }\,u_{k}(0,\,\cdot\,)=k^{\varkappa}\varphi.
\end{array}
\label{logLap.equ.k}%
\end{equation}

Since $\,v(t,x):=k^{\varkappa}S_{t}\varphi\,(x)$\thinspace\ is the solution of%
\begin{equation}
\frac{\partial}{\partial t}v(t,x)\;=\;\tfrac{1}{2}\Delta v(t,x)\quad\text{with
initial condition}\,\ v(0,\,\cdot\,)=k^{\varkappa}\varphi,
\end{equation}
we see that $\,f_{k}(t,x):=k^{\varkappa}S_{t}\varphi\,(x)-u_{k}(t,x)$%
\thinspace\ solves%
\begin{equation}%
\begin{array}
[c]{c}%
\displaystyle
\frac{\partial}{\partial t}f_{k}(t,x)\;=\;\tfrac{1}{2}\Delta f_{k}%
(t,x)\,+\,k^{2-d}\,\varrho\,\Gamma^{1}(kx)\left[  k^{\varkappa}S_{t}%
\varphi\,(x)-f_{k}(t,x)_{\!_{\!_{\,}}}\right]  ^{2}\vspace{6pt}\\
\text{with initial condition }\,f_{k}(0,\,\cdot\,)=0.
\end{array}
\label{equ.f}%
\end{equation}

Consider now the critical scaling $\varkappa=\varkappa_{\mathrm{c\,}}%
.$\thinspace\ By our claims in Theorem~\ref{T.fluct}, $\,f_{k}$\thinspace
\ should be asymptotically bounded in $\mathsf{P}$--law by solutions
$\,v$\thinspace\ of$\ $%
\begin{equation}%
\displaystyle
\frac{\partial}{\partial t}v(t,x)\;=\;\tfrac{1}{2}\Delta v(t,x)\,+\,c\,(S_{t}%
\varphi)^{1+\gamma}\,(x) \label{equ.lim'}%
\end{equation}
for different values of $\,c.$ Consequently, in a sense, we have to justify
the transition from equation (\ref{equ.f}) to the log-Laplace equation
(\ref{equ.lim'}) corresponding to the limiting fluctuations, recall
(\ref{equ.lim}). Here the $\,x\mapsto\Gamma^{1}(kx)$\thinspace\ entering into
equation (\ref{equ.f}) are random homogeneous fields with infinite overall
density, and the solutions $\,f_{k}$\thinspace\ depend on $\,\Gamma^{1}%
.$\thinspace\ But the most fascinating fact here seems to be the index jump
from $\,2$\thinspace\ to $\,1+\gamma$, which occurs when passing from
\eqref{equ.f} to \eqref{equ.lim'}. Unfortunately, we are unable to explain
this from an individual ergodic theorem acting on the (ergodic) underlying
random measure $\,\Gamma$.

We take another route. For the heuristic exposition, we simplify as follows.
First of all, we restrict our attention to the case $\,\varphi(x)\equiv\theta$
corresponding to total mass process fluctuations. Clearly, we have the
domination%
\begin{equation}
0\,\leq\,u_{k}(t,x)\,\leq\,k^{\varkappa}\,\theta. \label{domi0}%
\end{equation}
Replacing one of the $\,u_{k}(t,x)$\thinspace\ factors in the non-linear term
of (\ref{logLap.equ.k}) by $\,k^{\varkappa}S_{t}\varphi\,(x)\equiv
k^{\varkappa}\,\theta,$\thinspace\ and denoting the solution to the new
equation with the same initial condition by $\,w_{k\,},$\thinspace\ then
$\,u_{k}\geq w_{k\,},$\thinspace\ and we can explicitly calculate $\,w_{k}%
$\thinspace\ by the Feynman-Kac formula,%
\begin{equation}
w_{k}(t,x)\ =\ k^{\varkappa}\,\theta\,\mathcal{E}_{x}\exp\!\Big[-k^{2-d}%
\int_{0}^{t}\!{d}s\;\varrho\,\Gamma^{1}(kW_{s})\,k^{\varkappa}\,\theta\Big].
\label{linprob0}%
\end{equation}
For the upper bound (\ref{fluct1}), we may work with $\,w_{k}$\thinspace
\ instead of $\,u_{k\,}.$\thinspace\ It suffices to show that $\left\langle
\mu,k^{\varkappa}S_{t}\varphi-w_{k}(t,\cdot)\right\rangle $\thinspace
\ converges to $\left\langle \mu,v\right\rangle $ in $L^{2}(\mathsf{P}%
),$\thinspace\ where $v$ is the solution to (\ref{equ.lim'}) with constant
$c=\overline{c}.$ We therefore show that the $\mathsf{P}$--expectations
converge, and the $\mathsf{P}$--variances go to $0.$ In this heuristics we
concentrate on the convergence of $\mathsf{E}$--expectations only, and we
simplify by assuming $\,\mu=\delta_{x}$\thinspace\ (although formally excluded
in the theorem by (\ref{En}) if $\,d>3$\thinspace$).$ We then have to show
that%
\begin{equation}
\mathsf{E}k^{\varkappa}\,\theta\,\mathcal{E}_{x}\biggl(1-\exp
\!\Big[-k^{2-d+\varkappa}\,\theta\int_{0}^{t}\!{d}s\;\varrho\,\Gamma
^{1}(kW_{s})\Big]\biggr)\;\underset{k\uparrow\infty}{\longrightarrow
}\;t\,\overline{c}\,\theta^{1+\gamma}. \label{suff.E}%
\end{equation}
By definition (\ref{not.Gamma.eps}) of $\,\Gamma^{1}$\thinspace\ and
(\ref{logL.Gamma}) of $\,\Gamma,$\thinspace\ the left hand side of
(\ref{suff.E}) can be rewritten as%
\begin{equation}%
\begin{array}
[c]{l}%
\displaystyle
k^{\varkappa}\,\theta\,\mathcal{E}_{x}\biggl(1-\mathsf{E}\exp\!\Big[-\int
\!\Gamma(dz)\ k^{2-d+\varkappa}\,\varrho\theta\int_{0}^{t}\!{d}s\;\vartheta
_{1}(kW_{s}-z)\Big ]\biggr)\vspace{6pt}\\%
\displaystyle
\quad=\ k^{\varkappa}\,\theta\,\mathcal{E}_{x}\biggl(\!1-\exp
\!\bigg[\!-\!k^{(2-d+\varkappa)\gamma+d}\,(\varrho\theta)^{\gamma}%
\!\int\!\!dz\,\Big(\!\int_{0}^{t}\!{d}s\;\mathsf{1}_{B(z,\frac{1}{k})}%
(W_{s})\Big)^{\!\gamma}\bigg]\biggr).
\end{array}
\end{equation}
We may additionally introduce the indicator $\,\mathsf{1}_{\{\tau\leq t\}}%
$\thinspace\ where $\,\tau=\tau_{1/k}^{z}[W]$\thinspace\ denotes the
\emph{first hitting time}\/ of the ball $\,B(z,1/k)$\thinspace\ by the path
$\,W$\thinspace\ starting from $x,$ and we continue with%
\[
=\ k^{\varkappa}\,\theta\,\mathcal{E}_{x}\biggl(\!1-\exp
\!\bigg[\!-\!k^{(2-d+\varkappa)\gamma+d}\,(\varrho\theta)^{\gamma}%
\int\!\!dz\,\mathsf{1}_{\{\tau\leq t\}\,}\Big(\!\int_{0}^{t}\!{d}%
s\;\mathsf{1}_{B(z,\frac{1}{k})}(W_{s})\Big)^{\!\gamma}\bigg]\biggr).
\]
Now we look at the $\mathcal{E}_{x}$--expectation of the exponent term. As the
probability of hitting the small ball $\,B(z,1/k)$\thinspace\ is of order
$\,k^{2-d},$\thinspace\ and the time spent afterwards in the ball is of order
$\,k^{-2},$\thinspace\ the expectation of the exponent term is of order
$\,k^{(-d+\varkappa)\gamma+2}=k^{-\varkappa}$\thinspace\ converging to zero as
$\,k\uparrow\infty$. Heuristically this justifies the use of the approximation
$\,1-\mathrm{e}^{-x}\approx x$.\thinspace\ Note that then the leading factor
$\,k^{\varkappa}$\thinspace\ is cancelled, and we arrive at a constant
multiple of $\,\theta^{1+\gamma}.$\thinspace\ 

According to this simplified calculation, the index jump has its origin in an
averaging of exponential functionals of $\,\Gamma$\thinspace\ [as in
(\ref{logL.Gamma})], generating a transition from $\,\theta$\thinspace\ to
$\,\theta^{\gamma}.$\thinspace\ Note that the smallness of the exponent is
largely due to the presence of the indicator of $\,\{\tau\leq t\}.$%
\thinspace\ This fact is also behind our estimates of variances in
Section~\ref{varlin}.

We recall that the simplification $\,u_{k}\rightsquigarrow w_{k}$%
\thinspace\ which we used in the upper bound is basically a
\emph{linearization}\/ of the problem, that is we pass from the non-linear
log-Laplace equation (\ref{logLap.equ.k}) to the linear equation%
\begin{equation}%
\begin{array}
[c]{c}%
\displaystyle
\frac{\partial}{\partial t}w_{k}(t,x)\;=\;\tfrac{1}{2}\Delta w_{k}%
(t,x)\,-\,k^{2-d}\,\varrho\,\Gamma^{1}(kx)\,k^{\varkappa}\theta\,w_{k}%
(t,x)\vspace{6pt}\\
\text{with initial condition }\,u_{k}(0,\,\cdot\,)=k^{\varkappa}\theta.
\end{array}
\label{equ.w}%
\end{equation}
In the case of a catalyst with finite expectation as in
\cite{DawsonFleischmannGorostiza1989.hydro}, this linearization was a key step
for deriving the limiting fluctuations. The difference between $u_{k}$ and
$w_{k}$ was asymptotically negligible. But in the present model of a catalyst
of infinite overall density, this is \emph{no longer the case.}\/ In fact,
$\,u_{k}(t,x)-w_{k}(t,x)$\thinspace\ does not converge to $0$ in $\mathsf{P}%
$--probability. Therefore, our upper bound is not sharp.

For the lower bound, we replace $u_{k}^{2}$ in (\ref{logLap.equ.k}) by
$w_{k\,}^{2},$ and denoting the solution to the new equation with the same
initial condition by $\,m_{k\,}.$\thinspace\ Then%
\[
k^{\varkappa}\theta-u_{k}(t,x)\ \geq\ k^{\varkappa}\theta-m_{k}%
(t,x)\ =\ k^{2-d}\varrho\,\mathcal{E}_{x}\int_{0}^{t}\!ds\ \Gamma^{1}%
(kW_{s})\,w_{k}^{2}(t-s,W_{s}).
\]
Inserting for $w_{k}$ the Feynman-Kac representation (\ref{linprob0}) we
arrive at an explicit expression. Similarly as above, we then show that
$\left\langle _{\!_{\!_{\,}}}\mu,k^{\varkappa}S_{t}\varphi-m_{k}%
(t,\cdot)\right\rangle $ converges to $\left\langle \mu,v\right\rangle $ in
$L^{2}(\mathsf{P}),$\thinspace\ where $v$ is the solution to (\ref{equ.lim'})
with constant $c=\underline{c}.$

The \emph{structure of the remaining paper}\/ is as follows. After some basic
preparations, in Section~\ref{S.upper} we concentrate on the upper bound,
whereas the lower bound follows in Section~\ref{S.lower}.

\section{Preparation: Some basic estimates\label{S.basic}}

In this section we provide some simple but useful tools for the main body of
the proof. For basic facts on Brownian motion, see, for instance,
\cite{RevuzYor1991} or \cite{KaratzasShreve1991}.

\subsection{Simple estimates for the Brownian semigroup\label{SS.Brown}}

We frequently use the argument (based on the triangle inequality) that, for
$\eta>0$ and $s>0$, there exists $c_{\eqref{klaustrick}}%
=c_{\eqref{klaustrick}}(\eta,s)$ such that for all $x,$
\begin{equation}
\int\!dy\ \phi_{\eta}(y)\,p_{s}(x-y)\ \leq\ \phi_{\eta}(x)\int
\!\!dy\ \mathrm{e}^{\eta|x-y|}\,p_{s}(x-y)\ =\ c_{\eqref{klaustrick}}%
\,\phi_{\eta}(x). \label{klaustrick}%
\end{equation}

For a while, let $\,t>0\,\ $and $\,\varphi\in\mathcal{C}_{\exp\,}^{+}%
.$\thinspace\ Recall that $(s,x)\mapsto S_{s}\varphi\,(x)$ is uniformly
continuous, hence for any $\varepsilon>0$ one may choose $\delta>0$ such that,
for \thinspace$r,s\in\lbrack0,t]\,\ $and $\,\,x,y\in\mathbb{R}^{d}$,
\begin{equation}
\big|S_{r}\varphi\,(x)-S_{s}\varphi\,(y)\big|\,\leq\ \varepsilon\quad
\text{if}\,\ \,|r-s|\leq\delta,\,\ |x-y|\leq\delta. \label{unicont}%
\end{equation}
For convenience we expose the following simple fact.

\begin{lemma}
[\textbf{Brownian semigroup estimate}]\label{L.bo2}There is a $\,\lambda
_{\ref{L.bo2}}=\lambda_{\ref{L.bo2}}(t,\varphi)>0$\thinspace\ and a constant
$\,c_{\ref{L.bo2}}=c_{\ref{L.bo2}}(t,\varphi)$ such that, for every
$x\in\mathbb{R}^{d}$,
\begin{equation}
\tilde{\phi}(x)\;:=\;\sup_{0\leq s\leq t}\,\,\sup_{y\in B(x,1)}\,S_{s}%
\varphi\,(y)\;\leq\;c_{\ref{L.bo2}}\,\phi_{\lambda_{\ref{L.bo2}}}(x).
\label{bo2}%
\end{equation}

\end{lemma}

Note that in all dimensions, for each $\,\lambda>0,$%
\begin{equation}
\sup_{x\in\mathbb{R}^{d}}\,\int\!\!dz\ \phi_{\lambda}(z)\,|z-x|^{2-d}%
\ <\ \infty. \label{all.d}%
\end{equation}
In fact, on the unit ball $\,B(x,1),$\thinspace\ use that $\,\int_{|z|\leq
1}\!dz$\thinspace$|z|^{2-d}<\infty,$\thinspace\ and outside this ball, exploit
$\,|z-x|^{2-d}\leq1$.

We continue with the following observation.

\begin{lemma}
\label{L.sim.id}Let $\,d\geq5.$\thinspace\ Then, for some constant
$\,c_{\ref{L.sim.id}}$\thinspace\ and all $x,y\in\mathbb{R}^{d}$,
\begin{equation}
\int\!\!{d}z\ |z-x|^{2-d}\,|z-y|^{2-d}\ =\ c_{\ref{L.sim.id}}\,|x-y|^{4-d}%
\ =\ c_{\ref{L.sim.id}}\ \mathrm{en}(x-y). \label{da}%
\end{equation}

\end{lemma}%

\proof
Clearly, using the definition of the Green function as an integral of the
transition densities,%
\begin{equation}
\int\!\!{d}z\ |z-x|^{2-d}\,|z-y|^{2-d}\ =\ c\int\!\!{d}z\int_{0}^{\infty
}\!ds\ p_{s}(z-x)\int_{0}^{\infty}\!dt\ p_{t}(z-y). \label{Green}%
\end{equation}
Interchanging integrations, using Chapman-Kolmogorov, substituting, and
interchanging again gives%
\begin{equation}
=\ c\int_{0}^{\infty}\!dt\ t\,p_{t}(x-y)\ =\ c\,|x-y|^{4-d}\int_{0}^{\infty
}\!dt\ t\,p_{t}(\iota)
\end{equation}
with $\,\iota$\thinspace\ any point on the unit sphere. The latter integral is
finite since $\,d>4,$\thinspace\ finishing the proof.%
\endproof

In dimension four, the situation is slightly more involved.

\begin{lemma}
\label{L.sim.est}Let $\,d=4$\thinspace\ and $\,\lambda>0.$\thinspace\ Then,
for some constant $\,c_{\ref{L.sim.est}}=c_{\ref{L.sim.est}}(\lambda
)$\thinspace\ and all \thinspace$x,y\in\mathbb{R}^{4}$,
\begin{equation}
\int\!\!{d}z\ \phi_{\lambda}(z)\,|z-x|^{-2}\,|z-y|^{-2}\ \leq
\ c_{\ref{L.sim.est}}\left[  1+\log^{+}\!\left(  |x-y|^{-1}\right)  \right]
\!. \label{log}%
\end{equation}

\end{lemma}%

\proof
If $\,|x-y|\,>2,$\thinspace\ then the left hand side of (\ref{log}) is bounded
in $\,x,y.$\thinspace\ In fact, for $\,z$\thinspace\ in a unit sphere around a
singularity, say $\,x,$\thinspace\ we use $\,|z-y|\,\geq1$\thinspace\ and
(\ref{all.d}). Outside both unit spheres, the integrand is bounded by
$\,\phi_{\lambda\,}.$

Now suppose $\,|x-y|\,\leq2.$\thinspace\ We may also assume that $\,x\neq
y.$\thinspace\ As in the proof of Lemma~\ref{L.sim.id}, the left hand side of
(\ref{log}) leads to the integral%
\begin{equation}
\int_{0}^{\infty}\!ds\int_{0}^{\infty}\!dt\int\!\!{d}z\ \phi_{\lambda
}(z)\,p_{s}(z-x)\,p_{t}(z-y).
\end{equation}
First we additionally restrict the integrals to $\,s,t\leq|x-y|^{-1}%
.$\thinspace\ In this case, we drop $\,\phi_{\lambda}(z),$\thinspace\ use
Chapman-Kolmogorov, substitute, and interchange the order of integration to
get the bound%
\begin{equation}
\int_{0}^{2\,|x-y|^{-1}}\!dt\ t\,p_{t}(x-y)\ \leq\ \int_{0}^{2\,|x-y|^{-3}%
}\!dt\ t\,p_{t}(\iota)\ \leq\ c\,\left[  1+\log\left(  |x-y|^{-1}\right)
\right]  \!.
\end{equation}
To see the last step, split the integral at $\,t=1.$\thinspace\ To finish the
proof, by symmetry in $\,x,y,$\thinspace\ it suffices to consider%
\begin{equation}
\int_{0}^{\infty}\!ds\int_{|x-y|^{-1}}^{\infty}\!dt\int\!\!{d}z\ \phi
_{\lambda}(z)\,p_{s}(z-x)\,p_{t}(z-y). \label{suff.con}%
\end{equation}
Now, by a substitution,%
\begin{equation}
\int_{|x-y|^{-1}}^{\infty}\!dt\ p_{t}(z-y)\ \leq\ |z-y|^{-2}\int
_{|x-y|^{-1}\,|z-y|^{-2}}^{\infty}\!dt\ c\,t^{-2}\ =\ c\,|x-y|\ \leq\ 2c.
\label{41}%
\end{equation}
Plugging (\ref{41}) into (\ref{suff.con}) and using the Green's function again
gives the bound%
\begin{equation}
c\,\sup_{x\in\mathbb{R}^{4}}\int\!\!{d}z\ \phi_{\lambda}(z)\,|z-x|^{-2},
\end{equation}
which is finite by (\ref{all.d}).%
\endproof

\subsection{Brownian hitting and occupation time estimates\label{hittin}}

Further key tools are the asymptotics of the hitting times of small balls.
Recall that $\,\tau=\tau_{1/k}^{z}[W]$\thinspace\ denotes the \emph{first
hitting time}\/ of the closed ball $\,B(z,1/k)$\thinspace\ by the Brownian
motion $\,W$\thinspace\ started in $\,x.$ The following results are taken from
\cite{LeGall1986.Ann}, see formula (0a) and Lemma~2.1 there.

\begin{lemma}
[\textbf{Hitting time asymptotics and bounds}]\label{L.B.hitting.asy}Suppose
$\,d\geq3.$\thinspace\ Then the following results hold.

\begin{itemize}
\item[(a)] There is a constant $c_{\mathrm{(\ref{hitting})}}$, which depends
only on the dimension $d$, such that
\begin{equation}
\mathcal{P}_{x}(\tau<\infty)\;\leq\;c_{\mathrm{(\ref{hitting})}}%
\,k^{2-d}\,|z-x|^{2-d}\quad\text{for}\,\ x,z\in\mathbb{R}^{d}. \label{hitting}%
\end{equation}

\item[(b)] There are constants $c_{\eqref{bo1}}$ and $\lambda_{\eqref{bo1}}%
>0$, depending on $d$ and $t>0$, such that for $\,x,z\in\mathbb{R}^{d},$
\begin{equation}
k^{d-2}\,\mathcal{P}_{x}(\tau\leq t)\ \leq\ c_{(\ref{bo1})}\big[|z-x|^{2-d}%
+1\big]\exp\bigl[-\lambda_{\eqref{bo1}}|z-x|^{2}\bigr]. \label{bo1}%
\end{equation}

\item[(c)] The following convergence holds uniformly whenever $|x-z|$ is
bounded from zero,
\begin{equation}
\lim_{k\uparrow\infty}\,k^{d-2}\,\mathcal{P}_{x}(\tau\leq t)\ =\ c_{(\ref{ba1}%
)}\int_{0}^{t}\!ds\,\,p_{s}(z-x)\quad\text{for}\,\ z\not =x, \label{ba1}%
\end{equation}
where $\,c_{(\ref{ba1})}:=\frac{(d-2)\pi^{d/2}}{G(d/2)}$ (and $\,G$%
\thinspace\ is the Gamma function).

\item[(d)] Finally, writing $\,\tau^{i}:=\tau_{1/k}^{z_{i}}[W]$ for $i=1,2$,
there are constants $c_{\eqref{loc3}}$ and $\lambda_{\eqref{loc3}}>0$,
depending on $d$ and $t$, such that for $\,x,z\in\mathbb{R}^{d},$%
\begin{equation}%
\begin{array}
[c]{l}%
\mathcal{P}_{x}\!\left(  _{\!_{\!_{\,}}}\tau_{1}<\tau_{2}<{k^{2}t}\right)
\vspace{6pt}\\
\leq\ c_{(\ref{loc3})}\,k^{4-2d}\,\Big(\bigl|(z_{1}-x)/k\bigr|^{2-d}%
+1\Big)\exp\!\Big[-\lambda_{\eqref{loc3}}\bigl|(z_{1}-x)/k\bigr|^{2}%
\Big ]\vspace{6pt}\\
\qquad\qquad\times\ \Big(\big|(z_{2}-z_{1})/k\big|^{2-d}+1\Big)\exp
\!\Big[-\lambda_{\eqref{loc3}}\big|(z_{2}-z_{1})/k\big|^{2}\Big ].
\end{array}
\label{loc3}%
\end{equation}

\end{itemize}
\end{lemma}

The following lemmas are all consequences of Lemma~\ref{L.B.hitting.asy}.

\begin{lemma}
\label{LemmaC} Let $\,d\geq3$. Fix $\,\varphi\in\mathcal{C}_{\mathrm{exp}%
\,}^{+},$\thinspace\ $\eta\geq0,$\thinspace\ and $\,t>0$. Then there are
constants $\,$ $c_{\ref{LemmaC}}$ and $\lambda_{\ref{LemmaC}}$ such that for
$\,x,z\in\mathbb{R}^{d},$
\begin{align}
&  \mathcal{E}_{x}\varphi(W_{t}){\mathsf{1}}_{\{\tau\leq t\}}\Big(k^{2}%
\!\int_{0}^{t}\!{d}s\;\vartheta_{1}(kW_{s}-kz)\Big)^{\!\eta}\ \nonumber\\
&  \leq\ c_{\ref{LemmaC}}k^{2-d}\phi_{\lambda_{\ref{L.bo2}}}%
(z)\,\big[|z-x|^{2-d}+1\big]\exp\bigl[-\lambda_{\ref{LemmaC}}|z-x|^{2}\bigr].
\label{ausdenken}%
\end{align}

\end{lemma}

\proof Initially, let $\varphi$ be any non-negative function. Using the strong
Markov property at time $\,\tau$,
\begin{align}
\mathcal{E}_{x}\varphi(W_{t})\,  &  \,\mathsf{1}_{\{\tau\leq t\}\,}%
\Big(k^{2}\!\int_{0}^{t}\!{d}s\;\vartheta_{1}(kW_{s}-kz)\Big)^{\!\eta
}\nonumber\\
&  =\;\mathcal{E}_{x}\,\varphi(W_{t})\,\mathsf{1}_{\{\tau\leq t\}\,}%
\mathcal{E}_{x}\bigg\{\!\Big(k^{2}\!\int_{0}^{t}\!{d}s\;\vartheta_{1}%
(kW_{s}-kz)\Big)^{\!\eta}\;\bigg|\;\mathcal{F}_{\tau}\bigg\}\label{MP}\\[6pt]
&  =\;\mathcal{E}_{x}\,\varphi(W_{t})\,\mathsf{1}_{\{\tau\leq t\}\,}%
g(\tau,W_{\tau}),\nonumber
\end{align}
where%
\begin{equation}
g(r,y)\;:=\;\mathcal{E}_{y}\Big(k^{2}\!\int_{0}^{t-r}\!{d}s\;\vartheta
_{1}(kW_{s}-kz)\Big)^{\!\eta} \label{not.f}%
\end{equation}
for $\,0\leq r\leq t$\thinspace\ and $\,y\in\partial B(z,1/k).$\thinspace
\ But,
\begin{equation}
g(r,y)\leq\,\mathcal{E}_{y}\Big(\int_{0}^{\infty}\!{d}s\;\vartheta
_{1}(kW_{k^{-2}s}-kz)\Big)^{\!\eta}=\,\mathcal{E}_{ky}\Big(\int_{0}^{\infty
}\!{d}s\ \vartheta_{1}(W_{s}-kz)\Big)^{\!\eta}.
\end{equation}
Note that the right hand side is independent of $\,k,z,y$\thinspace\ (in the
considered range of~$y),$\thinspace\ and finite since in $\,d\geq3$%
\thinspace\ all such moments are finite. Consequently, there is a constant
$\,c$\thinspace\ such that $\,g(r,y)\leq c.$\thinspace\ If now $\,\varphi
\in\mathcal{C}_{\exp\,}^{+},$\thinspace\ by the strong Markov property at time
$\,\tau$,
\begin{equation}
\mathcal{E}_{x}\varphi(W_{t})\,\,\mathsf{1}_{\{\tau\leq t\}\,}=\;\mathcal{E}%
_{x}\,\mathsf{1}_{\{\tau\leq t\}\,}\mathcal{E}_{W_{\tau}}\varphi(\tilde
{W}_{t-\tau})\ \leq\ \mathcal{P}_{x}(\tau\leq t)\,\phi_{\lambda_{\ref{L.bo2}}%
}(z), \label{Markov1}%
\end{equation}
using (\ref{bo2}) in the second step. By \eqref{bo1},
\begin{equation}
\mathcal{P}_{x}(\tau\leq t)\ \leq\ c_{(\ref{bo1})}k^{2-d}\big[|z-x|^{2-d}%
+1\big]\exp\bigl[-\lambda_{(\ref{bo1})}|z-x|^{2}\bigr]. \label{hitt}%
\end{equation}
The result follows by combining \eqref{Markov1} and \eqref{hitt}. \endproof

\begin{lemma}
\label{L.Br.hit.occ}Let $\,d\geq3.$\thinspace\ Fix $\eta\geq0$, $\varphi
\in\mathcal{C}_{\mathrm{exp}}^{+}$, and $t>0$. Then there is a constant
$c_{\mathrm{\ref{L.Br.hit.occ}}}$ such that

\begin{itemize}
\item[(a)] $\displaystyle\mathcal{E}_{x}\Big(k^{2}\!\int_{0}^{t}%
\!{d}s\;\vartheta_{1}(kW_{s}-kz)\Big)^{\!\eta}\;\leq
\;c_{\mathrm{\ref{L.Br.hit.occ}}}\,k^{2-d}\,|z-x|^{2-d},$\newline for all
$x,z\in\mathbb{R}^{d}$ and $k\geq1.$\vspace{2pt}

\item[(b)] $\displaystyle\int\!dz\ \mathcal{E}_{x}\varphi(W_{t}){\mathsf{1}%
}_{\{\tau\leq t\}}\Big(k^{2}\!\int_{0}^{t}\!{d}s\;\vartheta_{1}(kW_{s}%
-kz)\Big)^{\!\eta}\ \leq\ c_{\mathrm{\ref{L.Br.hit.occ}}}k^{2-d}\phi
_{\lambda_{\mathrm{\ref{L.bo2}}}}(x),$\newline for all $x\in\mathbb{R}^{d}$
and $\,k\geq1.$
\end{itemize}
\end{lemma}

\proof The proof of (a) follows from \eqref{MP} for $\varphi\equiv1$ and
\eqref{hitting}, the proof of (b) by integrating (\ref{ausdenken}) and
applying (\ref{klaustrick}).\endproof

\section{Upper bound: Proof of (\ref{fluct1})\label{S.upper}}

\subsection{Anderson model with stable random potential\label{SS.Anderson}}

As motivated in Section~\ref{SS.heuristics}, we look at the mild solution to
the linear equation
\begin{equation}%
\begin{array}
[c]{c}%
\displaystyle
\frac{\partial}{\partial t}w_{k}(t,x)\;=\;\tfrac{1}{2}\Delta w_{k}%
(t,x)-k^{2-d}\,\varrho\,\Gamma^{1}(kx)\,k^{\varkappa}S_{t}\varphi
\,(x)\,w_{k}(t,x)\vspace{6pt}\\
\text{with initial condition }\,w_{k}(0,\,\cdot\,)=k^{\varkappa}\varphi.
\end{array}
\end{equation}
This is an \emph{Anderson model}\/ with the time-dependent scaled stable
random potential\/ $-k^{2-d}\,\varrho\,\Gamma^{1}(kx)\,k^{\varkappa}%
S_{t}\varphi\,(x).$\thinspace\ We study its fluctuation behaviour around the
heat flow:

\begin{proposition}
[\textbf{Limiting fluctuations of }$w_{k}$\thinspace]\label{P.Anderson}%
\hspace{-2.3pt}Under the assumptions of Theorem~\emph{\ref{T.fluct},} if
$\,\varkappa=\varkappa_{\mathrm{c}\,},$\thinspace\ then for any $\,\varphi
\in\mathcal{C}_{\exp}^{+}$\thinspace\ and $\,t\geq0,$\thinspace\ in
$\mathsf{P}$--probability,%
\begin{equation}
\left\langle _{\!_{\!_{\,}}}\mu,k^{\varkappa}S_{t}\varphi-w_{k}(t,\cdot
)\right\rangle \;\underset{k\uparrow\infty}{\longrightarrow}\;\overline{c}\,%
\Big\langle
\mu,\int_{0}^{t}\!dr\ S_{r}\!\left(  _{\!_{\!_{\,}}}(S_{t-r}\varphi
)^{1+\gamma}\right)  \!%
\Big\rangle
, \label{Anderson}%
\end{equation}
where the constant $\overline{c}=\overline{c}(\gamma,\varrho)$\thinspace\ is
given by%
\begin{equation}
\overline{c}\ :=\ \varrho^{\gamma}\,\frac{(d-2)\pi^{d/2}}{G(d/2)}%
\ \mathcal{E}_{\imath}\Big(\int_{0}^{\infty}\!ds\ \vartheta_{1}(W_{s}%
)\Big)^{\!\gamma}, \label{c.fluct}%
\end{equation}
where $\,\imath$\thinspace\ is any point on the unit sphere of \thinspace
$\mathbb{R}^{d}$.
\end{proposition}

To see how the case $n=1$ of (\ref{fluct1}) follows from
Proposition~\ref{P.Anderson}, we fix a sample~$\Gamma.$\thinspace\ For
$\varphi\in\mathcal{C}_{\exp\,}^{+}$, we use the abbreviation
\begin{equation}
\varphi_{k}(x)\;:=\;\varphi(x/k)\quad\text{for}\,\ k>0,\;\,x\in\mathbb{R}^{d}.
\label{not.phik}%
\end{equation}
Formulas (\ref{loglap}) and (\ref{scale}) give
\begin{equation}%
\begin{array}
[c]{l}%
\log\,\mathbb{E}_{\mu_{k}}\!\exp\!\left[  k^{\varkappa}\big(\langle
X_{t\,}^{k},-\varphi\rangle-\langle S_{t}\mu,-\varphi\rangle
\big)_{\!_{\!_{\,_{{}}}}}\right]  \vspace{4pt}\\
\quad=\;\log\mathbb{E}_{\mu_{k}}\!\exp\!\left[  \langle X_{k^{2}%
t\,},-k^{\varkappa-d}\varphi_{k}\rangle+k^{\varkappa}\langle S_{t}\mu
,\varphi\rangle_{\!_{\!_{\,_{{}}}}}\right]  \vspace{8pt}\\
\quad=\ -\left\langle \mu_{k},v_{k}(k^{2}t)_{\!_{\!_{\,}}}\right\rangle
+k^{\varkappa}\langle S_{t}\mu,\varphi\rangle\;=\;\langle\mu,\,k^{\varkappa
}S_{t}\varphi\rangle-\left\langle \mu_{k\,}^{k},\,k^{d}\,v_{k}(k^{2}%
t,k\,\cdot\,)\right\rangle \!,
\end{array}
\label{logLap.k}%
\end{equation}
with $\,v_{k}$\thinspace\ the mild solution to (\ref{logLap.equ}) with initial
condition $\,v_{k}(0)=k^{\varkappa-d}\varphi_{k\,}.$\thinspace\ Setting%
\begin{equation}
u_{k}(t,x)\;:=\;\,k^{d}\,v_{k}(k^{2}t,kx)\quad\text{for}\,\ t\geq
0,\;\,x\in\mathbb{R}^{d}, \label{not.u}%
\end{equation}
$u_{k}$\thinspace\ solves%
\begin{equation}
u_{k}(t,x)\;=\;k^{\varkappa}S_{t}\varphi\,(x)\,-\,k^{2-d}\varrho\int_{0}%
^{t}\!{d}s\;S_{s}\bigl(\Gamma^{1}(k\,\cdot\,)\,u_{k}^{2}(t-s,\,\cdot
\,)\bigr)(x). \label{new}%
\end{equation}
Recall that this can be rewritten in Feynman-Kac form as%
\begin{equation}%
\begin{array}
[c]{l}%
k^{\varkappa}S_{t}\varphi\,(x)-u_{k}(t,x)\vspace{6pt}\\%
\displaystyle
\quad=\;k^{\varkappa}\mathcal{E}_{x}\varphi(W_{t})\bigg(1-\exp\!\Big[-k^{2-d}%
\varrho\int_{0}^{t}\!{d}s\;\Gamma^{1}(kW_{s})\,u_{k}(t-s,W_{s})\Big]\bigg).
\end{array}
\label{FK}%
\end{equation}
Using $\,u_{k}(t-s,W_{s})\leq k^{\varkappa}S_{t-s}\varphi\,(W_{s})$%
\thinspace\ in (\ref{FK}), and the Feynman-Kac representation%
\begin{equation}
w_{k}(t,x)\ :=\ k^{\varkappa}\mathcal{E}_{x}\varphi(W_{t})\exp\!\Big[-k^{2-d}%
\varrho\int_{0}^{t}\!{d}s\;\Gamma^{1}(kW_{s})\,k^{\varkappa}S_{t-s}%
\varphi\,(W_{s})\Big], \label{linprob}%
\end{equation}
we arrive at
\begin{equation}
0\ \leq\ k^{\varkappa}S_{t}\varphi\,(x)-u_{k}(t,x)\ \leq\ k^{\varkappa}%
S_{t}\varphi\,(x)-w_{k}(t,x). \label{76}%
\end{equation}
Hence, the case $n=1$ of (\ref{fluct1}) follows from
Proposition~\ref{P.Anderson}.

Proposition~\ref{P.Anderson} is proved in two steps: In Section~\ref{explin}
we show that the expectations converge, and in Section~\ref{varlin} that the
variances vanish asymptotically.

\subsection{Convergence of expectations\label{explin}}

\begin{proposition}
[\textbf{Convergence of expectations}]\label{L.exp.conv}Let $\,\varkappa
=\varkappa_{\mathrm{c}\,}.$\thinspace\ There exists a $\,\lambda
_{\ref{L.exp.conv}}>0$\thinspace\ such that for every \thinspace
$\varepsilon>0$ there is a $k_{\ref{L.exp.conv}}=k_{\ref{L.exp.conv}%
}(\varepsilon)>0$ with
\begin{equation}
\bigg|\mathsf{E}\!\left(  _{\!_{\!_{\,}}}k^{\varkappa}S_{t}\varphi
\,(x)-w_{k}(t,x)\right)  -\ \overline{c}\int_{0}^{t}\!{d}r\;S_{r}%
(S_{t-r}\varphi)^{1+\gamma}(x)\bigg|\ \leq\ \varepsilon\phi_{\gamma
\lambda_{\ref{L.exp.conv}}}(x)
\end{equation}
for \thinspace$x\in\mathbb{R}^{d},\ \,k\geq k_{\ref{L.exp.conv}\,},$%
\thinspace\ where $\,\overline{c}$\thinspace\ is as in \eqref{c.fluct}.
\end{proposition}

Theorem~\ref{C.LLN} immediately follows from this proposition. Indeed, turning
back to the situation $\,\varkappa<\varkappa_{\mathrm{c}\,},$\thinspace\ note
from (\ref{logLap.k}) (which holds for general $\varkappa)$ that%
\begin{equation}
\log\,\mathbb{E}_{\mu_{k}}\!\exp\!\left[  k^{\varkappa}\langle X_{t\,}%
^{k}-S_{t}\mu,-\varphi\rangle_{\!_{\!_{\,}}}\right]  \;=\;\left\langle
_{\!_{\!_{\,}}}\mu,\,k^{\varkappa}S_{t}\varphi-u_{k}(t,\,\cdot\,)\right\rangle
\,\geq\,0.
\end{equation}
It suffices to show that the right hand side converges to zero in
$L^{1}(\mathsf{P}).$ Using (\ref{76}),
\begin{equation}
\mathsf{E}\!\left\langle _{\!_{\!_{\,}}}\mu,\,k^{\varkappa}S_{t}\varphi
-u_{k}(t,\,\cdot\,)\right\rangle \ \leq\ k^{\varkappa-\varkappa_{\mathrm{c}}%
}\,\mathsf{E}\!\left\langle _{\!_{\!_{\,}}}\mu,\,k^{\varkappa_{\mathrm{c}}%
}S_{t}\varphi-w_{k}(t,\,\cdot\,)\right\rangle \!,
\end{equation}
where $\,w_{k}$\thinspace\ from (\ref{linprob}) is defined using the critical
index $\varkappa_{\mathrm{c}\,}.$\thinspace\ By Proposition~\ref{L.exp.conv},
which does not require the finiteness of the energy, the expectation on the
right remains bounded, implying the statement.\smallskip

The rest of this section is devoted to the proof of this proposition. Recall
that $\,\varkappa\,\ $equals $\,\varkappa_{\mathrm{c}\,},$\thinspace\ which is
defined in (\ref{not.ski}). The proof is prepared by six lemmas. In all these
lemmas, $\tau=\tau_{1/k}^{y}[W]$ denotes the first hitting time of the ball
$B(y,1/k)$ by the Brownian motion $W,$ and $\,\pi_{x}$\thinspace\ the law of
$\,\tau_{1/k}^{y}[W]$\thinspace\ if $\,W$\thinspace\ is started in $\,x.$

\begin{lemma}
\label{highs} There exists a constant \thinspace$c_{\ref{highs}}>0$ such that
\begin{align}
k^{d-2}  &  \int\!\!dy\ \mathcal{E}_{x}1_{\tau\leq t}\,\mathcal{E}_{W_{\tau}%
}\varphi(\tilde{W}_{t-\tau})\Big(k^{2}\int_{M/k^{2}}^{\infty}ds\;\vartheta
_{1/k}(\tilde{W}_{s}-y)\,\phi_{\lambda_{\ref{L.bo2}}}(y)\Big)^{\!\gamma
}\nonumber\\
&  \leq\ c_{\ref{highs}}M^{\gamma(1-d/2)}\phi_{\gamma\lambda_{\ref{L.bo2}}%
}(x)\quad\text{for}\,\ M>1,\ \,k>0,\ \,x\in\mathbb{R}^{d}.
\end{align}

\end{lemma}

\proof Note that, for any $\,\iota\in\partial B(0,1),$\thinspace\ by Brownian
scaling,
\begin{align}
\mathcal{E}_{\iota/k}\,k^{2}  &  \int_{M/k^{2}}^{\infty}ds\;\vartheta
_{1/k}({W}_{s})\,=\ \mathcal{E}_{\iota}\int_{M}^{\infty}ds\;\vartheta_{1}%
({W}_{s})\nonumber\\
&  =\ \int_{M}^{\infty}ds\;\mathcal{P}_{\iota}\!\left(  |{W}_{s}%
|\leq1_{\!_{\!_{\,}}}\right)  \ \leq\ \int_{|y|\leq1}dy\int_{M}^{\infty
}ds\;p_{s}(y)\ \leq\ c_{\eqref{upta}}M^{1-d/2}. \label{upta}%
\end{align}
We now use $\,\varphi\leq c,$\thinspace\ Jensen's inequality, \eqref{upta},
\eqref{bo1}, and \eqref{klaustrick}, to get
\begin{align}
k^{d-2}  &  \int\!\!dy\ \mathcal{E}_{x}\mathsf{1}_{\tau\leq t}\,\mathcal{E}%
_{W_{\tau}}\varphi(\tilde{W}_{t-\tau})\Big(k^{2}\int_{M/k^{2}}^{\infty
}ds\;\vartheta_{1/k}(\tilde{W}_{s}-y)\,\phi_{\lambda_{\ref{L.bo2}}%
}(y)\Big)^{\!\gamma}\nonumber\\
&  \leq\ c\,k^{d-2}\int\!\!dy\ \phi_{\gamma\lambda_{\ref{L.bo2}}%
}(y)\,\mathcal{E}_{x}\mathsf{1}_{\tau\leq t}\,\mathcal{E}_{\iota/k}%
\Big(k^{2}\int_{M/k^{2}}^{\infty}ds\;\vartheta_{1/k}(\tilde{W}_{s}%
)\Big)^{\!\gamma}\label{li2}\\
&  \leq\ c\,M^{\gamma(1-d/2)}\int\!\!dy\ \phi_{\gamma\lambda_{\ref{L.bo2}}%
}(y)\left[  _{\!_{\!_{\,}}}|x-y|^{2-d}+1\right]  \exp\!\left[  -|x-y|^{2}%
/16\right] \nonumber\\[2pt]
&  \leq\ c_{\ref{highs}}\,M^{\gamma(1-d/2)}\phi_{\gamma\lambda_{\ref{L.bo2}}%
}(x).\nonumber
\end{align}
This is the required statement. \endproof

\begin{lemma}
\label{la2} For every $\delta>0$, there exists a constant $c_{\ref{la2}%
}=c_{\ref{la2}}(\delta)>0$ such that
\begin{equation}
\mathcal{E}_{x}\varphi(W_{t})\bigg[\int\!\!dy\,\Big(\int_{0}^{t}%
\!ds\ \vartheta_{1}(kW_{s}-y)\,S_{t-s}\varphi\,(W_{s})\Big)^{\!\gamma
}\bigg]^{2}\,\leq\ c_{\ref{la2}}\,k^{4-4\gamma+\delta}\phi_{\gamma
\lambda_{\ref{L.bo2}}}(x),
\end{equation}
for all $\,x\in\mathbb{R}^{d}\,\ $and $\,k\geq1$.
\end{lemma}

\proof Using Brownian scaling in the second, substitution and \eqref{bo2} in
the last step, we estimate,%
\begin{align}
\mathcal{E}_{x}  &  \varphi(W_{t})\bigg[\int\!\!dy\,\Big(\int_{0}%
^{t}\!ds\ \vartheta_{1}(kW_{s}-y)\,S_{t-s}\varphi\,(W_{s})\Big)^{\!\gamma
}\bigg]^{2}\nonumber\\
&  \leq\ \Vert\varphi\Vert_{\infty}%
{\displaystyle\iint}
\!dy_{1}dy_{2}\ \mathcal{E}_{x}\prod_{i=1}^{2}\Big(\int_{0}^{t}\!ds\ \vartheta
_{1}(kW_{s}-y_{i})\,S_{t-s}\varphi\,(W_{s})\Big)^{\!\gamma}\nonumber\\
&  =\ \Vert\varphi\Vert_{\infty}%
{\displaystyle\iint}
\!dy_{1}dy_{2}\ \mathcal{E}_{0}\prod_{i=1}^{2}\Big(\int_{0}^{t}\!ds\ \vartheta
_{1}(W_{k^{2}s}+kx-y_{i})\,S_{t-s}\varphi\,(\tfrac{1}{k}W_{k^{2}%
s})\Big)^{\!\gamma}\nonumber\\
&  \leq\ k^{-4\gamma}\,\Vert\varphi\Vert_{\infty}%
{\displaystyle\iint}
\!dy_{1}dy_{2}\ \mathcal{E}_{0}\prod_{i=1}^{2}\Big(\int_{0}^{k^{2}%
t}\!ds\ \vartheta_{1}(W_{s}-y_{i})\,\tilde{\phi}(y_{i}/k+x)\Big)^{\!\gamma}.
\label{loc1}%
\end{align}
To study the double integral, denote by $\tau_{1}$, $\tau_{2}$ the first
hitting times of the balls $B(y_{1},1)$ respectively $B(y_{2},1)$ by the
Brownian path $W$. Pick $p>1$ such that $2d+2(2-d)/p<4+\delta$, and $q$ such
that $1/p+1/q=1$. By H\"{o}lder's inequality,%
\begin{gather}
\mathcal{E}_{0}\prod_{i=1}^{2}\Big(\int_{0}^{k^{2}t}\!ds\ \vartheta_{1}%
(W_{s}-y_{i})\,\tilde{\phi}(y_{i}/k+x)\Big)^{\!\gamma}\ \leq\ \left[
\mathcal{P}_{0}\big(\tau_{1}<{k^{2}t},\ \tau_{2}<{k^{2}t}\big)\right]
^{1/p}\nonumber\\
\times\ \bigg[\mathcal{E}_{0}\prod_{i=1}^{2}\Big(\int_{0}^{\infty
}\!ds\ \vartheta_{1}(W_{s}-y_{i})\,\tilde{\phi}(y_{i}/k+x)\Big)^{\!\gamma
q}\bigg]^{\!1/q}.
\end{gather}
For the second factor on the right hand side we get, using Cauchy-Schwarz, and
the maximum principle to pass from $\,y_{i}$\thinspace\ to $\,0,$
\begin{align}
&  \bigg[\mathcal{E}_{0}\prod_{i=1}^{2}\Big(\int_{0}^{\infty}\!ds\ \vartheta
_{1}(W_{s}-y_{i})\,\tilde{\phi}(y_{i}/k+x)\Big)^{\!\gamma q}\bigg]^{\!1/q}%
\nonumber\\
&  \leq\ \prod_{i=1}^{2}\biggl(\mathcal{E}_{0}\Big(\int_{0}^{\infty
}\!ds\ \vartheta_{1}(W_{s}-y_{i})\,\tilde{\phi}(y_{i}/k+x)\Big)^{\!2\gamma
q}\biggr)^{\!1/2q}\nonumber\\
&  \leq\ \tilde{\phi}^{\gamma}(y_{1}/k+x)\,\tilde{\phi}^{\gamma}%
(y_{2}/k+x)\,\biggl(\mathcal{E}_{0}\!\Big(\int_{0}^{\infty}\!ds\ \vartheta
_{1}(W_{s})\Big)^{\!2\gamma q}\biggr)^{\!1/q}. \label{loc4}%
\end{align}
Recall from Lemma~\ref{L.Br.hit.occ}(a) that the total occupation times of
Brownian motion in the unit ball in $d\geq3$ have moments of all orders.
Hence, the latter expectation is finite.

By \eqref{loc3} using substitution in the $y$-variables,
\begin{align}
\iint\!  &  dy_{1}dy_{2}\ \tilde{\phi}^{\gamma}(y_{1}/k+x)\,\tilde{\phi
}^{\gamma}(y_{2}/k+x)\,\left[  \mathcal{P}_{0}\!\left(  \tau_{1}<k^{2}%
t,\ \tau_{2}<k^{2}t\right)  \right]  ^{\!1/p}\nonumber\\
&  \leq\ c_{(\ref{loc3})}^{1/p}\,k^{2d+2(2-d)/p}\int\!\!dy_{1}\,\tilde{\phi
}^{\gamma}(y_{1}+x)\big(|y_{1}|^{2-d}+1\big)^{\!1/p}\exp\bigl[-|y_{1}%
|^{2}/(16p)\bigr]\nonumber\\
&  \qquad\qquad\qquad\quad\quad\times\int\!\!dy_{2}\,\tilde{\phi}^{\gamma
}(y_{2}+x)\big(|y_{2}|^{2-d}+1\big)^{\!1/p}\exp\!\bigl[-|y_{2}|^{2}%
/(16p)\bigr]\nonumber\\
&  \leq\ c_{(\ref{loc2})}k^{4+\delta}\phi_{\gamma\lambda_{\ref{L.bo2}}}(x),
\label{loc2}%
\end{align}
using \eqref{klaustrick} in the last step. Plugging \eqref{loc2} into
\eqref{loc1} completes the proof. \endproof

\begin{lemma}
\label{la5} For all \thinspace$\varepsilon>0$ there exists \thinspace
$\delta=\delta(\varepsilon)>0$ and \thinspace$k_{\ref{la5}}=k_{\ref{la5}%
}(\varepsilon)>0$, such that
\begin{gather}
k^{d-2}\int\!\!dy\,\ \mathcal{E}_{x}\mathsf{1}_{t-\delta\leq\tau\leq
t}\,\mathcal{E}_{W_{\tau}}\!\Big(k^{2}\int_{0}^{t-\tau}\!{d}s\;\vartheta
_{1/k}(\tilde{W}_{s}-y)\,S_{t-\tau-s}\varphi\,(\tilde{W}_{s})\Big)^{\!\gamma
}\nonumber\\
\leq\ \varepsilon\phi_{\gamma\lambda_{\ref{L.bo2}}}(x)\quad for\text{
\thinspace$k\geq k$}_{\ref{la5}}\text{ and\thinspace\ $x\in\mathbb{R}^{d}$.}%
\end{gather}

\end{lemma}

\proof For any $\delta,M>0$ we have,
\begin{subequations}
\begin{align}
k^{d-2}  &  \int\!\!dy\,\ \mathcal{E}_{x}\mathsf{1}_{t-\delta\leq\tau\leq
t}\,\mathcal{E}_{W_{\tau}}\!\Big(k^{2}\int_{0}^{t-\tau}\!\!\!\!\!\!{d}%
s\;\vartheta_{1/k}(\tilde{W}_{s}-y)\,S_{t-\tau-s}\varphi\,(\tilde{W}%
_{s})\Big)^{\!\gamma}\nonumber\\
&  \leq\ k^{d-2}\int\!\!dy\ \phi_{\gamma\lambda_{\ref{L.bo2}}}(y)\,\mathcal{E}%
_{x}\mathsf{1}_{t-\delta\leq\tau\leq t}\,\mathcal{E}_{W_{\tau}}\!\Big(k^{2}%
\int_{0}^{M/k^{2}}\!{d}s\;\vartheta_{1/k}(\tilde{W}_{s}-y)\Big)^{\!\gamma
}\label{te1}\\
&  \qquad\ +\ k^{d-2}\int\!\!dy\ \phi_{\gamma\lambda_{\ref{L.bo2}}%
}(y)\,\mathcal{E}_{x}\mathsf{1}_{\tau\leq t}\,\mathcal{E}_{W_{\tau}%
}\!\Big(k^{2}\int_{M/k^{2}}^{\infty}\!{d}s\;\vartheta_{1/k}(\tilde{W}%
_{s}-y)\Big)^{\!\gamma}. \label{te2}%
\end{align}
We look at \eqref{te2} and choose $M$ such that this term is small. Indeed,
the inner expectation in \eqref{te2} can be made arbitrarily small
(simultaneously for all $k$ and $y$) by choice of $M$. Hence we can use
\eqref{bo1} to see that this term can be bounded by $\varepsilon\phi
_{\gamma\lambda_{\ref{L.bo2}}}(x)$, for all sufficiently large~$k$, by choice
of $M$ (and independently of $\delta$).

We look at \eqref{te1} and choose $\delta>0$ such that
\end{subequations}
\begin{equation}
c_{(\ref{ba1})}M^{\gamma}\int_{t-\delta}^{t}\!ds\int\!\!dy\ \phi
_{\gamma\lambda_{\ref{L.bo2}}}(y)\,p_{s}(y-x)\ <\ \varepsilon\phi
_{\gamma\lambda_{\ref{L.bo2}}}(x). \label{choicedelta}%
\end{equation}
The term~\eqref{te1} can be bounded from above by
\begin{equation}
M^{\gamma}k^{d-2}\int\!\!dy\ \phi_{\gamma\lambda_{\ref{L.bo2}}}(y)\,\pi
_{x}[t-\delta,t]. \label{simpl}%
\end{equation}
By \eqref{ba1} there exists $A\subset\mathbb{R}^{d}$ and $k_{\ref{la5}}\geq0$
such that, for all $x-y\in A$ and $k\geq k_{\ref{la5}}$,
\begin{equation}
k^{d-2}\pi_{x}[t-\delta,t]\,-\,c_{(\ref{ba1})}\int_{t-\delta}^{t}%
\!ds\,p_{s}(y-x)\ <\ \varepsilon\,\int_{0}^{t}ds\,p_{s}(y-x)
\end{equation}
and
\begin{equation}
\int_{A^{c}}dz\,\big[|z|^{2-d}+1\big]\exp\bigl[\lambda_{\ref{L.bo2}%
}|z|-|z|^{2}/16\bigr]\ <\ \varepsilon. \label{for}%
\end{equation}
We can thus bound \eqref{simpl}, for all $k\geq k_{\ref{la5}}$ and
$x\in\mathbb{R}^{d}$ by
\begin{align*}
&  M^{\gamma}k^{d-2}\int\!\!dy\ \phi_{\gamma\lambda_{\ref{L.bo2}}}(y)\,\pi
_{x}[t-\delta,t]\ \leq\ c_{(\ref{ba1})}M^{\gamma}\int_{x+A}\!\!dy\ \phi
_{\gamma\lambda_{\ref{L.bo2}}}(y)\,\int_{t-\delta}^{t}\!ds\,p_{s}(y-x)\\
&  \ \quad+\ \varepsilon\,M^{\gamma}\int\!\!dy\ \phi_{\gamma\lambda
_{\ref{L.bo2}}}(y)\,\int_{0}^{t}ds\,p_{s}(y-x)\,+\,M^{\gamma}\int_{x+A^{c}%
}\!\!dy\ \phi_{\gamma\lambda_{\ref{L.bo2}}}(y)\,\,k^{d-2}\pi_{x}[0,t].
\end{align*}
By \eqref{choicedelta} the first term is bounded by $\varepsilon\phi
_{\gamma\lambda_{\ref{L.bo2}}}(x)$, as is the second term. For the last term
we use the upper bound \eqref{bo1} for $k^{2-d}\pi_{x}[0,t]$ and then
\eqref{for} to see the upper bound of $\varepsilon\phi_{\gamma\lambda
_{\ref{L.bo2}}}(x)$. \endproof

\begin{lemma}
\label{replace1} For every $M>1$ and \thinspace$\varepsilon>0$, there exists a
\thinspace$k_{\ref{replace1}}=k_{\ref{replace1}}(M,\varepsilon)>0$%
\thinspace\ such that
\begin{align*}
&  k^{d-2}\varrho^{\gamma}\int\!\!dy\,\ \mathcal{E}_{x}\mathsf{1}_{\tau\leq
t}\,\mathcal{E}_{W_{\tau}}\varphi(\tilde{W}_{t-\tau})\bigg|\Big(k^{2}\int
_{0}^{M/k^{2}}ds\,\vartheta_{1/k}(\tilde{W}_{s}-y)S_{t-\tau-s}\varphi
\,(\tilde{W}_{s})\Big)^{\!\gamma}\\
&  -\Big(k^{2}\int_{0}^{M/k^{2}}ds\,\vartheta_{1/k}(\tilde{W}_{s}-y)S_{t-\tau
}\varphi\,(y)\Big)^{\!\gamma}\bigg|\ \leq\ \varepsilon\phi_{\gamma
\lambda_{\ref{L.bo2}}}(x)\quad\text{for }\,k\geq k_{\ref{replace1}\,}%
,\ \,x\in\mathbb{R}.
\end{align*}

\end{lemma}

\proof Recall that $|a^{\gamma}-b^{\gamma}|\leq|a-b|^{\gamma}$. We use
\eqref{unicont} to choose $k_{\ref{replace1}}>1/M$ such that
\begin{equation}
\big|S_{r}\varphi\,(x)-S_{s}\varphi\,(y)\big|\leq\varepsilon^{1/\gamma}%
\quad\text{if}\,\ |r-s|\leq M/k_{\ref{replace1}\,}^{2},\,\,\,|x-y|\leq
1/k_{\ref{replace1}\,}. \label{unicont2}%
\end{equation}
Hence, for all $k\geq k_{\ref{replace1}}$ and $x\in\mathbb{R}^{d}$,
\begin{align*}
k^{d-2}  &  \varrho^{\gamma}\int\!\!dy\,\ \mathcal{E}_{x}\mathsf{1}_{\tau\leq
t}\mathcal{E}_{W_{\tau}}\varphi(\tilde{W}_{t-\tau})\bigg|\Big(k^{2}\int
_{0}^{M/k^{2}}ds\,\vartheta_{1/k}(\tilde{W}_{s}-y)S_{t-\tau-s}\varphi
\,(\tilde{W}_{s})\Big)^{\!\gamma}\\
&  \qquad\qquad\qquad\qquad\qquad\quad\quad\quad-\Big(k^{2}\int_{0}^{M/k^{2}%
}ds\,\vartheta_{1/k}(\tilde{W}_{s}-y)S_{t-\tau}\varphi\,(y)\Big)^{\!\gamma
}\bigg|\\
&  \leq\ k^{d-2}\varrho^{\gamma}\int\!\!dy\,\ \mathcal{E}_{x}\mathsf{1}%
_{\tau\leq t}\,\mathcal{E}_{W_{\tau}}\varphi(\tilde{W}_{t-\tau})\Big(k^{2}%
\int_{0}^{M/k^{2}}ds\,\vartheta_{1/k}(\tilde{W}_{s}-y)\\
&  \qquad\qquad\qquad\qquad\qquad\qquad\quad\quad\quad\quad\,\ \times
\big|S_{t-\tau-s}\varphi\,(\tilde{W}_{s})-S_{t-\tau}\varphi\,(\tilde{W}%
_{s})\big|\Big)^{\!\gamma}\\
&  \leq\ \varepsilon k^{d-2}\varrho^{\gamma}\int\!\!dy\,\ \mathcal{E}%
_{x}\mathsf{1}_{\tau\leq t}\,\mathcal{E}_{W_{\tau}}\varphi(\tilde{W}_{t-\tau
})\Big(k^{2}\int_{0}^{M/k^{2}}ds\,\vartheta_{1/k}(\tilde{W}_{s}%
-y)\Big)^{\!\gamma}.
\end{align*}
To complete the proof use Cauchy-Schwarz, \eqref{bo2}, \eqref{bo1}, and
\eqref{klaustrick}, to get
\begin{subequations}
\begin{align}
k^{d-2}  &  \int\!\!dy\ \mathcal{E}_{x}\mathsf{1}_{\tau\leq t}\,\mathcal{E}%
_{W_{\tau}}\varphi(\tilde{W}_{t-\tau})\Big(k^{2}\int_{0}^{M/k^{2}%
}ds\;\vartheta_{1/k}(\tilde{W}_{s}-y)\Big)^{\!\gamma}\nonumber\\
&  \leq\ k^{d-2}\int\!\!dy\ \mathcal{E}_{x}\mathsf{1}_{\tau\leq t}%
\,\big(\mathcal{E}_{W_{\tau}}\varphi^{2}(\tilde{W}_{t-\tau})\big)^{1/2}%
\nonumber\\
&  \qquad\qquad\qquad\quad\ \ \times\Big[\mathcal{E}_{W_{\tau}}\Big(k^{2}%
\int_{0}^{M/k^{2}}ds\;\vartheta_{1/k}(\tilde{W}_{s}-y)\Big)^{\!2\gamma
}\Big ]^{1/2}\nonumber\\
&  \leq\ c_{\eqref{lin2}}k^{d-2}\int\!\!dy\ \phi_{\lambda_{\ref{L.bo2}}%
}(y)\,\mathcal{E}_{x}\mathsf{1}_{\tau\leq t}\,\Big[\mathcal{E}_{0}%
\Big(\int_{0}^{\infty}ds\;\vartheta_{1}({W}_{s})\Big)^{\!2\gamma}%
\Big ]^{1/2}\label{lin2}\\
&  \leq\ c_{\eqref{lin3}}\int\!\!dy\ \phi_{\lambda_{\ref{L.bo2}}}(y)\left[
_{\!_{\!_{\,}}}|x-y|^{2-d}+1\right]  \exp\!\left[  -|x-y|^{2}/16\right]
\label{lin3}\\[4pt]
&  \leq\ c_{\eqref{lin4}}\phi_{\gamma\lambda_{\ref{L.bo2}}}(x). \label{lin4}%
\end{align}
This gives the required statement. \endproof
\end{subequations}
\begin{lemma}
\label{replace2} Let $\,M>1$\thinspace\ and $\,c_{\ref{replace2}%
}=c_{\ref{replace2}}(M):=\mathcal{E}_{\iota}\,\big\{\!(\int_{0}^{M}%
\!ds\ \vartheta_{1}(W_{s}))^{\!\gamma}\big\}$ for $\iota\in\partial B(0,1)$.
For every $\varepsilon>0$ there exists a $k_{\ref{replace2}}=k_{\ref{replace2}%
}(\varepsilon,M)>0$\thinspace\ such that, for all $k\geq k_{\ref{replace2}}$
and $x\in\mathbb{R}^{d}$,
\begin{align*}
k^{d-2}  &  \varrho^{\gamma}\int\!\!dy\,\ \mathcal{E}_{x}\mathsf{1}_{\tau\leq
t}\,\!\left(  _{\!_{\!_{\,}}}S_{t-\tau}\varphi\,(y)\right)  ^{\!\gamma}\\
&  \times\bigg|\mathcal{E}_{W_{\tau}}\Big(k^{2}\int_{0}^{M/k^{2}}%
ds\,\vartheta_{1/k}(\tilde{W}_{s}-y)\Big)^{\!\gamma}\mathcal{E}_{\tilde
{W}_{M/k^{2}}}\varphi(\tilde{\tilde{W}}_{t-\tau-M/k^{2}})-c_{\ref{replace2}%
}S_{t-\tau}\varphi\,(y)\bigg|\\
\leq\ \varepsilon &  \phi_{\gamma\lambda_{\ref{L.bo2}}}(x),
\end{align*}

\end{lemma}

\proof In a first step we note that, by Brownian scaling,
\begin{align*}
\mathcal{E}_{W_{\tau}}  &  \Big(k^{2}\int_{0}^{M/k^{2}}ds\,\vartheta
_{1/k}(\tilde{W}_{s}-y)\Big)^{\!\gamma}\mathcal{E}_{\tilde{W}_{M/k^{2}}%
}\varphi(\tilde{\tilde{W}}_{t-\tau-M/k^{2}})\\
&  =\ \mathcal{E}_{0}\Big(k^{2}\int_{0}^{M/k^{2}}ds\,\vartheta_{1/k}(\tfrac
{1}{k}\tilde{W}_{sk^{2}}-y+W_{\tau})\Big)^{\!\gamma}\mathcal{E}_{W_{\tau
}+\frac{1}{k}\tilde{W}_{M}}\varphi(\tilde{\tilde{W}}_{t-\tau-M/k^{2}}).
\end{align*}
The main contribution to this expectation is coming from those $\tilde{W}$
with $\tilde{W}_{M}\leq\sqrt{k}$. Indeed, the remaining part of the integral
can be estimated by a constant multiple of $M^{\gamma}\mathcal{P}%
_{0}\big\{\tilde{W}_{M}>\sqrt{k}\big\}$, and we can estimate (with
$c_{\eqref{mue1}}$ depending on $M$)
\begin{align}
&  k^{d-2}\int\!dy\,\ \mathcal{E}_{x}\Big(\mathsf{1}_{\tau\leq t}\,\!\left(
_{\!_{\!_{\,}}}S_{t-\tau}\varphi\,(y)\right)  ^{\!\gamma}\Big)M^{\gamma
}\,\mathcal{P}_{0}\bigl(\tilde{W}_{M}>\sqrt{k}\bigr)\nonumber\\
&  \leq\ c_{\eqref{mue1}}\mathrm{e}^{-k/2M}\int\!dy\ \phi_{\gamma
\lambda_{\ref{L.bo2}}}(y)\big[|x-y|^{2-d}+1\big]\exp\!\big[\lambda
_{\ref{L.bo2}}|x-y|-|x-y|^{2}/16\big]\label{mue1}\\[4pt]
&  \leq\ \varepsilon\phi_{\gamma\lambda_{\ref{L.bo2}}}(x),\nonumber
\end{align}
for sufficiently large values of $k$, recalling \eqref{bo1} and \eqref{klaustrick}.

In the next step we use \eqref{unicont} to choose $k$ large enough such that
\begin{equation}
\big|S_{r}\varphi\,(w+z)-S_{s}\varphi\,(y)\big|\leq\varepsilon\quad\text{if
}\,|r-s|\leq M/k^{2},\ |z|\leq1/\sqrt{k},\ |w-y|\leq1/k.
\end{equation}
Using this,
\begin{align*}
k^{d-2}  &  \int\!dy\,\ \mathcal{E}_{x}\mathsf{1}_{\tau\leq t}\,\!\left(
_{\!_{\!_{\,}}}S_{t-\tau}\varphi\,(y)\right)  ^{\!\gamma}\mathcal{E}_{W_{\tau
}}\mathsf{1}_{\tilde{W}_{M}<\sqrt{K}}\,\Big(k^{2}\int_{0}^{M/k^{2}%
}ds\,\vartheta_{1/k}(\tilde{W}_{s}-y)\Big)^{\!\gamma}\\
&  \qquad\qquad\qquad\qquad\qquad\qquad\quad\times\Big|\mathcal{E}_{W_{\tau
}+\tfrac{1}{k}\tilde{W}_{M}}\varphi(\tilde{\tilde{W}}_{t-\tau-M/k^{2}%
})-\mathcal{E}_{y}\varphi(\tilde{\tilde{W}}_{t-\tau})\Big|\\
&  \leq\ \varepsilon k^{d-2}\int\!dy\,\ \phi_{\gamma\lambda_{\ref{L.bo2}}%
}(y)\,\mathcal{E}_{x}\mathsf{1}_{\tau\leq t}\,\mathcal{E}_{W_{\tau}%
}\!\Big(k^{2}\int_{0}^{M/k^{2}}ds\,\vartheta_{1/k}(\tilde{W}_{s}%
-y)\Big)^{\!\gamma}\\
&  \leq\ \varepsilon\int\!dy\,\ \phi_{\gamma\lambda_{\ref{L.bo2}}}%
(y)k^{d-2}\mathcal{P}_{x}\left(  \tau\leq t\right)  \,\mathcal{E}_{0}%
\Big(\int_{0}^{\infty}ds\,\vartheta_{1}(\tilde{W}_{s})\Big)^{\!\gamma},
\end{align*}
and the last line is $\ \leq\varepsilon\,\phi_{\gamma\lambda_{\ref{L.bo2}}%
}(x)$\thinspace\ by \eqref{bo1} and \eqref{klaustrick}.

Now it remains to observe that, by Brownian scaling,
\begin{gather*}
k^{d-2}\int\!dy\,\ \mathcal{E}_{x}\mathsf{1}_{\tau\leq t}\,\!\left(
_{\!_{\!_{\,}}}S_{t-\tau}\varphi\,(y)\right)  ^{\!\gamma}\mathcal{E}_{W_{\tau
}}\!\Big(k^{2}\int_{0}^{M/k^{2}}ds\,\vartheta_{1/k}(\tilde{W}_{s}%
-y)\Big)^{\!\gamma}\mathcal{E}_{y}\varphi(\tilde{\tilde{W}}_{t-\tau})\\
=\ k^{d-2}\int\!dy\ \mathcal{E}_{x}\mathsf{1}_{\tau\leq t}\left(
_{\!_{\!_{\,}}}S_{t-\tau}\varphi\,(y)\right)  ^{\!1+\gamma}\mathcal{E}%
_{kW_{\tau}}\Big(\int_{0}^{M}ds\,\vartheta_{1}(\tilde{W}_{s}-y)\Big)^{\!\gamma
}.
\end{gather*}
For $y\not \in B(x,1/k)$ the inner expectation is constant and equals
$c_{\ref{replace2}}$. The contribution coming from $y\in B(x,1/k)$ is very
easily seen to be bounded by a constant multiple of $k^{-2}\phi_{\gamma
\lambda_{\ref{L.bo2}}}(x)$. This completes the proof. \endproof

The following lemma is at the heart of our proof of
Proposition~\ref{L.exp.conv}. Recall that $\,\pi_{x}$\thinspace\ denotes the
law of $\,\tau=\tau_{1/k}^{y}[W]$\thinspace\ for $\,W$\thinspace\ starting in
$\,x.$

\begin{lemma}
[\textbf{A hitting time statement}]\label{la4}For every $\varepsilon>0$ there
exists a $k_{\ref{la4}}=k_{\ref{la4}}(\varepsilon)>0$\thinspace\ such that
\begin{gather}
\bigg|k^{d-2}\int\!dy\int_{0}^{t}\pi_{x}(ds)\big(S_{t-s}\varphi
\,(y)\big)^{1+\gamma}-c_{(\ref{ba1})}\int_{0}^{t}\!ds\ S_{s}(S_{t-s}%
\varphi)^{1+\gamma}(x)\bigg|\nonumber\\
\leq\ \varepsilon\phi_{\lambda_{\ref{L.bo2}}}(x)\quad\text{for\thinspace
\ }x\in\mathbb{R}^{d},\ \,k\geq k_{\ref{la4}\,}.
\end{gather}

\end{lemma}

\proof Fix $\varepsilon>0$. Recall that $(s,y)\mapsto S_{s}\varphi\,(y)$ is
uniformly continuous and bounded, and that there exists $R>0$ (dependent on
$\varepsilon$) such that $\left(  _{\!_{\!_{\,}}}S_{s}\varphi\,(y)\right)
^{1+\gamma}\leq\varepsilon\phi_{\lambda_{\ref{L.bo2}}}(y)$ for all $0\leq
s\leq t$, $|y|>R$. We can therefore choose $0=t_{0}\leq\dots\leq t_{n}=t$ such
that, for all $t_{j}\leq r,s\leq t_{j+1}$ and $y\in\mathbb{R}^{d}$,
\begin{equation}
\left\vert _{\!_{\!_{\,_{{}}}}}\!\left(  _{\!_{\!_{\,}}}S_{t-s}\varphi
\,(y)\right)  ^{1+\gamma}-\left(  _{\!_{\!_{\,}}}S_{t-r}\varphi\,(y)\right)
^{1+\gamma}\right\vert \ \leq\ \varepsilon\phi_{\lambda_{\ref{L.bo2}}}(y).
\label{unif}%
\end{equation}
Using \eqref{ba1} we may find $k_{\ref{la4}}$ such that, for all $k\geq
k_{\ref{la4}}$,
\begin{equation}
\Big|k^{d-2}\pi_{x}[t_{j},t_{j+1}]-\int_{t_{j}}^{t_{j+1}}\!ds\,p_{s}%
(x,y)\Big|\ <\ \varepsilon\,k^{d-2}\pi_{x}[t_{j},t_{j+1}] \label{pieces}%
\end{equation}
for all $\,0\leq j\leq n-1$\thinspace\ and all $\,x-y\in A,$\thinspace\ where
$A\subset\mathbb{R}^{d}$\thinspace\ is a set with
\begin{equation}
\int_{A^{c}}dz\,\big[|z|^{2-d}+1\big]\exp\bigl[\lambda_{\ref{L.bo2}%
}|z|-|z|^{2}/16\bigr]\ <\ \varepsilon. \label{forcrit}%
\end{equation}
Now we show that for all $x\in\mathbb{R}^{d}$, and $k\geq k_{\ref{la4}}$,
\begin{align}
k^{d-2}\int\!dy\int_{0}^{t}  &  \pi_{x}(ds)\big(S_{t-s}\varphi
\,(y)\big)^{1+\gamma}\nonumber\\
&  \leq\ c_{(\ref{ba1})}\int_{0}^{t}\!ds\ \big(S_{s}(S_{t-s}\varphi
)^{1+\gamma}\big)(x)+\varepsilon\,\phi_{\lambda_{\ref{L.bo2}}}(x).
\label{aim1}%
\end{align}
Indeed, using \eqref{unif} and \eqref{pieces}, we can estimate
\begin{subequations}
\begin{align}
k^{d-2}  &  \int\!dy\int_{0}^{t}\pi_{x}(ds)\big(S_{t-s}\varphi
\,(y)\big)^{1+\gamma}\nonumber\\
&  \leq\ k^{d-2}\int\!dy\sum_{j=0}^{n-1}\Big[\!\left(  _{\!_{\!_{\,}}%
}S_{t-t_{j}}\varphi\,(y)\right)  ^{\!1+\gamma}+\varepsilon\phi_{\lambda
_{\ref{L.bo2}}}(y)\Big ]\pi_{x}[t_{j},t_{j+1}]\nonumber\\
&  \leq\ \int\!\!dy\sum_{j=0}^{n-1}\Big[\!\left(  _{\!_{\!_{\,}}}S_{t-t_{j}%
}\varphi\,(y)\right)  ^{\!1+\gamma}+\varepsilon\phi_{\lambda_{\ref{L.bo2}}%
}(y)\Big ]\int_{t_{j}}^{t_{j+1}}\!ds\,p_{s}(x,y)\label{mainterm}\\
&  \quad\ \ +\ \varepsilon\int\!\!dy\sum_{j=0}^{n-1}\Big[\!\left(
_{\!_{\!_{\,}}}S_{t-t_{j}}\varphi\,(y)\right)  ^{\!1+\gamma}+\varepsilon
\phi_{\lambda_{\ref{L.bo2}}}(y)\Big ]\,k^{d-2}\pi_{x}[t_{j},t_{j+1}%
]\label{erterm1}\\
&  \quad\ \ +\ \int_{x+A^{c}}dy\sum_{j=0}^{n-1}\Big[\!\left(  _{\!_{\!_{\,}}%
}S_{t-t_{j}}\varphi\,(y)\right)  ^{\!1+\gamma}+\varepsilon\phi_{\lambda
_{\ref{L.bo2}}}(y)\Big ]\,k^{d-2}\pi_{x}[t_{j},t_{j+1}]. \label{erterm3}%
\end{align}
We give estimates for the two final summands, the error terms. The term
\eqref{erterm1} can be estimated, using \eqref{bo1}, by
\end{subequations}
\begin{align}
\varepsilon\!\!\int\!\!dy\sum_{j=0}^{n-1}  &  \Big[\!\left(  _{\!_{\!_{\,}}%
}S_{t-t_{j}}\varphi\,(y)\right)  ^{\!1+\gamma}+\varepsilon\phi_{\lambda
_{\ref{L.bo2}}}(y)\Big ]k^{d-2}\pi_{x}[t_{j},t_{j+1}]\,\leq\,2\varepsilon
\!\!\int\!\!dy\ \phi_{\lambda_{\ref{L.bo2}}}(y)k^{d-2}\pi_{x}[0,t]\nonumber\\
&  \leq\ 2\varepsilon\,c_{(\ref{bo1})}\int\!\!dy\ \phi_{\lambda_{\ref{L.bo2}}%
}(y)\big[|x-y|^{2-d}+1\big]\exp\!\big[-|x-y|^{2}/16\big]\nonumber\\[8pt]
&  \leq\ \varepsilon\,c_{\eqref{feh1}}\phi_{\lambda_{\ref{L.bo2}}}(x).
\label{feh1}%
\end{align}
The error term \eqref{erterm3} can be estimated as follows,
\begin{align}
\int_{x+A^{c}}  &  dy\,\sum_{j=0}^{n-1}\Big[\!\left(  _{\!_{\!_{\,}}%
}S_{t-t_{j}}\varphi\,(y)\right)  ^{\!1+\gamma}+\varepsilon\phi_{\lambda
_{\ref{L.bo2}}}(y)\Big ]\,k^{d-2}\pi_{x}[t_{j},t_{j+1}]\nonumber\\
&  \leq\ c_{(\ref{bo1})}\int_{x+A^{c}}dy\ \phi_{\lambda_{\ref{L.bo2}}%
}(y)\big[|x-y|^{2-d}+1\big]\exp\!\big[-|x-y|^{2}/16\big]\nonumber\\
&  \leq\ c_{\eqref{crit}}\,\phi_{\lambda_{\ref{L.bo2}}}(x)\int_{x+A^{c}%
}dy\,\big[|x-y|^{2-d}+1\big]\exp\!\big[\lambda_{\ref{L.bo2}}|x-y|-|x-y|^{2}%
/16\big], \label{crit}%
\end{align}
and the integral is smaller than $\varepsilon$ by \eqref{forcrit}. For the
first summand, the main term~\eqref{mainterm}, we argue that
\begin{align}
\int\!\!dy  &  \sum_{j=0}^{n}\Big[\!\left(  _{\!_{\!_{\,}}}S_{t-t_{j}}%
\varphi\,(y)\right)  ^{\!1+\gamma}+\varepsilon\phi_{\lambda_{\ref{L.bo2}}%
}(y)\Big ]\int_{t_{j}}^{t_{j+1}}\!ds\,p_{s}(x,y)\\
&  \leq\ \int\!\!dy\,\int_{0}^{t}ds\,\big[(S_{t-s}\varphi\,(y))^{1+\gamma
}+2\varepsilon\,\phi_{\lambda_{\ref{L.bo2}}}(y)\big]p_{s}(x,y)\nonumber\\
&  \leq\ \int_{0}^{t}ds\,S_{s}\big((S_{t-s}\varphi\,(x))^{1+\gamma
}\big)+2\varepsilon\int\!\!dy\ \phi_{\lambda_{\ref{L.bo2}}}(y)\,\int_{0}%
^{t}ds\,p_{s}(x,y).\nonumber
\end{align}
The last summand is again bounded by a constant multiple of $\varepsilon
\phi_{\lambda_{\ref{L.bo2}}}(x)$. Hence we have verified \eqref{aim1} and by
the analogous argument one can see that, for all $k\geq k_{\ref{la4}}$ and
$x\in\mathbb{R}^{d}$,
\[
k^{d-2}\int\!dy\int_{0}^{t}\pi_{x}(ds)\big(S_{t-s}\varphi\,(y)\big)^{1+\gamma
}\ \geq\ c_{(\ref{ba1})}\int_{0}^{t}\!ds\ \big(S_{s}(S_{t-s}\varphi
)^{1+\gamma}\big)(x)-\varepsilon\,\phi_{\lambda_{\ref{L.bo2}}}(x).
\]
This completes the proof.\endproof

\noindent\emph{Proof of Proposition}\/ \ref{L.exp.conv}. \thinspace Recall
from (\ref{linprob}) that%
\begin{align}
&  \mathsf{E}\!\left(  _{\!_{\!_{\,}}}k^{\varkappa}S_{t}\varphi\,(x)-w_{k}%
(t,x)\right)  \ \nonumber\\
&  =\ k^{\varkappa}\,\mathsf{E}\mathcal{E}_{x}\varphi(W_{t})\bigg(1-\exp
\!\Big [-k^{2-d+\varkappa}\varrho\int_{0}^{t}\!{d}s\;\Gamma^{1}(kW_{s}%
)\,S_{t-s}\varphi\,(W_{s})\Big]\bigg).
\end{align}
We use \eqref{logL.Gamma} to evaluate the expectation with respect to the
medium.
\begin{align}
\hspace*{5pt}  &  \mathsf{E}k^{\varkappa}\mathcal{E}_{x}\varphi(W_{t}%
)\bigg(1-\exp\!\Big [-k^{2-d+\varkappa}\varrho\int_{0}^{t}\!{d}s\;\Gamma
^{1}(kW_{s})\,S_{t-s}\varphi\,(W_{s})\Big]\bigg)\label{express}\\
&  =\,k^{\varkappa}\mathcal{E}_{x}\varphi(W_{t})\nonumber\\
&  \quad\ \times\bigg(1-\exp\!\Big [-k^{(2+\varkappa-d)\gamma}\varrho^{\gamma
}\int\!\!dy\Big(\int_{0}^{t}\!{d}s\;\vartheta_{1}(kW_{s}-y)\,S_{t-s}%
\varphi\,(W_{s})\Big)^{\!\gamma}\Big]\bigg).\nonumber
\end{align}
We now compare \eqref{express} to
\begin{equation}
k^{\varkappa}\mathcal{E}_{x}\varphi(W_{t})k^{(2+\varkappa-d)\gamma}%
\varrho^{\gamma}\int\!\!dy\Big(\int_{0}^{t}\!{d}s\;\vartheta_{1}%
(kW_{s}-y)\,S_{t-s}\varphi\,(W_{s})\Big)^{\!\gamma}. \label{exp1}%
\end{equation}
Clearly,
\begin{equation}
x-x^{2}\ \leq\ 1-\mathrm{e}^{-x}\ \leq\ x\quad\text{for }\,x\geq0.
\end{equation}
By the second inequality, the term (\ref{exp1}) is always an upper bound for
(\ref{express}). On the other hand, by the first inequality and
Lemma~\ref{la2}, the difference is bounded from above by a constant multiple
of
\begin{equation}%
\begin{array}
[c]{l}%
\displaystyle
k^{\varkappa}\mathcal{E}_{x}\varphi(W_{t})k^{2(2+\varkappa-d)\gamma}%
\Big[\int\!\!dy\Big(\int_{0}^{t}\!ds\ \vartheta_{1}(kW_{s}-y)\,S_{t-s}%
\varphi\,(W_{s})\Big)^{\!\gamma}\Big]^{2}\vspace{6pt}\\
\leq\ c_{\ref{la2}}k^{\varkappa+2(2+\varkappa-d)\gamma+4-4\gamma+\delta}%
\,\phi_{\gamma\lambda_{\ref{L.bo2}}}(x).
\end{array}
\end{equation}
Note that the exponent is negative iff $d\gamma>2+\delta\lbrack1+\gamma]$,
hence choosing $\delta>0$ sufficiently small justifies the approximation of
\eqref{express} by~\eqref{exp1}.

Recall that $\,\tau=\tau_{1/k}^{y}[W]$\thinspace\ denotes the \emph{first
hitting time}\/ of the ball $\,B(y,1/k)$\thinspace\ by our Brownian motion
$\,W$\thinspace\ started in $\,x.$ Now note that~\eqref{exp1} equals
\begin{equation}
k^{d-2}\varrho^{\gamma}\int\!\!dy\,\ \mathcal{E}_{x}\mathsf{1}_{\tau\leq
t}\,\mathcal{E}_{W_{\tau}}\,\varphi(\tilde{W}_{t-\tau})\Big(k^{2}\int
_{0}^{t-\tau}\!\!\!\!\!\!{d}s\;\vartheta_{1/k}(\tilde{W}_{s}-y)\,S_{t-\tau
-s}\varphi\,(\tilde{W}_{s})\Big)^{\!\gamma}, \label{92}%
\end{equation}
where the strong Markov property was used and the value for $\varkappa$ was
plugged in. By Lemma~\ref{highs} we may choose (and henceforth fix) a value
$M>1$ such that contributions to the innermost integral coming from
$s>M/k^{2}$, can be bounded by $\varepsilon\phi_{\lambda_{\ref{L.bo2}}%
}^{\gamma}(x)$, and additionally that
\begin{equation}
\mathcal{E}_{\iota}\Big(\int_{0}^{\infty}ds\,\vartheta_{1}(\tilde{W}%
_{s})\Big)^{\!\gamma}-\mathcal{E}_{\iota}\Big(\int_{0}^{M}ds\,\vartheta
_{1}(\tilde{W}_{s})\Big)^{\!\gamma}\,<\ \varepsilon. \label{Mchoose}%
\end{equation}
Moreover, by Lemma~\ref{la5}, if $k\geq k_{\ref{la5}}$, the contribution to
(\ref{92}) coming from $t-\delta\leq\tau\leq t$ can be made smaller than
$\varepsilon\phi_{\gamma\lambda_{\ref{L.bo2}}}(x)$ by choice of $\delta>0$,
which we also assume fixed from now on.

We let
\begin{equation}
k_{(\ref{k5})}\,:=\,\sqrt{M/\delta} \label{k5}%
\end{equation}
and note that $t-\tau\geq M/k^{2}$ whenever $t-\delta\geq\tau$ and $k\geq
k_{(\ref{k5})}$. Now let $\,k_{\ref{L.exp.conv}}\ :=\ k_{\ref{la5}}\vee
k_{\ref{replace1}}\vee k_{\ref{replace2}}\vee k_{\ref{la4}}\vee k_{(\ref{k5}%
)}.$\thinspace\ It remains to show that
\begin{align*}
\bigg|  &  k^{d-2}\varrho^{\gamma}\int\!\!dy\,\ \mathcal{E}_{x}\mathsf{1}%
_{\tau\leq t}\,\mathcal{E}_{W_{\tau}}\varphi(\tilde{W}_{t-\tau})\Big(k^{2}%
\int_{0}^{M/k^{2}}ds\,\vartheta_{1/k}(\tilde{W}_{s}-y)S_{t-\tau-s}%
\varphi\,(\tilde{W}_{s})\Big)^{\!\gamma}\\
&  \quad\quad-\ \overline{c}\int_{0}^{t}\!{d}r\;S_{r}(S_{t-r}\varphi
)^{1+\gamma}(x)\bigg|\ <\ \varepsilon\phi_{\gamma\lambda_{\ref{L.bo2}}%
}(x)\quad\text{ for }\,k\geq k_{\ref{L.exp.conv}\,},\ \,x\in\mathbb{R}^{d}.
\end{align*}
This will be done in three steps by the triangle inequality. The steps are
prepared in Lemmas~\ref{replace1} to \ref{la4}.

In the \emph{first} step note that by Lemma~\ref{replace1} we have, for all
$k\geq k_{\ref{L.exp.conv}}$ and $x\in\mathbb{R}^{d}$,
\begin{align}
k^{d-2}  &  \varrho^{\gamma}\int\!\!dy\,\ \mathcal{E}_{x}\mathsf{1}_{\tau\leq
t}\,\mathcal{E}_{W_{\tau}}\varphi(\tilde{W}_{t-\tau})\bigg|\Big(k^{2}\int
_{0}^{M/k^{2}}ds\,\vartheta_{1/k}(\tilde{W}_{s}-y)S_{t-\tau-s}\varphi
\,(\tilde{W}_{s})\Big)^{\!\gamma}\nonumber\\
&  \quad\quad-\Big(k^{2}\int_{0}^{M/k^{2}}ds\,\vartheta_{1/k}(\tilde{W}%
_{s}-y)S_{t-\tau}\varphi\,(y)\Big)^{\!\gamma}\bigg|\ \leq\ \varepsilon
\phi_{\gamma\lambda_{\ref{L.bo2}}}(x). \label{Arb1}%
\end{align}
We may therefore continue, using the Markov property,
\begin{align}
k^{d-2}  &  \varrho^{\gamma}\int\!\!dy\,\ \mathcal{E}_{x}\mathsf{1}_{\tau\leq
t}\,\mathcal{E}_{W_{\tau}}\varphi(\tilde{W}_{t-\tau})\Big(k^{2}\int
_{0}^{M/k^{2}}ds\,\vartheta_{1/k}(\tilde{W}_{s}-y)S_{t-\tau}\varphi
\,(y)\Big)^{\!\gamma}\nonumber\\
&  =\ k^{d-2}\varrho^{\gamma}\int\!\!dy\,\ \mathcal{E}_{x}\mathsf{1}_{\tau\leq
t}\,\!\left(  _{\!_{\!_{\,}}}S_{t-\tau}\varphi\,(y)\right)  ^{\!\gamma
}\mathcal{E}_{W_{\tau}}\Big(k^{2}\int_{0}^{M/k^{2}}ds\,\vartheta_{1/k}%
(\tilde{W}_{s}-y)\Big)^{\!\gamma}\\
&  \qquad\qquad\qquad\qquad\qquad\qquad\qquad\qquad\quad\quad\quad
\ \ \times\mathcal{E}_{\tilde{W}_{M/k^{2}}}\varphi(\tilde{\tilde{W}}%
_{t-\tau-M/k^{2}})\nonumber
\end{align}

As a \emph{second} step, by Lemma~\ref{replace2} we have, for all $k\geq
k_{\ref{L.exp.conv}}$ and $x\in\mathbb{R}^{d}$,
\begin{align*}
&  k^{d-2}\varrho^{\gamma}\int\!\!dy\,\ \mathcal{E}_{x}\mathsf{1}_{\tau\leq
t}\,\!\left(  _{\!_{\!_{\,}}}S_{t-\tau}\varphi\,(y)\right)  ^{\!\gamma}\\
&  \times\bigg|\mathcal{E}_{W_{\tau}}\Big(k^{2}\int_{0}^{M/k^{2}}%
ds\,\vartheta_{1/k}(\tilde{W}_{s}-y)\Big)^{\!\gamma}\mathcal{E}_{\tilde
{W}_{M/k^{2}}}\varphi(\tilde{\tilde{W}}_{t-\tau-M/k^{2}})-c_{\ref{replace2}%
}S_{t-\tau}\varphi\,(y)\bigg|\\[4pt]
&  \leq\ \varepsilon\,\phi_{\gamma\lambda_{\ref{L.bo2}}}(x).
\end{align*}
By \eqref{Mchoose} we have, using $\,c_{\eqref{Mchoose}}:=\mathcal{E}_{\iota
}\bigl(\int_{0}^{\infty}ds\,\vartheta_{1}(\tilde{W}_{s})\bigr)^{\!\gamma},$
\[
k^{d-2}\varrho^{\gamma}\!\int\!\!dy\,\mathcal{E}_{x}\mathsf{1}_{\tau\leq
t}\left(  _{\!_{\!_{\,}}}S_{t-\tau}\varphi\,(y)\right)  ^{\!\gamma
}\Big|c_{\ref{replace2}}S_{t-\tau}\varphi\,(y)-c_{\eqref{Mchoose}}S_{t-\tau
}\varphi\,(y)\Big|\,<\ \varepsilon\phi_{\gamma\lambda_{\ref{L.bo2}}}(x).
\]

In the \emph{third} step we recall that, by Lemma~\ref{la4}, for all $k\geq
k_{\ref{replace1}}$, and $x\in\mathbb{R}^{d}$,
\[
\bigg|\varrho^{\gamma}c_{\eqref{Mchoose}}k^{d-2}\!\!\int\!\!dy\!\!\int_{0}%
^{t}\!\!\!\pi_{x}(dr)\big(S_{t-r}\varphi\,(y)\big)^{\!1+\gamma}-c_{(\ref{ba1}%
)}\varrho^{\gamma}\!\!\!\int_{0}^{t}\!\!S_{r}(S_{t-r}\varphi)^{1+\gamma
}(x)\bigg|<\varepsilon\phi_{\gamma\lambda_{\ref{L.bo2}}}(x),
\]
and this completes the proof of Proposition~\ref{L.exp.conv}.\hfill$\square$

\subsection{Convergence of variances\label{varlin}}

In this section we establish that the variances with respect to the medium for
the solutions of the linearized integral equation vanish asymptotically.

\begin{proposition}
[\textbf{Convergence of variances}]\label{L.var}For every $\,\mu\in
\mathcal{M}_{\mathrm{tem}}\,\ $satisfying the assumption in Theorem\/
\emph{\ref{T.fluct},} for $\,\varphi\in\mathcal{C}_{\exp}^{+}$ and $\,t>0$,
\[
\lim_{k\rightarrow\infty}\mathsf{V}\mathrm{ar}\,k^{\varkappa}\!\!\int
\!\!\mu(dx)\,\mathcal{E}_{x}\varphi(W_{t})\exp\!\Big [\!-k^{2-d+\varkappa
}\varrho\!\int_{0}^{t}\!{d}s\;\Gamma^{1}(kW_{s})\,S_{t-s}\varphi
\,(W_{s})\Big]\ =\ 0.
\]

\end{proposition}

The remainder of this section is devoted to the proof of this proposition.
Recalling the definition (\ref{not.Gamma.eps}) of $\,\Gamma^{1},$ the variance
expression in Proposition~\ref{L.var} equals%
\begin{align}
&  k^{2\varkappa}\int\!\!\mu(dx)\!\int\!\!\mu(dy)\,\mathcal{E}_{x}%
\!\otimes\!\mathcal{E}_{y\,}\varphi(W_{t}^{1})\varphi(W_{t}^{2})\label{37}\\
&  \times\Bigg(\mathsf{E}\exp\!\bigg[-\int\!\Gamma({d}z)\,k^{2-d+\varkappa
}\varrho\sum_{i=1,2}\,\int_{0}^{t}\!{d}s\;\vartheta_{1}(kW_{s}^{i}%
-z)\,S_{t-s}\varphi\,(W_{s}^{i})\bigg ]\nonumber\\
&  \qquad-\;\prod_{i=1,2}\,\mathsf{E}\exp\!\Big[-\int\!\Gamma({d}%
z)\,k^{2-d+\varkappa}\varrho\int_{0}^{t}\!{d}s\;\vartheta_{1}(kW_{s}%
^{i}-z)\,S_{t-s}\varphi\,(W_{s}^{i})\Big]\Bigg ),\nonumber
\end{align}
where $\,(W^{1},W^{2})$\thinspace\ is distributed according to $\,\mathcal{P}%
_{x}\!\otimes\!\mathcal{P}_{y\,}.$\thinspace\ Exploiting the Laplace
functional (\ref{logL.Gamma}) of $\,\Gamma,$ (\ref{37}) can be rewritten as%
\begin{align}
&  k^{2\varkappa}\int\!\!\mu(dx)\!\int\!\!\mu(dy)\,\mathcal{E}_{x}%
\!\otimes\!\mathcal{E}_{y\,}\varphi(W_{t}^{1})\varphi(W_{t}^{2}%
)\label{written}\\
&  \times\Bigg(\exp\!\bigg[-\int\!\!{d}z\ k^{(\varkappa-d)\gamma}%
\varrho^{\gamma}\Big(\sum_{i=1,2}\,k^{2}\!\int_{0}^{t}\!{d}s\;\vartheta
_{1}(kW_{s}^{i}-z)\,S_{t-s}\varphi\,(W_{s}^{i})\Big )^{\!\gamma}%
\bigg]\nonumber\\
&  \ \quad\ -\ \exp\!\bigg[\!-\!\int\!\!{d}z\ k^{(\varkappa-d)\gamma}%
\varrho^{\gamma}\!\sum_{i=1,2}\Big(k^{2}\!\int_{0}^{t}\!{d}s\;\vartheta
_{1}(kW_{s}^{i}-z)\,S_{t-s}\varphi\,(W_{s}^{i})\Big)^{\!\gamma}%
\bigg]\Bigg).\nonumber
\end{align}
Note that by the elementary inequality%
\begin{equation}
(a+b)^{\gamma}\;\leq\;a^{\gamma}+b^{\gamma}\quad\text{for}\,\ a,b\geq0,
\label{elementary}%
\end{equation}
the argument in the first exponential expression is not smaller than the
argument in the second one. Therefore we may apply the elementary inequality%
\begin{equation}
\mathrm{e}^{-a}-\mathrm{e}^{-b}\ \leq\ b-a\quad\text{for}\,\ 0\leq a\leq b,
\label{alsoelementary}%
\end{equation}
and a $z$ substitution to get for the non-negative total expression in
(\ref{written}) the upper bound (we may drop from now on the factor
$\,\varrho^{\gamma})$
\begin{align}
&  k^{2\varkappa+(\varkappa-d)\gamma+d}\int\!\!{d}z\int\!\!\mu(dx)\!\int
\!\!\mu(dy)\,\mathcal{E}_{x}\!\otimes\!\mathcal{E}_{y\,}\varphi(W_{t}%
^{1})\varphi(W_{t}^{2})\label{121}\\
&  \qquad\times\bigg[\sum_{i=1,2}\Big(k^{2}\!\int_{0}^{t}\!{d}s\;\vartheta
_{1/k}(W_{s}^{i}-z)\,S_{t-s}\varphi\,(W_{s}^{i})\Big)^{\!\gamma}\nonumber\\
&  \qquad\qquad\qquad-\;\Big(\sum_{i=1,2}k^{2}\!\int_{0}^{t}\!{d}%
s\;\vartheta_{1/k}(W_{s}^{i}-z)\,S_{t-s}\varphi\,(W_{s}^{i})\Big)^{\!\gamma
}\bigg].\nonumber
\end{align}
It remains to show that (\ref{121}) converges to zero as $k\uparrow\infty.$
The proof rests solely on the fact that the square bracket expression vanishes
if \emph{one} of the motions does not hit the ball $\,B(z,1/k)$. For
simplification, write now $\tau\lbrack W^{i}]$ for the first hitting time
$\tau_{1/k}^{z}[W^{i}]$ of $B(z,1/k)$ by the Brownian motion $\,W^{i}%
$.\thinspace\ Hence, we get the bound
\begin{align}
&  k^{2\varkappa+(\varkappa-d)\gamma+d}\int\!\!{d}z\int\!\!\mu(dx)\!\int
\!\!\mu(dy)\,\mathcal{E}_{x}\!\otimes\!\mathcal{E}_{y}\,\mathsf{1}_{\left\{
\tau\lbrack W^{1}]\leq t\right\}  }\,\mathsf{1}_{\left\{  \tau\lbrack
W^{2}]\leq t\right\}  }\\
&  \qquad\times\varphi(W_{t}^{1})\,\varphi(W_{t}^{2})\sum_{i=1,2}%
\Big(k^{2}\!\int_{0}^{t}\!{d}s\;\vartheta_{1/k}(W_{s}^{i}-z)\,S_{t-s}%
\varphi\,(W_{s}^{i})\Big)^{\!\gamma},\nonumber
\end{align}
where we dropped the subtracted term. Interchanging expectation and summation,
and using independence, we obtain
\begin{align}
&  k^{2\varkappa+(\varkappa-d)\gamma+d}\int\!\!{d}z\int\!\!\mu(dx)\!\int
\!\!\mu(dy)\sum_{i=1,2}\mathcal{E}_{x}\mathsf{1}_{\left\{  \tau\lbrack
W^{i}]\leq t\right\}  }\,\varphi(W_{t}^{i})\label{get}\\[1pt]
&  \qquad\times\Big(k^{2}\!\int_{0}^{t}\!{d}s\;\vartheta_{1/k}(W_{s}%
^{i}-z)\,S_{t-s}\varphi\,(W_{s}^{i})\Big)^{\!\gamma}\,\mathcal{E}%
_{y}\mathsf{1}_{\left\{  \tau\lbrack W^{j}]\leq t\right\}  }\,\varphi
(W_{t}^{j}),\nonumber
\end{align}
where $j=3-i$. Then we may bound (\ref{get}) by%
\begin{gather}
c_{\ref{L.bo2}}^{\gamma}k^{2\varkappa+(\varkappa-d)\gamma+d}\int
\!\!{d}z\;\tilde{\phi}^{\gamma}(z)\int\!\!\mu(dx)\!\int\!\!\mu(dy)\sum
_{i=1,2}\mathcal{E}_{x}\mathsf{1}_{\left\{  \tau\lbrack W^{i}]\leq t\right\}
}\,\varphi(W_{t}^{i})\label{getalso}\\[1pt]
\times\Big(k^{2}\!\int_{0}^{t}\!{d}s\;\vartheta_{1/k}(W_{s}^{i}%
-z)\Big)^{\!\gamma}\,\mathcal{E}_{y}\mathsf{1}_{\left\{  \tau\lbrack
W^{j}]\leq t\right\}  }\,\varphi(W_{t}^{j}).\nonumber
\end{gather}
By Lemma~\ref{LemmaC}, there are constants $\,c_{(\ref{Hoe1})}$\thinspace\ and
$\,c_{(\ref{Hoe2})}$\thinspace\ such that, for all $\,x,y\in\mathbb{R}^{d}$,
\begin{equation}
\mathcal{E}_{y}\mathsf{1}_{\left\{  \tau\lbrack W^{j}]\leq t\right\}
}\,\varphi(W_{t}^{j})\ \leq\ c_{(\ref{Hoe1})}k^{2-d}\phi_{\lambda
_{\ref{L.bo2}}}(z)\,\big[|z-y|^{2-d}+1\big], \label{Hoe1}%
\end{equation}
and
\begin{align}
\mathcal{E}_{x}  &  \mathsf{1}_{\left\{  \tau\lbrack W^{i}]\leq t\right\}
}\,\Big(k^{2}\!\int_{0}^{t}\!{d}s\;\vartheta_{1/k}(W_{s}^{i}-z)\Big)^{\!\gamma
}\varphi(W_{t}^{i})\label{Hoe2}\\
&  \leq\ c_{(\ref{Hoe2})}\,\phi_{\lambda_{\ref{L.bo2}}}(x)\big[1+|z-x|^{2-d}%
\big]\,k^{2-d}.\nonumber
\end{align}
Assume for the moment that $\,d\geq5.$\thinspace\ Then, by (\ref{all.d}) and
Lemma~\ref{L.sim.id}, for each $\,\lambda>0$\thinspace\ there is a constant
$\,c_{(\ref{da2})}=c_{(\ref{da2})}(\lambda),$\thinspace\ such that%
\begin{equation}
\int\!\!{d}z\ \phi_{\lambda}(z)\,\big[1+|z-x|^{2-d}\big]\,\big[1+|z-y|^{2-d}%
\big]\ \leq\ c_{(\ref{da2})}\,\big[1+|x-y|^{4-d}\big]. \label{da2}%
\end{equation}
If $\,d=3,$\thinspace\ the left hand side of (\ref{da2}) is even bounded in
$\,x,y.$\thinspace\ In fact by (\ref{all.d}) and Cauchy-Schwarz, it suffices
to consider the singularity $\,\int_{|z|\leq1}\!dz$\thinspace$|z|^{2(2-d)}%
<\infty.$\thinspace\ Finally, if $\,d=4,$\thinspace\ by (\ref{all.d}) and
Lemma~\ref{L.sim.est}, estimate (\ref{da2}) holds if $\,|x-y|^{4-d}$%
\thinspace\ is replaced by $\,\log^{+}\left(  |x-y|^{-1}\right)  .$%
\thinspace\ If we extend definition (\ref{en}) by setting $\,\mathrm{en}%
(x):\equiv1$\thinspace\ in the case $\,d=3,$\thinspace\ then we can combine
the last three steps to obtain that for each $\,\lambda>0$\thinspace\ there is
a constant $\,c_{(\ref{dge3})}=c_{(\ref{dge3})}(\lambda),$\thinspace\ so that
for all $\,d\geq3,$%
\begin{equation}
\int\!\!{d}z\ \phi_{\lambda}(z)\,\big[1+|z-x|^{2-d}\big]\,\big[1+|z-y|^{2-d}%
\big]\ \leq\ c_{(\ref{dge3})}\,\big[1+\mathrm{en}(x-y)\big]. \label{dge3}%
\end{equation}

Based on \eqref{Hoe1}, \eqref{Hoe2} and (\ref{dge3}), from (\ref{getalso}) we
get the upper bound
\begin{equation}
c_{(\ref{hier})}\,k^{2\varkappa+(\varkappa-d)\gamma+d}\,k^{4-2d}\int
\!\!\mu(dx)\,\phi_{\lambda_{\ref{L.bo2}}}(x)\int\!\!\mu(dy)\,\phi
_{\lambda_{\ref{L.bo2}}}(y)\,\big[1+\mathrm{en}(x-y)\big]. \label{hier}%
\end{equation}
By our condition on $\mu$, the latter integral is finite. Moreover,
$2\varkappa+(\varkappa-d)\gamma+d+4-2d<2-d+\varkappa.$\thinspace\ But the last
expression is negative, finishing the proof.\hfill$\square$

\subsection{Upper bound for finite-dimensional distributions\label{findim}}

We use an induction argument to extend the result from the convergence of
one-dimensional distributions to all finite dimensional distributions. Recall
that we have to show that, for any $\varphi_{1},\ldots,\varphi_{n}$ and
$0=t_{0}<t_{1}<\cdots<t_{n}$, in $\mathsf{P}$--probability,
\begin{equation}%
\begin{array}
[c]{l}%
\displaystyle
\limsup_{k\rightarrow\infty}\mathbb{E}_{\mu_{k}}\!\exp\!\Big[\sum_{i=1}%
^{n}k^{\varkappa}\big\langle X_{t_{i}}^{k}-S_{t_{i}}\mu,\,-\varphi
_{i}\big\rangle\Big ]\vspace{6pt}\\
\qquad%
\displaystyle
\leq\ \exp\!\Bigg[\overline{c}\,%
\bigg\langle
\mu,\sum_{i=1}^{n}\int_{t_{i-1}}^{t_{i}}\!dr\ S_{r}\Big(\big(\sum_{j=i}%
^{n}S_{t_{j}-r}\varphi_{j}\big)^{1+\gamma}\Big)\!%
\bigg\rangle
\Bigg].
\end{array}
\end{equation}
The case $n=1$ was shown in the previous paragraphs, so we may assume that it
holds for $n-1$ and show that it also holds for $n$. By conditioning on
$\{X^{k}(t)\colon t\leq t_{n-1}\}$ and applying the transition functional we
get
\begin{align}
\mathbb{E}_{\mu_{k}}  &  \exp\!\Big[\sum_{i=1}^{n}k^{\varkappa}%
\big\langle X_{t_{i}}^{k}-S_{t_{i}}\mu,\,-\varphi_{i}%
\big\rangle\Big ]\nonumber\\
&  =\ \mathbb{E}_{\mu_{k}}\!\exp\!\bigg[\sum_{i=1}^{n-1}k^{\varkappa
}\big\langle X_{t_{i}}^{k}-S_{t_{i}}\mu,\,-\varphi_{i}\big\rangle\\
&  \qquad\qquad\qquad+\ k^{\varkappa}\big\langle S_{t_{n-1}}\mu,S_{t_{n}%
-t_{n-1}}\varphi_{n}\big\rangle-\big\langle X_{t_{n-1}}^{k},u_{k}%
(t_{n}-t_{n-1})\big\rangle\bigg],\nonumber
\end{align}
where $u_{k}$ is the solution of (\ref{new}) with $\varphi$ replaced by
$\,\varphi_{n\,}.$\thinspace\ Separating the non-random terms yields
\begin{align}
=\  &  \exp\!\Big[\big\langle S_{t_{n-1}}\mu,k^{\varkappa}S_{t_{n}-t_{n-1}%
}\varphi_{n}-u_{k}(t_{n}-t_{n-1})\big\rangle\Big ]\,\nonumber\\
&  \times\mathbb{E}_{\mu_{k}}\!\exp\!\bigg[\sum_{i=1}^{n-2}k^{\varkappa
}\big\langle X_{t_{i}}^{k}-S_{t_{i}}\mu,\,-\varphi_{i}\big\rangle\\
&  \qquad\quad\quad\qquad+\ k^{\varkappa}%
\Big\langle
X_{t_{n-1}}^{k}-S_{t_{n-1}}\mu,\,-\varphi_{n-1}-k^{-\varkappa}u_{k}%
(t_{n}-t_{n-1})%
\Big\rangle
\!\bigg].\nonumber
\end{align}
By Theorem~\ref{T.fluct} for $n=1$ with starting measure $\,S_{t_{n-1}}\mu
,$\thinspace\ in $\mathsf{P}$--probability,%
\begin{equation}%
\begin{array}
[c]{l}%
\displaystyle
\limsup_{k\uparrow\infty}\exp\!\Big[\big\langle S_{t_{n-1}}\mu,k^{\varkappa
}S_{t_{n}-t_{n-1}}\varphi_{n}-u_{k}(t_{n}-t_{n-1})\big\rangle\Big ]\\%
\displaystyle
\qquad\leq\ \exp\!\bigg[\overline{c}\,%
\Big\langle
S_{t_{n-1}}\mu,\int_{0}^{t_{n}-t_{n-1}}\!dr\ S_{r}(S_{t_{n}-t_{n-1}-r}%
\varphi_{n})^{1+\gamma}%
\Big\rangle
\bigg]\vspace{2pt}\\%
\displaystyle
\qquad=\ \exp\!\bigg[\overline{c}\,%
\Big\langle
\mu,\int_{t_{n-1}}^{t_{n}}\!dr\ S_{r}(S_{t_{n}-r}\varphi_{n})^{1+\gamma}%
\Big\rangle
\bigg]\vspace{2pt}.
\end{array}
\label{erschterterm}%
\end{equation}
The remaining expectation can be written as
\begin{align}
\mathbb{E}_{\mu_{k}}\!\exp\!\bigg[  &  \sum_{i=1}^{n-2}k^{\varkappa
}\big\langle X_{t_{i}}^{k}-S_{t_{i}}\mu,\,-\varphi_{i}\big\rangle\nonumber\\
&  +\ k^{\varkappa}\big\langle X_{t_{n-1}}^{k}-S_{t_{n-1}}\mu,\,-\varphi
_{n-1}-S_{t_{n}-t_{n-1}}\varphi_{n}\big\rangle\label{remaining}\\
&  +\ k^{\varkappa}%
\Big\langle
X_{t_{n-1}}^{k}-S_{t_{n-1}}\mu,\ S_{t_{n}-t_{n-1}}\varphi_{n}-k^{-\varkappa
}u_{k}(t_{n}-t_{n-1})%
\Big\rangle
\bigg].\nonumber
\end{align}
To dominate this term observe that, by the induction assumption, in
{$\mathsf{P}$--probability},
\begin{align}
\limsup_{k\uparrow\infty}  &  \,\mathbb{E}_{\mu_{k}}\!\exp\!\Big[\sum
_{i=1}^{n-2}k^{\varkappa}\big\langle X_{t_{i}}^{k}-S_{t_{i}}\mu,\,-\varphi
_{i}\big\rangle\nonumber\\
\qquad &  \qquad\qquad+\ k^{\varkappa}\big\langle X_{t_{n-1}}^{k}-S_{t_{n-1}%
}\mu,\,-\varphi_{n-1}-S_{t_{n}-t_{n-1}}\varphi_{n}\big\rangle\Big ]\nonumber\\
&  \leq\ \exp\!\bigg[\overline{c}\,\Big\langle\mu,\sum_{i=1}^{n-2}%
\int_{t_{i-1}}^{t_{i}}\!dr\ S_{r}\Big(\big(\sum_{j=i}^{n-1}S_{t_{j}-r}%
\varphi_{j}\big)^{1+\gamma}\Big)\label{ind.ass}\\
&  \qquad\qquad+\int_{t_{n-2}}^{t_{n-1}}\!dr\ S_{r}\Big(\big(S_{t_{n-1}%
-r}\varphi_{n-1}+S_{t_{n-1}-r}S_{t_{n}-t_{n-1}}\varphi_{n}\big)^{1+\gamma
}\Big)\Big\rangle\!\bigg].\nonumber
\end{align}
We show below that in {$\mathsf{P}$--probability} the following convergence in
law holds:%
\begin{equation}
\exp\!\bigg[k^{\varkappa}\left\langle _{\!_{\!_{\,_{{}}}}}X_{t_{n-1}}%
^{k}-S_{t_{n-1}}\mu,S_{t_{n}-t_{n-1}}\varphi_{n}-k^{-\varkappa}\,u_{k}%
(t_{n}-t_{n-1})\right\rangle \!\bigg]\ \underset{k\uparrow\infty
}{\Longrightarrow}\ 1. \label{to.show}%
\end{equation}
Observe that $\,\limsup_{m\uparrow\infty}\xi_{m}\leq a$\thinspace\ in
probability for some $a,$ and $\zeta_{m}\Rightarrow1$ in law implies
$\limsup_{m\uparrow\infty}\xi_{m}\zeta_{m}\leq a$\thinspace\ in probability.
Hence \eqref{ind.ass} and \eqref{to.show} together imply that
\eqref{remaining} is asymptotically bounded from above by
\begin{equation}%
\begin{array}
[c]{l}%
\displaystyle
\exp\!\bigg[\overline{c}\,\Big\langle\mu,\sum_{i=1}^{n-2}\int_{t_{i-1}}%
^{t_{i}}\!dr\ S_{r}\Big(\big(\sum_{j=i}^{n-1}S_{t_{j}-r}\varphi_{j}%
\big)^{1+\gamma}\Big)\vspace{6pt}\\
\qquad\qquad\qquad%
\displaystyle
+\ \int_{t_{n-2}}^{t_{n-1}}\!dr\ S_{r}\Big(\big(S_{t_{n-1}-r}\varphi
_{n-1}+S_{t_{n}-r}\varphi_{n}\big)^{1+\gamma}\Big)\Big\rangle\bigg].
\end{array}
\label{zwoterterm}%
\end{equation}
Putting together \eqref{erschterterm} and \eqref{zwoterterm} yields the
claimed statement subject to the proof of \eqref{to.show}.

To prove \eqref{to.show} it suffices to show that, for any $a\geq1,$
\begin{equation}
\mathbb{E}_{\mu_{k}}\!\exp\!\bigg[ak^{\varkappa}\left\langle _{\!_{\!_{\,_{{}%
}}}}X_{t_{n-1}}^{k}-S_{t_{n-1}}\mu,\ -S_{t_{n}-t_{n-1}}\varphi_{n}%
+k^{-\varkappa}\,u_{k}(t_{n}-t_{n-1})\right\rangle \!\bigg] \label{to.show2}%
\end{equation}
converges in $\mathsf{P}$--probability to $1$. Using the Feynman-Kac
representation (\ref{FK}), the expectation in (\ref{to.show2}) equals
\begin{equation}%
\begin{array}
[c]{l}%
\exp\!%
\bigg\langle
\mu,ak^{\varkappa}\mathcal{E}_{x}\big(S_{t_{n}-t_{n-1}}\varphi_{n}(W_{t_{n-1}%
})-k^{-\varkappa}\,u_{k}(t_{n}-t_{n-1},W_{t_{n-1}})\big)\vspace{6pt}\\%
\displaystyle
\qquad\quad\times\Big(1-\exp\!\Big[-k^{2-d}\int_{0}^{t_{n-1}}\!dr\ \Gamma
^{1}(kW_{r})\,U_{k}(t_{n-1}-r,W_{r})\Big ]\Big)\!%
\bigg\rangle
,
\end{array}
\end{equation}
where $U_{k}$ is the solution of (\ref{new}) with $\varphi$ replaced by
$\,a\bigl(S_{t_{n}-t_{n-1}}\varphi_{n}-k^{-\varkappa}\,u_{k}(t_{n}%
-t_{n-1})\bigr).$\thinspace\ It therefore suffices to show that
\begin{equation}%
\begin{array}
[c]{l}%
\bigg\langle
\mu,\,ak^{\varkappa}\mathcal{E}_{x}\big(S_{t_{n}-t_{n-1}}\varphi
_{n}(W_{t_{n-1}})-k^{-\varkappa}\,u_{k}(t_{n}-t_{n-1},W_{t_{n-1}}%
)\big)\vspace{6pt}\\
\quad%
\displaystyle
\times\Big(1-\exp\!\Big[-k^{2-d}\int_{0}^{t_{n-1}}\!dr\ \Gamma^{1}%
(kW_{r})\,U_{k}(t_{n-1}-r,W_{r})\Big ]\Big)\!%
\bigg\rangle
\end{array}
\end{equation}
converges in $L^{1}(\mathsf{P})$ to zero. As this term is non-negative and as
\begin{equation}
U_{k}(t_{n-1}-r)\ \leq\ ak^{\varkappa}S_{t_{n}-r}\Big(S_{t_{n}-t_{n-1}}%
\varphi_{n}-k^{-\varkappa}\,u_{k}(t_{n}-t_{n-1})\Big)\ \leq\ ak^{\varkappa
}S_{t_{n}-r}\varphi_{n}\,,
\end{equation}
it finally suffices to show that
\begin{equation}%
\begin{array}
[c]{l}%
\mathsf{E}\,%
\bigg\langle
\mu,\,ak^{\varkappa}\mathcal{E}_{x}\big(S_{t_{n}-t_{n-1}}\varphi
_{n}(W_{t_{n-1}})-k^{-\varkappa}\,u_{k}(t_{n}-t_{n-1},W_{t_{n-1}}%
)\big)\vspace{6pt}\\%
\displaystyle
\qquad\times\ \Big(1-\exp\!\Big[-ak^{2-d+\varkappa}\int_{0}^{t_{n-1}%
}\!dr\ \Gamma^{1}(kW_{r})S_{t_{n}-r}\varphi_{n}(W_{r})\Big ]\Big)\!%
\bigg\rangle
\end{array}
\label{to.show3}%
\end{equation}
converges to zero. The first factor in the expectation can be expressed using
the Feynman-Kac representation~\eqref{FK} of $\,u_{k\,},$\thinspace\ which
gives
\begin{align}
a\int\!\!  &  \mu(dx)\,\mathsf{E}\mathcal{E}_{x}\biggl(k^{\varkappa
}\mathcal{E}_{W_{t_{n-1}}}\varphi_{n}(\tilde{W}_{t_{n}-t_{n-1}})\nonumber\\
&  \times\Big(1-\exp\!\Big[-k^{2-d}\int_{0}^{t_{n}-t_{n-1}}\!dr\ \Gamma
^{1}(k\tilde{W}_{r})u_{k}(t_{n}-t_{n-1}-r,\tilde{W}_{r}%
)\Big ]\Big)\label{exp01}\\
&  \times\Big(1-\exp\!\Big[-ak^{2-d+\varkappa}\int_{0}^{t_{n-1}}%
\!dr\ \Gamma^{1}(kW_{r})S_{t_{n}-r}\varphi_{n}(W_{r}%
)\Big ]\Big)\biggr),\nonumber
\end{align}
which again is dominated by
\begin{align}
&  a\int\!\!\mu(dx)\,\mathsf{E}\mathcal{E}_{x}\biggl(k^{\varkappa}%
\mathcal{E}_{W_{t_{n-1}}}\varphi_{n}(\tilde{W}_{t_{n}-t_{n-1}})\nonumber\\
&  \times\Big(1-\exp\!\Big[-ak^{2-d+\varkappa}\int_{0}^{t_{n}-t_{n-1}%
}\!dr\ \Gamma^{1}(k\tilde{W}_{r})S_{t_{n}-t_{n-1}-r}\varphi_{n}(\tilde{W}%
_{r})\Big ]\Big)\label{exp02}\\
&  \times\Big(1-\exp\!\Big[-ak^{2-d+\varkappa}\int_{0}^{t_{n-1}}%
\!dr\ \Gamma^{1}(kW_{r})S_{t_{n}-r}\varphi_{n}(W_{r}%
)\Big ]\Big)\biggr).\nonumber
\end{align}
We can now multiply the factors out and obtain
\begin{subequations}
\label{line}%
\begin{align}
a\int\!\!\mu(dx)\,{}  &  {}\mathsf{E}\mathcal{E}_{x}k^{\varkappa}{}%
\mathcal{E}_{W_{t_{n-1}}}\varphi_{n}(\tilde{W}_{t_{n}-t_{n-1}})\nonumber\\
\times\,  &  \bigg[\Big(1-\exp\!\Big[-ak^{2-d+\varkappa}\int_{0}%
^{t_{n}-t_{n-1}}\!dr\ \Gamma^{1}(k\tilde{W}_{r})S_{t_{n}-t_{n-1}-r}\varphi
_{n}(\tilde{W}_{r})\Big ]\Big)\label{line1}\\
&  \qquad+\Big(1-\exp\!\Big[-ak^{2-d+\varkappa}\int_{0}^{t_{n-1}}%
\!dr\ \Gamma^{1}(kW_{r})S_{t_{n}-r}\varphi_{n}(W_{r})\Big ]\Big) \label{line2}%
\\
&  \qquad-\Big(1-\exp\!\Big[-ak^{2-d+\varkappa}\int_{0}^{t_{n}-t_{n-1}%
}\!dr\ \Gamma^{1}(k\tilde{W}_{r})S_{t_{n}-t_{n-1}-r}\varphi_{n}(\tilde{W}%
_{r})\label{line3}\\
&  \qquad\qquad\qquad\qquad-ak^{2-d+\varkappa}\int_{0}^{t_{n-1}}%
\!dr\ \Gamma^{1}(kW_{r})S_{t_{n}-r}\varphi_{n}(W_{r}%
)\Big ]\Big)\bigg].\nonumber
\end{align}
We can now determine the limit in each of the three summands \eqref{line1} to
\eqref{line3} separately. For the first one we obtain from
Proposition~\ref{L.exp.conv}, as $k\uparrow\infty$,
\end{subequations}
\begin{equation}%
\begin{array}
[c]{l}%
a%
\displaystyle
\int\mu(dx)\,{}{}\mathsf{E}\mathcal{E}_{x}k^{\varkappa}{}\mathcal{E}%
_{W_{t_{n-1}}}\varphi_{n}(\tilde{W}_{t_{n}-t_{n-1}})\vspace{6pt}\\
\quad%
\displaystyle
\times\,\Big(1-\exp\!\Big[-ak^{2-d+\varkappa}\int_{0}^{t_{n}-t_{n-1}%
}\!dr\ \Gamma^{1}(k\tilde{W}_{r})S_{t_{n}-t_{n-1}-r}\varphi_{n}(\tilde{W}%
_{r})\Big ]\Big)\vspace{6pt}\\%
\displaystyle
\;\underset{k\uparrow\infty}{\Longrightarrow}\;a\int\mu(dx)\,\mathcal{E}%
_{x}\,\overline{c}\int_{0}^{t_{n}-t_{n-1}}\!dr\ S_{r}\big(aS_{t_{n}-t_{n-1}%
-r}\varphi_{n}\big)^{1+\gamma}(W_{t_{n-1}})\vspace{6pt}\\
\quad%
\displaystyle
=\ \overline{c}\,a^{2+\gamma}\Big\langle\mu,\int_{t_{n-1}}^{t_{n}}%
\!dr\ S_{r}\big(S_{t_{n}-r}\varphi_{n}\big)^{1+\gamma}\Big\rangle.
\end{array}
\label{lim1}%
\end{equation}
Similarly, the second one, \eqref{line2}, converges by
Proposition~\ref{L.exp.conv}, as $k\uparrow\infty$,
\begin{equation}%
\begin{array}
[c]{l}%
\displaystyle
a\int\!\!\mu(dx)\,{}{}\mathsf{E}\mathcal{E}_{x}k^{\varkappa}{}\mathcal{E}%
_{W_{t_{n-1}}}\varphi_{n}(\tilde{W}_{t_{n}-t_{n-1}})\vspace{6pt}\\
\quad%
\displaystyle
\times\,\Big(1-\exp\!\Big[-ak^{2-d+\varkappa}\int_{0}^{t_{n-1}}\!dr\ \Gamma
^{1}(kW_{r})S_{t_{n-1}-r}\big(S_{t_{n}-t_{n-1}}\varphi_{n}\big)(W_{r}%
)\Big ]\Big)\vspace{6pt}\\%
\displaystyle
\;\underset{k\uparrow\infty}{\Longrightarrow}\;a\int\mu(dx)\,\overline{c}%
\int_{0}^{t_{n-1}}\!dr\ S_{r}\big(S_{t_{n-1}-r}(aS_{t_{n}-t_{n-1}}\varphi
_{n})\big)^{1+\gamma}(x)\Big\}\vspace{6pt}\\
\quad%
\displaystyle
=\ \overline{c}\,a^{2+\gamma}\Big\langle\mu,\int_{0}^{t_{n-1}}\!dr\ S_{r}%
\big(S_{t_{n}-r}\varphi_{n}\big)^{1+\gamma}\Big\rangle.
\end{array}
\label{lim2}%
\end{equation}
Finally, the last expression \eqref{line3} equals, using
Proposition~\ref{L.exp.conv} to take the limit as $k\uparrow\infty$,
\begin{align}
-a\int &  \mu(dx)\,{}\mathsf{E}\mathcal{E}_{x}k^{\varkappa}{}\varphi
_{n}(\tilde{W}_{t_{n}})\nonumber\\
&  \times\,\Big(1-\exp\!\Big[-ak^{2-d+\varkappa}\int_{0}^{t_{n}}%
\!dr\ \Gamma^{1}(kW_{r})S_{t_{n}-r}\varphi_{n}(W_{r})\Big ]\Big)\label{lim3}\\
\;\underset{k\uparrow\infty}{\Longrightarrow}\;  &  -\,\overline
{c}\,a^{2+\gamma}\Big\langle\mu,\int_{0}^{t_{n}}\!dr\ S_{r}\big(S_{t_{n}%
-r}\varphi_{n}\big)^{1+\gamma}\Big\rangle.\nonumber
\end{align}
Comparing the right hand sides of \eqref{lim1} to \eqref{lim3} shows that they
cancel completely, which proves \eqref{to.show3} and completes the argument.

\section{Lower bound: Proof of (\ref{fluct1'})\label{S.lower}}

\subsection{A heat equation with random inhomogeneity\label{SS.heat.equ}}

As motivated in Section~\ref{SS.heuristics}, we look at the mild solution
$m_{k}$ to the linear equation
\begin{equation}%
\begin{array}
[c]{c}%
\displaystyle
\frac{\partial}{\partial t}m_{k}(t,x)\;=\;\tfrac{1}{2}\Delta m_{k}%
(t,x)-k^{2-d}\,\varrho\,\Gamma^{1}(kx)\,w_{k}^{2}(t,x)\vspace{6pt}\\
\text{with initial condition }\,m_{k}(0,\,\cdot\,)=k^{\varkappa}\varphi.
\end{array}
\label{heat}%
\end{equation}
This is a heat equation with the time-dependent scaled random inhomogeneity\/
$-k^{2-d}\,\varrho\,\Gamma^{1}(kx)\,w_{k}^{2}(t,x).$\thinspace\ We study its
asymptotic fluctuation behaviour around the heat flow:

\begin{proposition}
[\textbf{Limiting fluctuations of }$m_{k}$\thinspace]\label{P.heat}\hfill
Under the assumptions of \newline Theorem~\emph{\ref{T.fluct},} if
$\,\varkappa=\varkappa_{\mathrm{c}\,},$\thinspace\ then for any $\,\varphi
\in\mathcal{C}_{\exp}^{+}$\thinspace\ and $\,t\geq0,$\thinspace\ in
$\mathsf{P}$--probability,%
\begin{equation}
\liminf_{k\uparrow\infty}\left\langle _{\!_{\!_{\,}}}\mu,k^{\varkappa}%
S_{t}\varphi-m_{k}(t,\cdot)\right\rangle \;\geq\;\underline{c}\,%
\Big\langle
\mu,\int_{0}^{t}\!dr\ S_{r}\!\left(  _{\!_{\!_{\,}}}(S_{t-r}\varphi
)^{1+\gamma}\right)  \!%
\Big\rangle
, \label{heat.conv}%
\end{equation}
where the constant $\underline{c}=\underline{c}(\gamma,\varrho)$\thinspace\ is
given by
\begin{equation}
\underline{c}\ :=\ \gamma\,\varrho^{\gamma}\,\frac{2\,\pi^{d/2}}%
{d\,G(d/2)}\,\mathcal{E}_{0}\otimes\mathcal{E}_{0}\Big[\int_{0}^{\infty
}\!dr\ \vartheta_{2}(W_{r}^{1})+\int_{0}^{\infty}\!dr\ \vartheta_{2}(W_{r}%
^{2})\Big ]^{\gamma-1}. \label{c.u}%
\end{equation}

\end{proposition}

To see how the case $n=1$ of (\ref{fluct1'}) follows from
Proposition~\ref{P.heat}, we fix a sample~$\Gamma.$\thinspace\ Recall that
\begin{equation}
\log\,\mathbb{E}_{\mu_{k}}\!\exp\!\left[  k^{\varkappa}\big(\langle
X_{t\,}^{k},-\varphi\rangle-\langle S_{t}\mu,-\varphi\rangle
\big)_{\!_{\!_{\,_{{}}}}}\right]  =\;\left\langle \mu,\,k^{\varkappa}%
S_{t}\varphi-u_{k}(t,\,\cdot\,)_{\!_{\!_{\,}}}\right\rangle \!,
\end{equation}
where $\,u_{k}$\thinspace\ solves
\begin{equation}
k^{\varkappa}S_{t}\varphi\,(x)\,-\,u_{k}(t,x)\ =\ k^{2-d}\varrho\int_{0}%
^{t}\!{d}s\;S_{s}\bigl(\Gamma^{1}(k\,\cdot\,)\,u_{k}^{2}(t-s,\,\cdot
\,)\bigr)(x).
\end{equation}
As $u_{k}^{2}\geq w_{k\,}^{2},$\thinspace\ we obtain from (\ref{heat}),%
\begin{equation}
k^{\varkappa}S_{t}\varphi\,(x)-u_{k}(t,x)\ \geq\ k^{\varkappa}S_{t}%
\varphi\,(x)-m_{k}(t,x).
\end{equation}
Hence, the case $n=1$ of (\ref{fluct1'}) follows from Proposition~\ref{P.heat}.

Proposition~\ref{P.heat} is proved in two steps: In Section~\ref{SS.exp} we
show that the right hand side of (\ref{heat.conv}) is an asymptotic lower
bound of the expectations of the left hand side, and in Section~\ref{SS.var}
that the variances vanish asymptotically.

\subsection{Convergence of expectations\label{SS.exp}}

Fix again $t\geq0\,\ $and $\,\varphi\in\mathcal{C}_{\exp\,}^{+}.$

\begin{proposition}
[\textbf{Convergence of expectations}]\label{L.conv.exp}For $\underline{c}$ as
in \emph{(\ref{c.u}),}%
\begin{equation}
\liminf_{k\uparrow\infty}\mathsf{E}\left\langle _{\!_{\!_{\,}}}\mu
,k^{\varkappa}S_{t}\varphi-m_{k}(t,\cdot)\right\rangle \;\geq\;\underline{c}\,%
\Big\langle
\mu,\int_{0}^{t}\!dr\ S_{r}\!\left(  _{\!_{\!_{\,}}}(S_{t-r}\varphi
)^{1+\gamma}\right)  \!%
\Big\rangle
.
\end{equation}

\end{proposition}

The remainder of this section is devoted to the proof of this proposition. Set%
\begin{equation}
M_{1}(x)\ :=\ \mathsf{E}\left(  _{\!_{\!_{\,}}}k^{\varkappa}S_{t}%
\varphi(x)-m_{k}(t,x)\right)  \quad\text{for }\,x\in\mathbb{R}^{d},
\end{equation}
and for $\,y\in\mathbb{R}^{d},\ \,0\leq s\leq t,$%
\begin{equation}
I_{s}(y,W)\ :=\ \int_{0}^{t-s}\!dr\ \vartheta_{1}(kW_{r}-y)\,S_{t-s-r}%
\varphi\,(W_{r})\ \geq\ 0. \label{not.I}%
\end{equation}

\begin{lemma}
[\textbf{Dropping the exponential}]\label{L.dropping}For each $\,\delta
>0$\thinspace\ and for $\,c_{\ref{la2}}$\thinspace\ from Lemma~\emph{\ref{la2}%
},%
\begin{align}
&  \bigg|M_{1}(x)\,-\,k^{2\gamma-2}\,\gamma\,\varrho^{\gamma}\int
\!\!dz\ \mathcal{E}_{x}\int_{0}^{t}\!{d}s\;\vartheta_{1}(kW_{s}%
-z)\,\mathcal{E}_{W_{s}}\varphi(W_{t-s}^{1})\,\mathcal{E}_{W_{s}}%
\varphi(W_{t-s}^{2})\nonumber\\
&  \qquad\qquad\times\bigl(I_{s}(z,W^{1})+I_{s}(z,W^{2})\bigr)^{\!\gamma
-1}\bigg|\ \leq\ c_{\ref{la2}}\ 2\varrho^{2\gamma}\,k^{\delta-\varkappa}%
\,\phi_{\gamma\lambda_{\ref{L.bo2}}}(x),
\end{align}
for all $\,x\in\mathbb{R}^{d}$\thinspace\ and $\,k\geq1.$
\end{lemma}

%

\proof
By (\ref{heat}) and the Feynman-Kac representation (\ref{linprob}),%
\begin{align}
&  \mathsf{E}\!\left(  _{\!_{\!_{\,}}}k^{\varkappa}S_{t}\varphi\,(x)-m_{k}%
(t,x)\right)  \ =\ k^{2-d}\,\varrho\,\mathsf{E}\mathcal{E}_{x}\int_{0}%
^{t}\!{d}s\;\Gamma^{1}(kW_{s})\,w_{k}^{2}(t-s,W_{s})\\
&  =\ k^{2-d+2\varkappa}\,\varrho\,\mathsf{E}\mathcal{E}_{x}\int_{0}^{t}%
\!{d}s\;\Gamma^{1}(kW_{s})\,\mathcal{E}_{W_{s}}\varphi(W_{t-s}^{1}%
)\,\mathcal{E}_{W_{s}}\varphi(W_{t-s}^{2})\nonumber\\
&  \qquad\times\exp\!\bigg[-k^{2-d+\varkappa}\varrho\int_{0}^{t-s}%
\!{d}r\;\Gamma^{1}(kW_{r}^{1})\,S_{t-s-r}\varphi\,(W_{r}^{1})\nonumber\\
&  \qquad\qquad\quad\ -k^{2-d+\varkappa}\varrho\int_{0}^{t-s}\!{d}%
r\;\Gamma^{1}(kW_{r}^{2})\,S_{t-s-r}\varphi\,(W_{r}^{2})\bigg],\nonumber
\end{align}
where $\,W^{1}$\thinspace\ and $\,W^{2}$\thinspace\ are independent Brownian
motions starting from $\,W_{s\,}.$\thinspace\ By the definition
(\ref{not.Gamma.eps}) of $\,\Gamma^{1}$\thinspace\ this equals%
\begin{gather}
k^{2-d+2\varkappa}\,\varrho\,\mathsf{E}\int\!\!\Gamma(dz)\,\mathcal{E}_{x}%
\int_{0}^{t}\!{d}s\;\vartheta_{1}(kW_{s}-z)\,\mathcal{E}_{W_{s}}%
\varphi(W_{t-s}^{1})\,\mathcal{E}_{W_{s}}\varphi(W_{t-s}^{2}%
)\label{monstrum.new}\\
\times\exp\!\bigg[-\int\!\!\Gamma(dy)\,k^{2-d+\varkappa}\varrho\,\Big(I_{s}%
(y,W^{1})+I_{s}(y,W^{2})\Big)\bigg].\nonumber
\end{gather}
Recall that for measurable $\,\varphi,\psi\geq0,$%
\begin{equation}
\mathsf{E}\left\langle \Gamma,\varphi\right\rangle \,\mathrm{e}^{-\left\langle
\Gamma,\psi\right\rangle }\ =\ \gamma\int\!\!dz\ \varphi(z)\,\psi^{\gamma
-1}(z)\,\exp\!\Big[-\int\!\!dy\ \psi^{\gamma}(y)\Big ] \label{Gamma.formula}%
\end{equation}
(cf.\ \cite[Section~4]{DawsonFleischmann1992.equ}) and $\,k^{2-d+2\varkappa
}k^{(2-d+\varkappa)(\gamma-1)}=k^{2\gamma-2}$\thinspace\ for $\,\varkappa
=\varkappa_{\mathrm{c}\,}.$\thinspace\ Applying this to (\ref{monstrum.new})
yields%
\begin{align}
&  k^{2\gamma-2}\,\gamma\,\varrho^{\gamma}\int\!\!dz\ \mathcal{E}_{x}\int
_{0}^{t}\!{d}s\;\vartheta_{1}(kW_{s}-z)\,\mathcal{E}_{W_{s}}\varphi
(W_{t-s}^{1})\,\mathcal{E}_{W_{s}}\varphi(W_{t-s}^{2})\nonumber\\
&  \times\bigl(I_{s}(z,W^{1})+I_{s}(z,W^{2})\bigr)^{\!\gamma-1}\,\exp
\!\Big[-\int\!\!dy\ k^{(2-d+\varkappa)\gamma}\varrho^{\gamma}\bigl(I_{s}%
(y,W^{1})+I_{s}(y,W^{2})\bigr)^{\!\gamma}\Big ].\nonumber
\end{align}
By the inequality $\,1-\mathrm{e}^{-a}\leq a$\thinspace\ we have%
\begin{align*}
&  \bigg|M_{1}(x)\,-\,k^{2\gamma-2}\,\gamma\,\varrho^{\gamma}\!\int
\!\!dz\ \mathcal{E}_{x}\int_{0}^{t}\!{d}s\;\vartheta_{1}(kW_{s}%
-z)\,\mathcal{E}_{W_{s}}\varphi(W_{t-s}^{1})\,\mathcal{E}_{W_{s}}%
\varphi(W_{t-s}^{2})\\
&  \qquad\qquad\qquad\qquad\qquad\qquad\qquad\qquad\qquad\quad\times
\bigl(I_{s}(z,W^{1})+I_{s}(z,W^{2})\bigr)^{\!\gamma-1}\bigg|\\
&  \leq\ k^{2\gamma-2}k^{(2-d+\varkappa)\gamma}\,\gamma\,\varrho^{\gamma
}\!\int\!\!dz\ \mathcal{E}_{x}\int_{0}^{t}\!{d}s\;\vartheta_{1}(kW_{s}%
-z)\,\mathcal{E}_{W_{s}}\varphi(W_{t-s}^{1})\,\mathcal{E}_{W_{s}}%
\varphi(W_{t-s}^{2})\\
&  \qquad\qquad\quad\quad\times\bigl(I_{s}(z,W^{1})+I_{s}(z,W^{2}%
)\bigr)^{\!\gamma-1}\!\int\!\!dy\ \varrho^{\gamma}\bigl(I_{s}(y,W^{1}%
)+I_{s}(y,W^{2})\bigr)^{\!\gamma}.
\end{align*}
Applying (\ref{elementary}) to the last integrand and using the symmetry in
$W^{1},W^{2},$\thinspace\ we see that the right hand side in the former
display does not exceed%
\begin{align*}
&  2k^{2\gamma-2}k^{(2-d+\varkappa)\gamma}\,\gamma\,\varrho^{2\gamma}%
\int\!\!dz\int\!\!dy\ \ \mathcal{E}_{x}\int_{0}^{t}\!{d}s\;\vartheta
_{1}(kW_{s}-z)\,\mathcal{E}_{W_{s}}\varphi(W_{t-s}^{1})\,\mathcal{E}_{W_{s}%
}\varphi(W_{t-s}^{2})\\
&  \qquad\qquad\qquad\qquad\qquad\qquad\qquad\quad\quad\quad\times
\bigl(I_{s}(z,W^{1})+I_{s}(z,W^{2})\bigr)^{\!\gamma-1}I_{s}(y,W^{1})^{\gamma}.
\end{align*}
We now drop $\,I_{s}(z,W^{2})$\thinspace\ and evaluate the expectation with
respect to $\,W^{2},$\thinspace\ obtaining the upper bound%
\begin{align}
&  2k^{2\gamma-2}k^{(2-d+\varkappa)\gamma}\,\gamma\,\varrho^{2\gamma}%
\int\!\!dz\int\!\!dy\ \ \mathcal{E}_{x}\int_{0}^{t}\!{d}s\;\vartheta
_{1}(kW_{s}-z)\,S_{t-s}\varphi\,(W_{s})\\
&  \qquad\qquad\qquad\qquad\qquad\quad\ \times\mathcal{E}_{W_{s}}%
\varphi(W_{t-s}^{1})\,I_{s}(z,W^{1})^{\gamma-1}I_{s}(y,W^{1})^{\gamma
}.\nonumber
\end{align}
Applying the Markov property at time $s$ and time-homogeneity, this equals%
\begin{align*}
&  2k^{2\gamma-2}k^{(2-d+\varkappa)\gamma}\,\gamma\,\varrho^{2\gamma}%
\int\!\!dz\int\!\!dy\ \ \mathcal{E}_{x}\int_{0}^{t}\!{d}s\;\vartheta
_{1}(kW_{s}-z)\,S_{t-s}\varphi\,(W_{s})\\
&  \times\varphi(W_{t})\Big(\int_{s}^{t}\!dr\ \vartheta_{1}(kW_{r}%
-z)\,S_{t-r}\varphi\,(W_{r})\Big)^{\!\gamma-1}\Big(\int_{s}^{t}\!dr\ \vartheta
_{1}(kW_{r}-y)\,S_{t-r}\varphi\,(W_{r})\Big)^{\!\gamma}\!.
\end{align*}
The last factor can be bounded by $\,\left(  \,I_{0}(y,W)\right)  ^{\!\gamma
}.$\thinspace\ Then we integrate with respect to $s$ and obtain%
\begin{align*}
&  2k^{2\gamma-2}k^{(2-d+\varkappa)\gamma}\,\varrho^{2\gamma}\int
\!\!dz\int\!\!dy\ \ \mathcal{E}_{x}\varphi(W_{t})\Big(\int_{0}^{t}%
\!dr\ \vartheta_{1}(kW_{r}-z)\,S_{t-r}\varphi\,(W_{r})\Big)^{\!\gamma}\\
&  \times\left(  \,I_{0}(y,W)\right)  ^{\!\gamma}\ =\ 2k^{2\gamma
-2}k^{(2-d+\varkappa)\gamma}\,\varrho^{2\gamma}\mathcal{E}_{x}\varphi
(W_{t})\Big[\int\!\!dy\ \ I_{0}(y,W)^{\gamma}\Big ]^{2}.
\end{align*}
Using now Lemma~\ref{la2}, we arrive at%
\begin{align*}
&  \bigg|M_{1}(x)\,-\,k^{2\gamma-2}\,\gamma\,\varrho^{\gamma}\int
\!\!dz\ \mathcal{E}_{x}\int_{0}^{t}\!{d}s\;\vartheta_{1}(kW_{s}%
-z)\,\mathcal{E}_{W_{s}}\varphi(W_{t-s}^{1})\,\mathcal{E}_{W_{s}}%
\varphi(W_{t-s}^{2})\\
&  \qquad\qquad\qquad\qquad\qquad\qquad\qquad\qquad\qquad\quad\times
\bigl(I_{s}(z,W^{1})+I_{s}(z,W^{2})\bigr)^{\!\gamma-1}\bigg|\\
&  \leq\ c_{\ref{la2}}\ 2\varrho^{2\gamma}\,k^{2\gamma-2}k^{(2-d+\varkappa
)\gamma}\,k^{4-4\gamma+\delta}\,\phi_{\gamma\lambda_{\ref{L.bo2}}%
}(x)\ =\ c_{\ref{la2}}\ 2\varrho^{2\gamma}\,k^{\delta-\varkappa}\,\phi
_{\gamma\lambda_{\ref{L.bo2}}}(x),
\end{align*}
finishing the proof.%
\endproof

It remains to find the limit of%
\begin{align}
&  k^{2\gamma-2}\,\gamma\,\varrho^{\gamma}\int\!\!dz\ \mathcal{E}_{x}\int
_{0}^{t}\!{d}s\;\vartheta_{1}(kW_{s}-z)\,\mathcal{E}_{W_{s}}\varphi
(W_{t-s}^{1})\,\mathcal{E}_{W_{s}}\varphi(W_{t-s}^{2})\\
&  \qquad\qquad\qquad\qquad\qquad\qquad\qquad\ \,\times\bigl(I_{s}%
(z,W^{1})+I_{s}(z,W^{2})\bigr)^{\!\gamma-1}.\nonumber
\end{align}
Substituting $\,z\rightsquigarrow kz$\thinspace\ gives%
\begin{align*}
&  k^{2\gamma-2+d}\,\gamma\,\varrho^{\gamma}\int\!\!dz\ \mathcal{E}_{x}%
\int_{0}^{t}\!{d}s\;\vartheta_{1/k}(W_{s}-z)\,\mathcal{E}_{W_{s}}%
\varphi(W_{t-s}^{1})\,\mathcal{E}_{W_{s}}\varphi(W_{t-s}^{2})\\
&  \times\biggl(\int_{0}^{t-s}\!dr\left[  _{\!_{\!_{\,_{{}}}}}\vartheta
_{1/k}(W_{r}^{1}-z)\,S_{t-s-r}\varphi\,(W_{r}^{1})+\vartheta_{1/k}(W_{r}%
^{2}-z)\,S_{t-s-r}\varphi\,(W_{r}^{2})\right]  \biggr)^{\!\!\gamma-1}.
\end{align*}
Fix $\,x,z\in\mathbb{R}^{d}$\thinspace\ and $\,0<s<t$\thinspace\ for a while
and consider%
\begin{align*}
&  g_{k}(s,x,z)\ :=\ k^{2\gamma-2+d}\,\gamma\,\mathcal{E}_{x}\,\vartheta
_{1/k}(W_{s}-z)\,\mathcal{E}_{W_{s}}\varphi(W_{t-s}^{1})\,\mathcal{E}_{W_{s}%
}\varphi(W_{t-s}^{2})\\
&  \times\biggl(\int_{0}^{t-s}\!dr\left[  _{\!_{\!_{\,_{{}}}}}\vartheta
_{1/k}(W_{r}^{1}-z)\,S_{t-s-r}\varphi\,(W_{r}^{1})+\vartheta_{1/k}(W_{r}%
^{2}-z)\,S_{t-s-r}\varphi\,(W_{r}^{2})\right]  \biggr)^{\!\!\gamma-1}.
\end{align*}

\begin{lemma}
\label{L.bound.of.g}%
\[
\liminf_{k\uparrow\infty}\,g_{k}(s,x,z)\ \geq\ \underline{c}\,p_{s}%
(x-z)\left(  S_{t-s}\varphi\,(z)\right)  ^{\!1+\gamma}.
\]

\end{lemma}

The lemma immediately implies Proposition~\ref{L.conv.exp}. Indeed, applying
Fatou's lemma we get
\begin{align*}
\liminf_{k\uparrow\infty}\int\!\!\mu(dx)\int\!\!dz\int_{0}^{t}\!ds\ g_{k}%
(s,x,z)\  &  \geq\ \int\!\!\mu(dx)\int\!\!dz\int_{0}^{t}\!ds\ \liminf
_{k\uparrow\infty}g_{k}(s,x,z)\\
&  \geq\ \underline{c}\int_{0}^{t}\!ds\left\langle _{\!_{\!_{\,}}}\mu
,\,S_{s}(S_{t-s}\varphi)^{1+\gamma}\right\rangle .
\end{align*}

\noindent\emph{Proof of Lemma}~\ref{L.bound.of.g}. \thinspace Shifting the
Brownian motions,%
\begin{align*}
g_{k}(s,x,z)\,=\,\  &  k^{2\gamma-2+d}\,\gamma\,\mathcal{E}_{x}\vartheta
_{1/k}(W_{s}-z)\,\mathcal{E}_{z}\varphi(W_{t-s}^{1}+W_{s}-z)\,\mathcal{E}%
_{z}\varphi(W_{t-s}^{2}+W_{s}-z)\\
&  \qquad\qquad\,\ \times\Big(\int_{0}^{t-s}\!dr\ \vartheta_{1/k}(W_{r}%
^{1}+W_{s}-2z)\,S_{t-s-r}\varphi\,(W_{r}^{1}+W_{s}-z)\\
&  +\int_{0}^{t-s}\!dr\ \vartheta_{1/k}(W_{r}^{2}+W_{s}-2z)\,S_{t-s-r}%
\varphi\,(W_{r}^{2}+W_{s}-z)\Big)^{\!\gamma-1}.
\end{align*}
By the uniform continuity of $\varphi,$%
\begin{equation}
\lim_{k\uparrow\infty}\,\sup_{\left\vert W_{s}-z\right\vert \leq
1/k}\,\left\vert _{\!_{\!_{\,}}}\varphi(W_{t-s}^{i}+W_{s}-z)-\varphi
(W_{t-s}^{i})\right\vert \ =\ 0,
\end{equation}
and by (\ref{unicont}),%
\begin{equation}
\lim_{k\uparrow\infty}\,\sup_{\left\vert W_{s}-z\right\vert \leq
1/k}\,\left\vert _{\!_{\!_{\,}}}S_{t-s-r}\varphi(W_{t-s}^{i}+W_{s}%
-z)-S_{t-s-r}\varphi(W_{t-s}^{i})\right\vert \ =\ 0, \label{unicont'}%
\end{equation}
we get%
\begin{align*}
g_{k}(s,x,z)\ =\  &  k^{2\gamma-2+d}\,\left(  _{\!_{\!_{\,}}}\gamma
+o(1)\right)  \,\mathcal{E}_{x}\,\vartheta_{1/k}(W_{s}-z)\,\mathcal{E}%
_{z}\varphi(W_{t-s}^{1})\,\mathcal{E}_{z}\varphi(W_{t-s}^{2})\\
&  \qquad\qquad\ \,\times\Big(\int_{0}^{t-s}\!dr\ \vartheta_{1/k}(W_{r}%
^{1}+W_{s}-2z)\,S_{t-s-r}\varphi\,(W_{r}^{1})\\
&  +\int_{0}^{t-s}\!dr\ \vartheta_{1/k}(W_{r}^{2}+W_{s}-2z)\,S_{t-s-r}%
\varphi\,(W_{r}^{2})\Big)^{\!\gamma-1}.
\end{align*}
By the triangle inequality, $\,\vartheta_{1/k}(W_{r}^{i}+W_{s}-2z)\leq
\vartheta_{2/k}(W_{s}^{i}-z).$\thinspace\ Hence,%
\begin{align*}
g_{k}(s,x,z)\ \geq\  &  k^{2\gamma-2+d}\,\left(  _{\!_{\!_{\,}}}%
\gamma+o(1)\right)  \,\mathcal{E}_{x}\,\vartheta_{1/k}(W_{s}-z)\,\mathcal{E}%
_{z}\varphi(W_{t-s}^{1})\,\mathcal{E}_{z}\varphi(W_{t-s}^{2})\\
&  \quad\,\ \times\biggl(\int_{0}^{t-s}\!dr\ \vartheta_{2/k}(W_{r}%
^{1}-z)\,S_{t-s-r}\varphi\,(W_{r}^{1})\\
&  \quad\quad\qquad+\int_{0}^{t-s}\!dr\ \vartheta_{2/k}(W_{r}^{2}%
-z)\,S_{t-s-r}\varphi\,(W_{r}^{2})\biggr)^{\!\!\gamma-1}.
\end{align*}
Calculating the expectation with respect to $W$ gives%
\begin{equation}
\mathcal{E}_{x}\,\vartheta_{1/k}(W_{s}-z)\ =\ \frac{\pi^{d/2}}{G(1+d/2)}%
\,k^{-d}\,p_{s}(x-z)\left(  _{\!_{\!_{\,}}}1+o(1)\right)  \!.
\end{equation}
Using (\ref{unicont'}) once more we obtain%
\begin{align*}
&  g_{k}(s,x,z)\ \geq\ k^{2\gamma-2}\,\left(  _{\!_{\!_{\,}}}\gamma
+o(1)\right)  \,\frac{\pi^{d/2}}{G(1+d/2)}\,p_{s}(x-z)\,\mathcal{E}_{z}%
\varphi(W_{t-s}^{1})\,\mathcal{E}_{z}\varphi(W_{t-s}^{2})\\
&  \qquad\qquad\qquad\times\Big(\int_{0}^{t-s}\!dr\ \bigl[\vartheta
_{2/k}(W_{r}^{1}-z)+\vartheta_{2/k}(W_{r}^{2}-z)\bigr]\,S_{t-s-r}%
\varphi\,(z)\Big)^{\!\gamma-1}.
\end{align*}
Define events%
\begin{equation}
A_{k}^{i}(z)\ :=\ \left\{  _{\!_{\!_{\,}}}|W_{r}^{i}-z|>1/k\ \forall
r>1/k\right\}  \!.
\end{equation}
Evidently,%
\begin{align*}
&  \mathcal{E}_{z}\varphi(W_{t-s}^{1})\,\mathcal{E}_{z}\varphi(W_{t-s}%
^{2})\Big(\int_{0}^{t-s}\!dr\ \bigl[\vartheta_{2/k}(W_{r}^{1}-z)+\vartheta
_{2/k}(W_{r}^{2}-z)\bigr]\,S_{t-s-r}\varphi\,(z)\Big)^{\!\gamma-1}\\
&  \geq\mathcal{E}_{z}\varphi(W_{t-s}^{1})\,\mathcal{E}_{z}\varphi(W_{t-s}%
^{2})\\
&  \quad\ \times\Big(\int_{0}^{1/k}\!dr\ \bigl[\vartheta_{2/k}(W_{r}%
^{1}-z)+\vartheta_{2/k}(W_{r}^{2}-z)\bigr]\,S_{t-s-r}\varphi
\,(z)\Big)^{\!\gamma-1}\mathsf{1}_{A_{k}^{1}(z)}\mathsf{1}_{A_{k}^{2}(z)}\\
&  \geq\mathcal{E}_{z}\varphi(W_{t-s}^{1})\,\mathcal{E}_{z}\varphi(W_{t-s}%
^{2})\Big(\int_{0}^{1/k}\!\!dr\,\bigl[\vartheta_{2/k}(W_{r}^{1}-z)+\vartheta
_{2/k}(W_{r}^{2}-z)\bigr]\,S_{t-s-r}\varphi(z)\Big)^{\!\gamma-1}\\
&  \quad\ -2\mathcal{E}_{z}\varphi(W_{t-s}^{1})\,\mathcal{E}_{z}%
\varphi(W_{t-s}^{2})\Big(\int_{0}^{1/k}\!dr\ \vartheta_{2/k}(W_{r}%
^{1}-z)\,S_{t-s-r}\varphi\,(z)\Big)^{\!\gamma-1}(1-\mathsf{1}_{A_{k}^{1}(z)}).
\end{align*}
We calculate the expressions in the last two lines separately. For the first
line we get, by the Markov property at time $1/k,$%
\begin{align*}
&  \mathcal{E}_{z}\varphi(W_{t-s}^{1})\,\mathcal{E}_{z}\varphi(W_{t-s}%
^{2})\Big(\int_{0}^{1/k}\!dr\ \bigl[\vartheta_{2/k}(W_{r}^{1}-z)+\vartheta
_{2/k}(W_{r}^{2}-z)\bigr]\,S_{t-s-r}\varphi\,(z)\Big)^{\!\gamma-1}\\
&  =\mathcal{E}_{0}\otimes\mathcal{E}_{0}\Big(\int_{0}^{1/k}%
\!dr\ \bigl[\vartheta_{2/k}(W_{r}^{1})+\vartheta_{2/k}(W_{r}^{2}%
)\bigr]\,S_{t-s-r}\varphi\,(z)\Big)^{\!\gamma-1}\\
&  \qquad\qquad\qquad\qquad\qquad\qquad\qquad\ \ \times S_{t-s-1/k}%
\varphi\,(z+W_{1/k}^{1})\,S_{t-s-1/k}\varphi\,(z+W_{1/k}^{2}).
\end{align*}
By (\ref{unicont}), this equals asymptotically
\begin{align*}
&  \left(  _{\!_{\!_{\,}}}S_{t-s}\varphi\,(z)\right)  ^{1+\gamma}%
\mathcal{E}_{0}\otimes\mathcal{E}_{0}\Big(\int_{0}^{1/k}\!dr\ \bigl[\vartheta
_{2/k}(W_{r}^{1})+\vartheta_{2/k}(W_{r}^{2})\bigr]\,\Big)^{\!\gamma-1}\\
&  =\,\left(  _{\!_{\!_{\,}}}S_{t-s}\varphi\,(z)\right)  ^{1+\gamma
}k^{2-2\gamma}\,\mathcal{E}_{0}\otimes\mathcal{E}_{0}\Big(\int_{0}%
^{k}\!dr\ \bigl[\vartheta_{2}(W_{r}^{1})+\vartheta_{2}(W_{r}^{2}%
)\bigr]\,\Big)^{\!\gamma-1},
\end{align*}
where in the last step Brownian scaling was used. Therefore the first line is
asymptotically equivalent to%
\begin{equation}
\left(  _{\!_{\!_{\,}}}S_{t-s}\varphi\,(z)\right)  ^{1+\gamma}k^{2-2\gamma
}\mathcal{E}_{0}\otimes\mathcal{E}_{0}\Big(\int_{0}^{\infty}\!\!dr\,\left[
\vartheta_{2}(W_{r}^{1})+\vartheta_{2}(W_{r}^{2})\right]  \Big)^{\!\gamma-1}.
\label{first.line}%
\end{equation}
Turning now to the second line,%
\begin{align*}
&  2\mathcal{E}_{z}\varphi(W_{t-s}^{1})\,\mathcal{E}_{z}\varphi(W_{t-s}%
^{2})\Big(\int_{0}^{1/k}\!dr\ \vartheta_{2/k}(W_{r}^{1}-z)\,S_{t-s-r}%
\varphi\,(z)\Big)^{\!\gamma-1}(1-\mathsf{1}_{A_{k}^{1}(z)})\\
&  =\ 2\left(  _{\!_{\!_{\,}}}S_{t-s}\varphi\,(z)\right)  ^{\!\gamma}\,\left(
_{\!_{\!_{\,}}}1+o(1)\right)  \mathcal{E}_{z}\varphi(W_{t-s}^{1})\Big(\int
_{0}^{1/k}\!dr\ \vartheta_{2/k}(W_{r}^{1}-z)\Big)^{\!\gamma-1}(1-\mathsf{1}%
_{A_{k}^{1}(z)}),
\end{align*}
where the expectation with respect to $W^{2}$ was evaluated, and
(\ref{unicont}) was used. Recalling that $\varphi$ is bounded and applying
Cauchy-Schwarz we obtain an upper bound%
\begin{align}
&  c_{(\ref{second.line})}\left(  \mathcal{P}_{0}\!\left(  _{\!_{\!_{\,}}%
}A_{k}^{1}(0)^{\mathrm{c}}\right)  \right)  ^{1/2}\bigg[\mathcal{E}%
_{0}\Big(\int_{0}^{1/k}\!dr\ \vartheta_{2/k}(W_{r}^{1})\Big)^{\!2\gamma
-2}\bigg]^{1/2}\label{second.line}\\
&  =\ c_{(\ref{second.line})}\,k^{2-2\gamma}\left[  _{\!_{\!_{\,_{{}}}}%
}\mathcal{P}_{0}\!\left(  _{\!_{\!_{\,}}}\exists r>k:\ |W_{r}|\leq1\right)
\right]  ^{1/2}\bigg[\mathcal{E}_{0}\Big(\int_{0}^{k}\!dr\ \vartheta_{2}%
(W_{r}^{1})\Big)^{\!2\gamma-2}\bigg]^{1/2}.\nonumber
\end{align}
Since the expectation is bounded and the probability goes to zero,
(\ref{second.line}) is $o(k^{2-2\gamma}).$\thinspace\ Together with
(\ref{first.line}) this proves the lemma.\hfill$\square$

\subsection{Convergence of variances\label{SS.var}}

\begin{proposition}
[\textbf{Convergence of variances}]\label{L.var'}For every $\,\mu
\in\mathcal{M}_{\mathrm{tem}}\,\ $satisfying the assumption in Theorem\/
\emph{\ref{T.fluct},} for $\,\varphi\in\mathcal{C}_{\exp}^{+}$ and $\,t>0$,
\begin{equation}
\lim_{k\rightarrow\infty}\mathsf{V}\mathrm{ar}\int\!\!\mu
(dx)\,\Big [k^{\varkappa}S_{t}\varphi(x)-m_{k}(t,x)\Big]\ =\ 0.
\end{equation}

\end{proposition}

The remainder of this section is devoted to the proof of this proposition. We
may set $\,\varrho=1.$\thinspace\ Recall that%
\begin{align*}
&  k^{\varkappa}S_{t}\varphi(x)-m_{k}(t,x)\ =\ k^{2-d+2\varkappa}%
\,\mathcal{E}_{x}\int_{0}^{t}\!ds\int\!\!\Gamma(dz)\ \vartheta_{1}%
(kW_{s}-z)\,\\
&  \times\mathcal{E}_{W_{s}}\otimes\mathcal{E}_{W_{s}}\varphi(W_{t-s}%
^{1})\,\varphi(W_{t-s}^{2})\exp\Big[-k^{2-d+\varkappa}\int\!\!\Gamma
(dz)\left(  I_{s}(z,W^{1})+I_{s}(z,W^{2})\right)  \Big ].
\end{align*}
Define%
\begin{equation}
M_{2}(x,\tilde{x})\ :=\ \mathsf{E}\left(  _{\!_{\!_{\,}}}k^{\varkappa}%
S_{t}\varphi(x)-m_{k}(t,x)\right)  \left(  _{\!_{\!_{\,}}}k^{\varkappa}%
S_{t}\varphi(\tilde{x})-m_{k}(t,\tilde{x})\right)  \!.
\end{equation}
Similarly to (\ref{Gamma.formula}), for measurable $\,\varphi_{1},\varphi
_{2},\psi\geq0,$%
\begin{align*}
&  \mathsf{E}\left\langle \Gamma,\varphi_{1}\right\rangle \left\langle
\Gamma,\varphi_{2}\right\rangle \,\mathrm{e}^{-\left\langle \Gamma
,\psi\right\rangle }\ =\ \gamma(1-\gamma)\int\!\!dz\ \varphi_{1}%
(z)\,\varphi_{2}(z)\,\psi^{\gamma-2}(z)\,\exp\!\Big[-\int\!\!dy\ \psi^{\gamma
}(y)\Big ]\\
&  \qquad\qquad\ \ +\ \gamma^{2}\int\!\!dz_{1}\ \varphi_{1}(z_{1}%
)\,\psi^{\gamma-1}(z_{1})\int\!\!dz_{2}\ \varphi_{2}(z_{2})\,\psi^{\gamma
-1}(z_{2})\,\exp\!\Big[-\int\!\!dy\ \psi^{\gamma}(y)\Big ].
\end{align*}
Applying this formula, we get%
\begin{equation}
M_{2}(x,\tilde{x})\ =\ M_{21}(x,\tilde{x})+M_{22}(x,\tilde{x}),
\end{equation}
where%
\begin{align*}
&  M_{21}(x,\tilde{x})\ :=\ \gamma(1-\gamma)\,k^{4-2d+4\varkappa
}\,k^{(2-d+\varkappa)(\gamma-2)}\,\mathcal{E}_{x}\otimes\mathcal{E}_{\tilde
{x}}\int_{0}^{t}\!ds\int_{0}^{t}\!d\tilde{s}\int\!\!dz\ \\
&  \vartheta_{1}(kW_{s}-z)\,\vartheta_{1}(k\tilde{W}_{\tilde{s}}%
-z)\,\mathcal{E}_{W_{s}}\!\otimes\mathcal{E}_{W_{s}}\varphi(W_{t-s}%
^{1})\varphi(W_{t-s}^{2})\,\mathcal{E}_{\tilde{W}_{\tilde{s}}}\!\otimes
\mathcal{E}_{\tilde{W}_{\tilde{s}}}\varphi(\tilde{W}_{t-\tilde{s}}^{1}%
)\varphi(\tilde{W}_{t-\tilde{s}}^{2})\\
&  \qquad\times\left(  _{\!_{\!_{\,}}}I_{s}(z,W^{1})+I_{s}(z,W^{2}%
)+I_{\tilde{s}}(z,\tilde{W}^{1})+I_{\tilde{s}}(z,\tilde{W}^{2})\right)
^{\!\gamma-2}\,\\
&  \qquad\times\exp\!\bigg[-k^{\gamma(2-d+\varkappa)}\int\!\!dy\,\left(
_{\!_{\!_{\,}}}I_{s}(y,W^{1})+I_{s}(y,W^{2})+I_{\tilde{s}}(y,\tilde{W}%
^{1})+I_{\tilde{s}}(y,\tilde{W}^{2})\right)  ^{\!\gamma}\bigg]
\end{align*}
and%
\begin{align*}
&  M_{22}(x,\tilde{x})\ :=\ \gamma^{2}\,k^{4-2d+4\varkappa}\,k^{(2-d+\varkappa
)(2\gamma-2)}\,\mathcal{E}_{x}\otimes\mathcal{E}_{\tilde{x}}\int_{0}%
^{t}\!ds\int_{0}^{t}\!d\tilde{s}\int\!\!dz\int\!\!d\tilde{z}\\
&  \vartheta_{1}(kW_{s}-z)\,\vartheta_{1}(k\tilde{W}_{\tilde{s}}-\tilde
{z})\,\mathcal{E}_{W_{s}}\!\otimes\mathcal{E}_{W_{s}}\varphi(W_{t-s}%
^{1})\varphi(W_{t-s}^{2})\,\mathcal{E}_{\tilde{W}_{\tilde{s}}}\!\otimes
\mathcal{E}_{\tilde{W}_{\tilde{s}}}\varphi(\tilde{W}_{t-\tilde{s}}^{1}%
)\varphi(\tilde{W}_{t-\tilde{s}}^{2})\\
&  \qquad\times\left(  _{\!_{\!_{\,}}}I_{s}(z,W^{1})+I_{s}(z,W^{2}%
)+I_{\tilde{s}}(z,\tilde{W}^{1})+I_{\tilde{s}}(z,\tilde{W}^{2})\right)
^{\gamma-1}\\
&  \qquad\times\left(  _{\!_{\!_{\,}}}I_{s}(\tilde{z},W^{1})+I_{s}(\tilde
{z},W^{2})+I_{\tilde{s}}(\tilde{z},\tilde{W}^{1})+I_{\tilde{s}}(\tilde
{z},\tilde{W}^{2})\right)  ^{\gamma-1}\\
&  \qquad\times\exp\!\bigg[-\!k^{\gamma(2-d+\varkappa)}\!\int\!\!dy\,\left(
_{\!_{\!_{\,}}}I_{s}(y,W^{1})+I_{s}(y,W^{2})+I_{\tilde{s}}(y,\tilde{W}%
^{1})+I_{\tilde{s}}(y,\tilde{W}^{2})\right)  ^{\!\gamma}\bigg].
\end{align*}
The following Lemmas~\ref{L.M21} and \ref{L.M22} together directly imply
Proposition~\ref{L.var'}.

\begin{lemma}
\label{L.M21}%
\[
\lim_{k\uparrow\infty}\int\!\!\mu(dx)\int\!\!\mu(d\tilde{x})\,M_{21}%
(x,\tilde{x})\ =\ 0
\]

\end{lemma}

\begin{lemma}
\label{L.M22}%
\[
\limsup_{k\uparrow\infty}\int\!\!\mu(dx)\int\!\!\mu(d\tilde{x})\,\left[
M_{22}(x,\tilde{x})-M_{1}(x)M_{1}(\tilde{x})\right]  \ \leq\ 0.
\]

\end{lemma}

\noindent\emph{Proof of Lemma}~\ref{L.M21}. \thinspace By definition of the
critical index $\varkappa=\varkappa_{\mathrm{c}\,},$%
\begin{equation}
4-2d+4\varkappa+(\gamma-2)(2-d+\varkappa)\,=\,\varkappa-2+2\gamma.
\end{equation}
Dropping the exponential in $M_{21}(x,\tilde{x})$ and $I_{s}(z,W^{2}%
)+I_{s}(z,W^{2})$ gives%
\begin{align*}
&  M_{21}(x,\tilde{x})\ \leq\ k^{\varkappa-2+2\gamma}\gamma(1-\gamma
)\,\mathcal{E}_{x}\otimes\mathcal{E}_{\tilde{x}}\int_{0}^{t}\!ds\int_{0}%
^{t}\!d\tilde{s}\int\!\!dz\ \vartheta_{1}(kW_{s}-z)\,\vartheta_{1}(k\tilde
{W}_{\tilde{s}}-z)\,\\
&  \mathcal{E}_{W_{s}}\!\!\otimes\mathcal{E}_{W_{s}}\varphi(W_{t-s}%
^{1})\varphi(W_{t-s}^{2})\,\mathcal{E}_{\tilde{W}_{\tilde{s}}}\!\!\otimes
\mathcal{E}_{\tilde{W}_{\tilde{s}}}\varphi(\tilde{W}_{t-\tilde{s}}^{1}%
)\varphi(\tilde{W}_{t-\tilde{s}}^{2})\!\bigl(I_{s}(z,W^{1})\!+\!I_{\tilde{s}%
}(z,\tilde{W}^{1})\bigr)^{\!\gamma-2}\!.
\end{align*}
By independence of all Brownian paths, it follows that the expression in the
second line in the previous formula is bounded by
\[
S_{t-s}\varphi(W_{s})\,S_{t-\tilde{s}}\varphi(\tilde{W}_{\tilde{s}%
})\mathcal{E}_{W_{s}}\varphi(W_{t-s}^{1})\mathcal{E}_{\tilde{W}_{\tilde{s}}%
}\varphi(\tilde{W}_{t-\tilde{s}}^{1})\left(  \!_{\!_{\!_{\,}}}I_{s}%
(z,W^{1})+I_{\tilde{s}}(z,\tilde{W}^{1})\right)  ^{\!\gamma-2}\!.
\]
By the Markov property,
\begin{align*}
&  M_{21}(x,\tilde{x})\ \leq\ k^{\varkappa-2+2\gamma}\gamma(1-\gamma
)\,\mathcal{E}_{x}\otimes\mathcal{E}_{\tilde{x}}\varphi(W_{t})\varphi
(\tilde{W}_{t})\\
&  \qquad\times\int\!\!dz\int_{0}^{t}\!ds\int_{0}^{t}\!d\tilde{s}%
\ \vartheta_{1}(kW_{s}-z)\,\vartheta_{1}(k\tilde{W}_{\tilde{s}}-z)\,S_{t-s}%
\varphi(W_{s})\,S_{t-\tilde{s}}\varphi(\tilde{W}_{\tilde{s}})\\
&  \qquad\times\left(  \int_{s}^{t}\!dr\ \vartheta_{1}(kW_{r}-z)\,S_{t-r}%
\varphi(W_{r})+\int_{\tilde{s}}^{t}\!dr\ \vartheta_{1}(k\tilde{W}%
_{r}-z)\,S_{t-r}\varphi(\tilde{W}_{r})\right)  ^{\!\!\gamma-2}\!.
\end{align*}
Carrying out the integration over $s$ and $\tilde{s}$ gives%
\begin{align*}
&  k^{\varkappa-2+2\gamma}\,\mathcal{E}_{x}\otimes\mathcal{E}_{\tilde{x}%
}\varphi(W_{t})\varphi(\tilde{W}_{t})\\
&  \qquad\qquad\qquad\int\!\!dz\Big(I_{0}(z,W)^{\gamma}+I_{0}(z,\tilde
{W})^{\gamma}-\bigl(I_{0}(z,W)+I_{0}(z,\tilde{W})\bigr)^{\!\gamma}\Big).
\end{align*}
Changing the integration variable $k\rightsquigarrow kz$, we obtain%
\begin{align*}
&  \int\!\!\mu(dx)\int\!\!\mu(d\tilde{x})\,M_{21}(x,\tilde{x})\ \leq
\ k^{\varkappa-2+d}\int\!\!dz\int\!\!\mu(dx)\int\!\!\mu(d\tilde{x}%
)\ \mathcal{E}_{x}\otimes\mathcal{E}_{\tilde{x}}\varphi(W_{t})\varphi
(\tilde{W}_{t})\\
&  \times\Bigg[\Big(k^{2}\!\int_{0}^{t}\!{d}s\;\vartheta_{1/k}(W_{s}%
-z)\,S_{t-s}\varphi\,(W_{s})\Big)^{\!\gamma}+\Big(k^{2}\!\int_{0}^{t}%
\!{d}s\;\vartheta_{1/k}(\tilde{W}_{s}-z)\,S_{t-s}\varphi\,(\tilde{W}%
_{s})\Big)^{\!\gamma}\\
&  -\biggl(k^{2}\!\int_{0}^{t}\!{d}s\;\left[  \vartheta_{1/k}(W_{s}%
-z)\,S_{t-s}\varphi\,(W_{s})+\vartheta_{1/k}(\tilde{W}_{s}-z)\,S_{t-s}%
\varphi\,(\tilde{W}_{s})\right]  \biggr)^{\!\gamma}\Bigg].
\end{align*}
The right hand side of this inequality coincides with (\ref{121}), since
$\,2\varkappa+(\varkappa-d)\gamma+d=\varkappa-2+d,$\thinspace\ hence converges
to zero. \hfill$\square$\bigskip

\noindent\emph{Proof of Lemma}~\ref{L.M22}. \thinspace Dropping some
non-negative summands, we get%
\begin{align*}
&  M_{22}(x,\tilde{x})\ \leq\ k^{4\gamma-4}\,\gamma^{2}\,\mathcal{E}%
_{x}\otimes\mathcal{E}_{\tilde{x}}\int_{0}^{t}\!ds\int_{0}^{t}\!d\tilde{s}%
\int\!\!dz\int\!\!d\tilde{z}\\
&  \vartheta_{1}(kW_{s}-z)\,\vartheta_{1}(k\tilde{W}_{\tilde{s}}-\tilde
{z})\,\mathcal{E}_{W_{s}}\otimes\mathcal{E}_{W_{s}}\varphi(W_{t-s}^{1}%
)\varphi(W_{t-s}^{2})\,\mathcal{E}_{\tilde{W}_{\tilde{s}}}\otimes
\mathcal{E}_{\tilde{W}_{\tilde{s}}}\varphi(\tilde{W}_{t-\tilde{s}}^{1}%
)\varphi(\tilde{W}_{t-\tilde{s}}^{2})\\[4pt]
&  \qquad\times\left(  _{\!_{\!_{\,}}}I_{s}(z,W^{1})+I_{s}(z,W^{2})\right)
^{\gamma-1}\bigl(I_{\tilde{s}}(\tilde{z},\tilde{W}^{1})+I_{\tilde{s}}%
(\tilde{z},\tilde{W}^{2})\bigr)^{\!\gamma-1}\\
&  \qquad\times\exp\!\Big[-k^{\gamma(2-d+\varkappa)}\int\!\!dy\,\left(
_{\!_{\!_{\,}}}I_{s}(y,W^{1})+I_{s}(y,W^{2})+I_{\tilde{s}}(y,\tilde{W}%
^{1})+I_{\tilde{s}}(y,\tilde{W}^{2})\right)  ^{\!\gamma}\Big ].
\end{align*}
On the other hand,%
\begin{align}
&  M_{1}(x)\ =\ k^{2\gamma-2}\,\gamma\,\mathcal{E}_{x}\int_{0}^{t}\!{d}%
s\int\!\!dz\;\vartheta_{1}(kW_{s}-z)\,\mathcal{E}_{W_{s}}\varphi(W_{t-s}%
^{1})\,\mathcal{E}_{W_{s}}\varphi(W_{t-s}^{2})\nonumber\\
&  \times\bigl(I_{s}(z,W^{1})+I_{s}(z,W^{2})\bigr)^{\!\gamma-1}\exp
\!\bigg[-\!\int\!\!dy\ k^{(2-d+\varkappa)\gamma}\bigl(I_{s}(y,W^{1}%
)+I_{s}(y,W^{2})\bigr)^{\!\gamma}\bigg].\nonumber
\end{align}
Taking the difference, applying inequality (\ref{alsoelementary}) and using
symmetry, we get%
\begin{align*}
&  M_{22}(x,\tilde{x})-M_{1}(x)M_{1}(\tilde{x})\ \leq\ k^{4\gamma-4}%
\,\gamma^{2}\,\mathcal{E}_{x}\otimes\mathcal{E}_{\tilde{x}}\int_{0}%
^{t}\!ds\int_{0}^{t}\!d\tilde{s}\int\!\!dz\int\!\!d\tilde{z}\\
&  \vartheta_{1}(kW_{s}-z)\,\vartheta_{1}(k\tilde{W}_{\tilde{s}}-\tilde
{z})\,\mathcal{E}_{W_{s}}\!\otimes\mathcal{E}_{W_{s}}\varphi(W_{t-s}%
^{1})\varphi(W_{t-s}^{2})\,\mathcal{E}_{\tilde{W}_{\tilde{s}}}\!\otimes
\mathcal{E}_{\tilde{W}_{\tilde{s}}}\varphi(\tilde{W}_{t-\tilde{s}}^{1}%
)\varphi(\tilde{W}_{t-\tilde{s}}^{2})\\[4pt]
&  \qquad\qquad\qquad\qquad\qquad\times\left(  _{\!_{\!_{\,}}}I_{s}%
(z,W^{1})+I_{s}(z,W^{2})\right)  ^{\gamma-1}\bigl(I_{\tilde{s}}(\tilde
{z},\tilde{W}^{1})+I_{\tilde{s}}(\tilde{z},\tilde{W}^{2})\bigr)^{\!\gamma-1}\\
&  \qquad\qquad\qquad\qquad\qquad\qquad\qquad\qquad\qquad\qquad\quad\times
k^{(2-d+\varkappa)\gamma}\,4\int\!\!dy\ \left(  _{\!_{\!_{\,}}}I_{s}%
(y,W^{1})\right)  ^{\!\gamma}.
\end{align*}
Dropping further non-negative summands and using again the Markov property, we
get the bound%
\begin{align*}
&  4\,k^{4\gamma-4}k^{(2-d+\varkappa)\gamma}\,\gamma^{2}\,\mathcal{E}%
_{x}\otimes\mathcal{E}_{\tilde{x}}\int_{0}^{t}\!ds\int_{0}^{t}\!d\tilde{s}%
\int\!\!dz\int\!\!d\tilde{z}\\
&  \vartheta_{1}(kW_{s}-z)\,\vartheta_{1}(k\tilde{W}_{\tilde{s}}-\tilde
{z})\,S_{t-s}\varphi(W_{s})\,\,S_{t-\tilde{s}}\varphi(\tilde{W}_{\tilde{s}%
})\,\varphi(W_{t})\,\varphi(\tilde{W}_{t})\\
&  \times\Big(\int_{s}^{t}\!dr\ \vartheta_{1}(kW_{r}-z)\,S_{t-r}\varphi
(W_{r})\Big)^{\!\gamma-1}\Big(\int_{\tilde{s}}^{t}\!dr\ \vartheta_{1}%
(k\tilde{W}_{r}-\tilde{z})\,S_{t-r}\varphi(\tilde{W}_{r})\Big)^{\!\gamma-1}\\
&  \qquad\qquad\qquad\qquad\qquad\qquad\qquad\quad\times\int\!\!dy\ \Big(\int
_{0}^{t}\!dr\ \vartheta_{1}(kW_{r}-y)\,S_{t-r}\varphi(W_{r})\Big)^{\!\gamma}.
\end{align*}
Carrying out the $s$ and $\tilde{s}$ integration, we obtain%
\begin{align*}
&  4\,k^{4\gamma-4}k^{(2-d+\varkappa)\gamma}\,\mathcal{E}_{x}\otimes
\mathcal{E}_{\tilde{x}}\varphi(W_{t})\,\varphi(\tilde{W}_{t})\int
\!\!dz\int\!\!d\tilde{z}\\
&  \times\Big(\int_{0}^{t}\!dr\ \vartheta_{1}(kW_{r}-z)\,S_{t-r}\varphi
(W_{r})\Big)^{\!\gamma}\Big(\int_{0}^{t}\!dr\ \vartheta_{1}(k\tilde{W}%
_{r}-\tilde{z})\,S_{t-r}\varphi(\tilde{W}_{r})\Big)^{\!\gamma}\\
&  \qquad\qquad\qquad\qquad\qquad\qquad\,\ \times\int\!\!dy\ \Big(\int_{0}%
^{t}\!dr\ \vartheta_{1}(kW_{r}-y)\,S_{t-r}\varphi(W_{r})\Big)^{\!\gamma}.
\end{align*}
We collect identical terms, use the boundedness of $\,\varphi,$ and obtain, up
to a constant factor, the bound%
\begin{align}
&  k^{4\gamma-4}k^{(2-d+\varkappa)\gamma}\,\mathcal{E}_{x}\varphi
(W_{t})\bigg[\int\!\!dz\,\Big(\int_{0}^{t}\!dr\ \vartheta_{1}(kW_{r}%
-z)S_{t-r}\varphi(W_{r})\Big)^{\!\gamma}\bigg]^{2}\\
&  \qquad\qquad\qquad\qquad\qquad\qquad\times\mathcal{E}_{\tilde{x}}%
\,\varphi(\tilde{W}_{t})\int\!\!d\tilde{z}\,\Big(\int_{0}^{t}\!dr\ \vartheta
_{1}(k\tilde{W}_{r}-\tilde{z})\Big)^{\!\gamma}.\nonumber
\end{align}
By Lemma~\ref{L.Br.hit.occ}(b) with $\,\eta=\gamma$,%
\begin{align*}
&  \mathcal{E}_{\tilde{x}}\,\varphi(\tilde{W}_{t})\int\!\!d\tilde
{z}\,\Big(\int_{0}^{t}\!dr\ \vartheta_{1}(k\tilde{W}_{r}-\tilde{z}%
)\Big)^{\!\gamma}\ \\
&  =\ k^{d-2\gamma}\int\!\!d\tilde{z}\ \mathcal{E}_{\tilde{x}}\,\varphi
(\tilde{W}_{t})\Big(k^{2}\int_{0}^{t}\!dr\ \vartheta_{1}(k\tilde{W}%
_{r}-k\tilde{z})\Big)^{\!\gamma}\ \leq\ c_{\ref{L.Br.hit.occ}}\,k^{2-2\gamma
}\,\phi_{\lambda_{\ref{L.bo2}}}(\tilde{x}),
\end{align*}
and by Lemma~\ref{la2},%
\[
\mathcal{E}_{x}\varphi(W_{t})\bigg[\int\!\!dz\,\Big(\int_{0}^{t}%
\!dr\ \vartheta_{1}(kW_{r}-z)S_{t-r}\varphi(W_{r})\Big)^{\!\gamma}%
\bigg]^{2}\,\leq\ c_{\ref{la2}}\,k^{4-4\gamma+\delta}\,\phi_{\lambda
_{\ref{L.bo2}}}(x).
\]
Noting that $\,4\gamma-4+(2-d+\varkappa)\gamma+6-6\gamma=-\varkappa$%
\thinspace\ and choosing $\,\delta<\varkappa,$\thinspace\ we obtain, up to a
constant factor, the upper bound $\,k^{\delta-\varkappa}\phi_{\lambda
_{\ref{L.bo2}}}(x)\,\phi_{\lambda_{\ref{L.bo2}}}(\tilde{x}).$\thinspace\ The
proof is completed by integration.\hfill$\square$

\subsection{Lower bound for finite-dimensional distributions\label{SS.fdd}}

The proof is analogous to the upper bound in Section~\ref{findim}. Again we
use induction to show that, for any $\varphi_{1},\ldots,\varphi_{n}$ and
$0=t_{0}<t_{1}<\cdots<t_{n}$, in $\mathsf{P}$--probability,
\begin{equation}%
\begin{array}
[c]{l}%
\displaystyle
\liminf_{k\rightarrow\infty}\mathbb{E}_{\mu_{k}}\!\exp\!\Big[\sum_{i=1}%
^{n}k^{\varkappa}\big\langle X_{t_{i}}^{k}-S_{t_{i}}\mu,\,-\varphi
_{i}\big\rangle\Big ]\vspace{6pt}\\
\qquad%
\displaystyle
\geq\ \exp\!\Bigg[\underline{c}\,%
\bigg\langle
\mu,\sum_{i=1}^{n}\int_{t_{i-1}}^{t_{i}}dr\ S_{r}\Big(\big(\sum_{j=i}%
^{n}S_{t_{j}-r}\varphi_{j}\big)^{1+\gamma}\Big)\!%
\bigg\rangle
\Bigg].
\end{array}
\label{130}%
\end{equation}
For the case $n=1$ this was shown in the previous paragraphs, so we may assume
that it holds for $n-1$ and show that it also holds for $n$. By conditioning
on $\{X^{k}(t)\colon t\leq t_{n-1}\}$ and applying the transition functional
we get
\begin{align}
\mathbb{E}_{\mu_{k}}  &  \exp\!\Big[\sum_{i=1}^{n}k^{\varkappa}%
\big\langle X_{t_{i}}^{k}-S_{t_{i}}\mu,\,-\varphi_{i}%
\big\rangle\Big ]\nonumber\\
&  =\ \exp\!\Big[\big\langle S_{t_{n-1}}\mu,k^{\varkappa}S_{t_{n}-t_{n-1}%
}\varphi_{n}-u_{k}(t_{n}-t_{n-1})\big\rangle\Big ]\,\\
&  \qquad\times\mathbb{E}_{\mu_{k}}\!\exp\!\bigg[\sum_{i=1}^{n-2}k^{\varkappa
}\big\langle X_{t_{i}}^{k}-S_{t_{i}}\mu,\,-\varphi_{i}\big\rangle\nonumber\\
&  \qquad\quad\quad\qquad+\ k^{\varkappa}%
\Big\langle
X_{t_{n-1}}^{k}-S_{t_{n-1}}\mu,\,-\varphi_{n-1}-k^{-\varkappa}u_{k}%
(t_{n}-t_{n-1})%
\Big\rangle
\!\bigg],\nonumber
\end{align}
where $u_{k}$ is the solution of (\ref{new}) with $\varphi$ replaced by
$\,\varphi_{n\,}.$\thinspace\ By Theorem~\ref{T.fluct} for $n=1$, in
$\mathsf{P}$--probability,%
\begin{equation}%
\begin{array}
[c]{l}%
\displaystyle
\liminf_{k\uparrow\infty}\exp\!\Big[\big\langle S_{t_{n-1}}\mu,k^{\varkappa
}S_{t_{n}-t_{n-1}}\varphi_{n}-u_{k}(t_{n}-t_{n-1})\big\rangle\Big ]\\%
\displaystyle
\qquad\geq\ \exp\!\bigg[\underline{c}\,%
\Big\langle
\mu,\int_{t_{n-1}}^{t_{n}}\!dr\ S_{r}(S_{t_{n}-r}\varphi_{n})^{1+\gamma}%
\Big\rangle
\bigg].
\end{array}
\end{equation}
The remaining expectation can be written as
\begin{align}
\mathbb{E}_{\mu_{k}}\!\exp\!\bigg[  &  \sum_{i=1}^{n-2}k^{\varkappa
}\big\langle X_{t_{i}}^{k}-S_{t_{i}}\mu,\,-\varphi_{i}\big\rangle\nonumber\\
&  +\ k^{\varkappa}\big\langle X_{t_{n-1}}^{k}-S_{t_{n-1}}\mu,\,-\varphi
_{n-1}-S_{t_{n}-t_{n-1}}\varphi_{n}\big\rangle\\
&  +\ k^{\varkappa}%
\Big\langle
X_{t_{n-1}}^{k}-S_{t_{n-1}}\mu,\ S_{t_{n}-t_{n-1}}\varphi_{n}-k^{-\varkappa
}u_{k}(t_{n}-t_{n-1})%
\Big\rangle
\bigg].\nonumber
\end{align}
To estimate this term from below observe that, by the induction assumption, in
{$\mathsf{P}$--probability},
\begin{align}
\liminf_{k\uparrow\infty}  &  \,\mathbb{E}_{\mu_{k}}\!\exp\!\Big[\sum
_{i=1}^{n-2}k^{\varkappa}\big\langle X_{t_{i}}^{k}-S_{t_{i}}\mu,\,-\varphi
_{i}\big\rangle\nonumber\\
\qquad &  \qquad\qquad+\ k^{\varkappa}\big\langle X_{t_{n-1}}^{k}-S_{t_{n-1}%
}\mu,\,-\varphi_{n-1}-S_{t_{n}-t_{n-1}}\varphi_{n}\big\rangle\Big ]\nonumber\\
&  \geq\ \exp\!\bigg[\underline{c}\,\Big\langle\mu,\sum_{i=1}^{n-2}%
\int_{t_{i-1}}^{t_{i}}\!dr\ S_{r}\Big(\big(\sum_{j=i}^{n-1}S_{t_{j}-r}%
\varphi_{j}\big)^{1+\gamma}\Big)\\
&  \qquad\qquad+\int_{t_{n-2}}^{t_{n-1}}\!dr\ S_{r}\Big(\big(S_{t_{n-1}%
-r}\varphi_{n-1}+S_{t_{n-1}-r}S_{t_{n}-t_{n-1}}\varphi_{n}\big)^{1+\gamma
}\Big)\Big\rangle\!\bigg].\nonumber
\end{align}
In (\ref{to.show}) it was shown that in {$\mathsf{P}$--probability},
\begin{equation}
\exp\!\bigg[k^{\varkappa}\left\langle _{\!_{\!_{\,_{{}}}}}X_{t_{n-1}}%
^{k}-S_{t_{n-1}}\mu,S_{t_{n}-t_{n-1}}\varphi_{n}-k^{-\varkappa}\,u_{k}%
(t_{n}-t_{n-1})\right\rangle \!\bigg]\ \underset{k\uparrow\infty
}{\Longrightarrow}\ 1.
\end{equation}
As $\,\liminf_{m\uparrow\infty}\xi_{m}\geq a$\thinspace\ in probability, for
some $a\geq0,$ and $\zeta_{m}\Rightarrow1$ in law, implies $\liminf
_{m\uparrow\infty}\xi_{m}\zeta_{m}\geq a$\thinspace\ in probability, this
completes the proof.\hfill$\square$\medskip

\noindent\emph{Acknowledgement.}\thinspace\ The authors are grateful for the
hospitality at the University of Bath repectively the Weierstrass Institute.

{\small
\bibliographystyle{alpha}
\bibliography{bibtex,bibtexmy}
}


\end{document}